\documentclass[11pt]{article}

\usepackage{latexsym,amsmath,amsfonts,amssymb,amstext,
array,enumerate,mathrsfs}

\newtheorem{Theorem}{Theorem}[section]
\newtheorem{Definition}[Theorem]{Definition}
\newtheorem{Proposition}[Theorem]{Proposition}

\newtheorem{Lemma}[Theorem]{Lemma}

\newtheorem{Remark}[Theorem]{Remark}

\numberwithin{equation}{section}

\allowdisplaybreaks

\def\qed{{\hfill\hbox{\enspace${ \square}$}} \smallskip} 
 
\def\dis{\displaystyle }

\def \N{\mathbb{N}}
\def \R{\mathbb{R}}

\def \E{\mathbb{E}}
\def \F{\mathbb{F}}
\def \G{\mathbb{G}}
\def \P{\mathbb{P}}
\def \Q{\mathbb{Q}}

\def \H{\mathbb{H}}

\def \Ac{{\cal A}}
\def \Bc{{\cal B}}

\def \Fc{{\cal F}}

\def \Ic{{\cal I}}
\def \Kc{{\cal K}}

\def \Pc{{\cal P}}

\def \Nc{{\cal N}}
\def \Rc{{\cal R}}
\def \Sc{{\cal S}}
\def \Tc{{\cal T}}

\def \Vc{{\cal V}}

\def \Vc{{\cal V}}

\def\cala{{\cal A}}
\def\calb{{\cal B}}

\def\calf{{\cal F}}
\def\calg{{\cal G}}
\def\calh{{\cal H}}

\def\caln{{\cal N}}

\def\calp{{\cal P}}

\def\calv{{\cal V}}
\def\call{{\cal L}}
\def\cals{{\cal S}}

\def\bfB{{\bf B}}

\def\esssup_#1{\underset{#1}{\mathrm{ess\,sup\, }}}
\def\essinf_#1{\underset{#1}{\mathrm{ess\,inf\, }}}

\begin{document}

\title{The randomization method in stochastic optimal control}

\author{
Marco Fuhrman\thanks{Universit\`a degli Studi di Milano, Dipartimento di Matematica, via Saldini 50, 20133 Milano, Italy;
\newline
e-mail: \texttt{marco.fuhrman@unimi.it} The author is a member of INDAM-GNAMPA.}
}

\maketitle

\begin{abstract}
In this paper we make a survey on the so called randomization method, a recent methodology to study   stochastic optimization problems. It allows to represent the value function of an optimal control problem by a suitable
  backward stochastic differential equation (BSDE), by means of an
auxiliary optimization problem having the same value as the starting one. 
This method works for a large class of control problems and provides a BSDE representation to many related PDEs of  Hamilton-Jacobi-Bellman type, even in the fully non linear case. 

After a general informal introduction 
we explain the method giving full details in a basic case. Then we try to give a complete picture of the existing applications and we present some related open problems.
\end{abstract}

\vspace{5mm}


\vspace{5mm}


\maketitle

\date{}

\section{Introduction}

The aim of this paper is to introduce the reader to the randomization method. This name is given to one possible approach to stochastic optimization problems that has been developed in the last few years. It allows to represent the value of an optimal control problem by a suitably formulated backward  stochastic differential equation (BSDE). This is achieved by introducing an auxiliary optimization problem, called randomized problem, which is equivalent to the original one and admits a BSDE representation. 

So far the randomization method has been applied to a great variety of different situations and new applications are still being developed. We hope that the present exposition will encourage new researchers to contribute to its further development. 

In fact, our main aim is to provide the reader with a complete exposition of the method, including technical aspects and less known prerequisites, so that he/she may be able to access the current literature on the subject and pursue new research autonomously. To this end we have tried to simplify the exposition as much as possible. Although this paper does not contain new results, a large part of the material presented here is not available in the literature since it is the outcome of an effort of synthesis and simplification of many results spread in the literature, corresponding to varying  formulations and often presented without details.

The present exposition also aims at giving a complete account of the present literature  on the randomization method (at least, the part I am aware of) and provide the reader with a picture of the state of the art.

We have tried to present some essential facts about the randomization method to a non-specialized audience, requiring only some knowledge of stochastic optimal control theory and BSDEs. A general description of the method is contained in Section \ref{s;informalintro}, which we tried to write in an informal and hopefully friendly   style. Section
\ref{S:appli}, that describes various applications of the method, is also easily accessible. 
Readers who are not interested in technical details can limit themselves to these two sections.
In the rest of the
paper we expect the audience to have a solid background in stochastic calculus. General references on the prerequisites are given in the next section.  Whenever  specific results are needed, precise references are given to make the exposition as accessible as possible.

The plan of the paper is as follows. In 
Section \ref{s;informalintro} we give a general introduction to the randomization method. 
For the exposition we have chosen  the context of a classical basic optimal control problem, that we assume to be familiar to most readers. At the beginning we quickly recall some of the possible approaches in order to compare the randomization method with them and provide some motivation. In the end of this section we discuss some advantages and disadvantages of this method in comparison with other possible methodologies.

In
Section \ref{S:Pointprocesses} we introduce some notation and recall some prerequisites on marked point processes, that play an important role in the sequel. We also present some results  that have been developed in connection with the randomization method but that may have an interest in themselves and that are anyway  needed later.

In Section  \ref{S:Formulation} we start a rigorous exposition of the method, still in the context of  the  classical control problem discussed in Section \ref{s;informalintro}. We introduce the auxiliary, ``randomized'' optimal control problem and we state the equivalence with the original one. 
Section \ref{proofthm} is devoted to the proof of this equivalence result. Most technicalities occur here and some of them are postponed to the Appendix.

In Section \ref{Sec:separandom} we introduce  rigorously the BSDE that allows to  represent the value of the addressed control problem.  We prove its well-posedness and the representation formula for the value.

In the final  Section \ref{S:appli}   we   give an overview of the  many optimization problems to which the randomization method has been applied and we also indicate some open problems and possible future research directions. The style of this last section will be rather informal, as inclusion of precise assumptions and statements   would take too long.

\section{An informal introduction to the randomization method}
\label{s;informalintro}

The aim of this section is to give an overview of the main steps of the approach to stochastic optimal control problems by means of the randomization method. We will present the main ideas in the familiar context of a basic stochastic optimal control problem: a controlled  stochastic differential equation of It\^o type over a finite horizon, with Markovian coefficients (i.e., deterministic and without memory effects), driven by the Brownian motion, with a classical expected reward functional consisting  of a running   and a final reward. In this Section we will not formulate precise assumptions and results, while starting from Section
\ref{S:Formulation} we will consider the same problem in a fully rigorous manner.

 \subsection{A quick overview of a classical  stochastic optimal control problem}
 
Let $T>0$ be a given deterministic and finite
terminal time. Consider a controlled stochastic equation and a reward functional  of the form
\begin{eqnarray}  \label{dynX1intro}
\left\{
\begin{array}{lll}
  dX_t^\alpha & =&  b(t, X_t^\alpha, \alpha_t)\,dt +
\sigma( t,X_t^\alpha, \alpha_t)\,dW_t, \qquad   t\in [0,T],
\\
 X_0^\alpha&=&x_0, 
\end{array}
\right.
\end{eqnarray} 
\begin{eqnarray} \label{rewardintro}
J(\alpha) & = & \E\Big[\int_0^Tf(t,X_t^\alpha,\alpha_t)\,dt+g(X^\alpha_T)\Big].
\end{eqnarray} 
The controlled process $X^\alpha$ takes values in $\R^n$, 
the initial condition   $x_0$ is a given point in $\R^n$,   while $W$ is a  Wiener process in $\R^{d}$, defined in some probability space $(\Omega,\calf,\P)$. The control process, denoted by $\alpha$, takes values in an abstract space   $A$ of control actions, 
%%%%%%%%%%%%%%%
which is assumed to be a Borel space (see Assumption
{\bf (A1)}
in
Section \ref{SubS:Notation} and the few lines before its statement).
%%%%%%%%%%%%%%%
The class of admissible controls, denoted by $\cala^W$, contains processes  $\alpha$ which are progressive with respect to the filtration generated by $W$. 
The functions $b,\sigma,f$ are defined on $[0,T]\times \R^n\times  A$ with
values in $\R^n$, $\R^{n\times d }$ and $\R$ respectively,  while $g$ is a real function defined on $\R^n$. They are assumed to satisfy appropriate Lipschitz and growth conditions that ensure that, for every choice of $\alpha\in\cala^W$,  the state equation \eqref{dynX1intro} is well posed and the reward functional \eqref{rewardintro} is finite. 
The stochastic optimal control problem  consists in maximizing $J(\alpha)$ over all $\alpha\in\Ac^W$, namely in characterizing the value:
\begin{equation}\label{primalvalueintro}
{\text{\Large$\upsilon$}}_0 = \sup_{\alpha\in\Ac^W} J(\alpha),
\end{equation}
and possibly finding an optimal control, i.e. a control  that attains the supremum.

Loosely speaking, there exist two classical approaches to this problem, that we briefly describe, referring the reader to the treatises \cite{FleSon06}
\cite{YongZhou99}, 
\cite{Pham09} and   \cite{touzibook13} for more information.

The first one is the dynamic programming approach. In a first step it consists in embedding the original control problem into a family of problems parametrized by an initial time $t\in[0,T]$ and an initial condition $x\in\R^n$, namely considering the solution $X_s^{t,x,\alpha}$, $s\in [t,T]$,
to
\begin{eqnarray}  \label{dynX1introdinprog}
\left\{
\begin{array}{lll}
  dX_s^{t,x,\alpha} & =&  b(s, X_s^{t,x,\alpha}, \alpha_s)\,ds +
\sigma( s,X_s^{t,x,\alpha}, \alpha_s)\,dW_s, \qquad   s\in [t,T],
\\
 X_t^{t,x,\alpha}&=&x, 
\end{array}
\right.
\end{eqnarray} 
and defining the value function
\begin{eqnarray} \label{rewardintroprogdin}
v(t,x) & = & \sup_{\alpha\in\Ac^W} \E\Big[\int_t^Tf(s,X_s^{t,x,\alpha},\alpha_s)\,ds+g(X^{t,x,\alpha}_T)\Big].
\end{eqnarray} 
The second step is to prove that the value function is related to the solution to the Hamilton-Jacobi-Bellman equation (HJB): introducing the Kolmogorov operator (depending on a control parameter $a\in A$) acting on a regular function $\phi:\R^n\to\R$
as
\[
\call^a\phi(t,x)= \frac{1}{2}\,{\rm Trace} \left[ D^2\phi(x)\,\sigma\sigma^T(t,x,a)\right] + 
D\phi(x)\cdot b(t,x,a)
\]
where $D\phi$ denotes the gradient and $D^2\phi$ the Hessian matrix, the HJB equation is
\begin{eqnarray}  \label{hjbintro}
\left\{
\begin{array}{rcll}
 -\partial_t v(t,x) & =&  \sup_{a\in A }\left[
\call^av(t,x) + f(t,x,a)
\right], &   t\in [0,T),\,x\in\R^n,
\\
v(T,x)&=&g(x), & x\in\R^n.
\end{array}
\right.
\end{eqnarray} 
Typical results state that, under appropriate assumptions, the value function is the unique solution to HJB in a weak sense, namely in the sense of viscosity solutions. Morevoer, if HJB admits a solution regular enough - and additional conditions are satisfied - one can prove the existence of an optimal control in feedback form.

The second approach is called the Stochastic Maximum Principle. We only sketch it roughly: it consists in associating to the state equation a second one - called dual equation - which is solved backward over the same time interval (in the most general case two dual equations are needed, see \cite{YongZhou99}). Then one can write down necessary and/or sufficient conditions for the optimality of a control process in terms of the solutions to the original and the dual equations.

 \subsection{BSDE representation of the value function}

Starting from the pioneering paper by E. Pardoux and S. Peng \cite{PardouxPeng90}, the theory of Backward Stochastic Differential Equations (BSDEs) started to play a major role in stochastic optimal control theory. We refer for instance to \cite{zhangbook17} or  \cite{pardouxrascanubook14} for  systematic expositions. The dual equations in the stochastic maximum principle are BSDEs of special type (they are linear) and were introduced much earlier 
by J.-M. Bismut 
in the paper
\cite{Bis76}. In many situations BSDEs can also be used to characterize the value function. 
One typical situation is as follows: suppose we are interested in controlled equations of the form
\begin{eqnarray}  \label{dynX1introweak}
\left\{
\begin{array}{lll}
  dX_s^\alpha  & =& \big[b( s,X_s^\alpha)+ \sigma( s,X_s^\alpha)\,r(s, X_s^\alpha, \alpha_s)\big]\,ds +
\sigma( s,X_s^\alpha)\,dW_s, \qquad   s\in [t,T],
\\
 X_t^\alpha&=&x, 
\end{array}
\right.
\end{eqnarray} 
namely where the volatility $\sigma$ is not controlled and the drift has the special form written above, for a suitable function $r: [0,T]\times\R^n\times A\to \R^d$.
One can formulate the controlled problem in a weak form: instead of solving directly \eqref{dynX1introweak} one starts from the (uncontrolled) equation
\[
 dX_s    =   b( s,X_s)\,ds +
\sigma( s,X_s)\,dW_s, \qquad   s\in [t,T],
\]
on a probability space $(\Omega,\calf,\P)$. Given a control process $\alpha$, one makes a change of probability of Girsanov type introducing the martingale
\[
L_s=\exp\left( \int_t^s 
r(u, X_u, \alpha_u)\,dW_u-\frac{1}{2} \int_t^s 
|r(u, X_u, \alpha_u)|^2\,du\right),\quad s\in [t,T],
\]
and setting $\P^\alpha(d\omega)=L_T(\omega)\,\P(d\omega)$. By the Girsanov theorem, the process
\[
W^\alpha_s=W_s-\int_t^s r(u, X_u, \alpha_u)\,du, \quad s\in [t,T],
\]
is a Brownian motion under $\P^\alpha$ and the process $X$ satisfies
\begin{eqnarray}  \label{dynX1introweakdue}
\left\{
\begin{array}{lll}
  dX_s   & =& \big[b( s,X_s )+ \sigma( s,X_s )\,r(s, X_s , \alpha_s)\big]\,ds +
\sigma( s,X_s)\,dW^\alpha_s, \qquad   s\in [t,T],
\\
 X_t&=&x.
\end{array}
\right.
\end{eqnarray} 
Therefore, $X$ is a weak solution to \eqref{dynX1introweak} in $(\Omega,\calf,\P^\alpha)$ and we formulate an optimal control problem writing the value function as
\begin{eqnarray} %\label{rewardintroprogdin}
v(t,x) & = & \sup_{\alpha\in\Ac^W} \E^\alpha\Big[\int_t^Tf(s,X_s,\alpha_s)\,ds+g(X_T)\Big],
\end{eqnarray} 
where $\E^\alpha$ denotes expectation under $\P^\alpha$.

In this case one can introduce a BSDE to represent the value function. One defines  the Hamiltonian
\begin{align}\label{hhamiltbsdecl}
    h(t,x,z)=\sup_{a\in A} \big[r(t,x,a)\cdot z + f(t,x,a)\big],
\end{align}
for $t\in [0,T]$, $x\in\R^n$, $z\in\R^d$.  
%%%%%%%%%%%%%%%%%%%
Under appropriate conditions (for instance, when $r$ is bounded and $f$ has linear growth in $x$ uniformly in $(t,a)$) the function $h$ in \eqref{hhamiltbsdecl} is Lipschitz in $z$ uniformly in $(t,x)$ and has linear growth in $(x,z)$ uniformly in $t$, so one can solve the BSDE
%%%%%%%%%%%%%%%
\begin{eqnarray}   
\left\{
\begin{array}{lll}
  -dY_s   & =& -Z_s\,dW_s +h(s,X_s,Z_s)\,ds,  \qquad   s\in [t,T],
\\
Y_T&=&g(X_T),
\end{array}
\right.
\end{eqnarray} 
for the pair of unknown processes $(Y,Z)$ with values in $\R\times \R^d$. Then it can be proved that $Y_t$ is almost surely constant (i.e., deterministic), that
\begin{equation}\label{FeynmanKac}
    Y_t=v(t,x)
\end{equation}
and that the control $\hat a(s,X_s,Z_s)$, $s\in [t,T]$, is optimal, if a suitable maximizer $\hat a(t,x,z)$ in the Hamiltonian    exists.

An easy computation shows  that the HJB equation in the present case has the form
\begin{eqnarray}  \label{hjbintrosemilinear}
\left\{
\begin{array}{rcll}
 -\partial_t v(t,x) & =&  
\call v(t,x) +h\big(t,x,\sigma(t,x)^T Dv(t,x)\big), &   t\in [0,T),\,x\in\R^n,
\\
v(T,x)&=&g(x), & x\in\R^n,
\end{array}
\right.
\end{eqnarray} 
where, for smooth $\phi:\R^n\to\R$,
\[
\call \phi(x)= \frac{1}{2}\,{\rm Trace}\left[ D^2\phi(x)\,\sigma\sigma^T(t,x)\right] + 
D\phi(x)\cdot b(t,x).
\]
The formula 
\eqref{FeynmanKac} is called non-linear Feynman-Kac formula. Using BSDEs it is therefore possible to give a probabilistic representation (or in some case, to construct solutions) for non-linear PDEs of parabolic type of the form \eqref{hjbintrosemilinear}.
This connection with PDEs, first explicitly described in \cite{PardouxPeng92}, has led to many results and has been one of the motivations for the great development of the theory of BSDEs (see e.g. \cite{pardouxrascanubook14}).

In spite of that, this method based on BSDEs has some limitations: on the analytic side, in contrast to the fully non linear HJB equation \eqref{hjbintro}, the PDE \eqref{hjbintrosemilinear} (as well as similar variants) is semilinear, that is, it is linear with respect to higher order derivatives; moreover, the nonlinearity with respect to the gradient $Dv$ has a specific form (due to the occurrence of $\sigma$ in the nonlinearity). 
On the control-theoretic side, the addressed optimization problem \eqref{dynX1introweak} has a rather particular form. 

In order to remove these restrictions, many efforts were done during the last decades. In particular, two entire theories were developed.

The first one is the theory of Second-Order BSDEs,   developed mainly by H.M. Soner and N. Touzi.
They are a special class of BSDEs that allow to represent solutions to fully non linear PDEs. They are formulated looking for unknown processes that solve the BSDE under a family of mutually singular probabilities.
Relevant papers are \cite{ChSoToVi2007}, \cite{STZ} and we also refer the reader  to the book \cite{zhangbook17}.

The second one is the theory of $G$-expectation, invented by S. Peng and developed by him and his collaborators. It consists in introducing new basic definitions in stochastic analysis starting from a concept of non-linear expectation and developing variants of the main notions and results of probability theory and stochastic calculus. The theory is particularly suitable for describing uncertainty and for handling nonlinear models. In particular, it provides probabilistic-type representations for solutions of fully non linear PDEs. Some important papers are
\cite{Peng07}, \cite{Peng08}, \cite{HuJiPeng2014} and see also the    book \cite{zhangbook17}.

More recently, a new method has been devised, that is able to extend the Feynman-Kac formula to fully non linear cases and to give a BSDE representation for the the value function of general optimal control problem. It has been called the randomization method. It was introduced in \cite{bou09} by B. Bouchard in the context of optimal switching problems, already in a very precise form. Among other early basic contributions we must cite \cite{KMPZ10}, 
\cite{EK2010}  and \cite{KP12}. 
By now, it has become a flexible tool to cope with a variety of optimization problems for stochastic processes.
In the next subsection we are going to present its application to the classical control problem introduced above.

\subsection{The randomization method}

This method consists of two steps. In the first step a new, auxiliary control problem is introduced, which is equivalent to the original one in the sense that the value (or the value function) is the same. In the second step, a suitable BSDE is introduced and its solution is shown to represent the value function.
Full details of the two steps will be given in Sections 
\ref{S:Formulation} and \ref{Sec:separandom}, while here we proceed informally as before.

\subsubsection{The randomized control problem}\label{subsub;randomized}

Let us come back to the control problem 
\eqref{dynX1intro}-\eqref{rewardintro}
and its value \eqref{primalvalueintro}.
In order to formulate the auxiliary control problem we first 
formally replace the control process $\alpha$ with another (uncontrolled) process $I$ with values in the control actions space $A$, independent of the Brownian motion $W$ and defined on the same probability space. Then we consider the corresponding state equation
\begin{eqnarray}  \label{dynXrandomintro}
\left\{
\begin{array}{lll}
  d  X_t &=&  b(t,   X_t,  I_t)\,dt + \sigma(t,  X_t,  I_t)\,d  W_t, \qquad   t\in [0,T],
\\
 X_0 &=&x_0, 
\end{array}
\right.
\end{eqnarray} 
Thus, we say that the control is randomized and we call equation \eqref{dynXrandomintro} the randomized state equation;  more generally this explains the terminology ``randomization method''.
As the process $I$ it is convenient to choose a piecewise-constant process associated with a Poisson random  measure. Indeed, this can be done even for a very general control action space $A$, although in specific situations different choices are possible: see Remark \ref{randomconb}. Thus, we define
\[
%\begin{equation}\label{I}
  I_t \ = \ \sum_{n\ge 0}  \eta_n\,1_{[  S_n,  S_{n+1})}(t), \qquad t\ge 0,
%\end{equation}
\]
where $(\eta_n)$ is an iid sequence of $A$-valued random variables, and  $(S_n)$ is an independent Poisson process on the real line ($S_0=0$, and $(S_n-S_{n-1})_{n\ge1}$ is an iid sequence of exponentially distributed random variables).

Next we define an optimization problem introducing an index set $\calv$ - to be described later - and, for every $\nu\in\calv$, a probability $\P^\nu$ equivalent to $\P$ and such that $W$ remains a Brownian motion under each $\P^\nu$. 
We finally introduce the reward functional 
\begin{eqnarray}  \label{defJrandomizedintro}
J^\Rc( \nu) &=&   \E^{ \nu}
\left[\int_0^Tf(t,  X_t,  I_t)\,dt+g(  X_T)\right],
\end{eqnarray} 
where $\E^{ \nu}$ denotes the expectation under $\P^{ \nu}$,
which we try to  maximize
over all $ \nu\in \calv$. The value of this ``randomized'' control problem
is defined as
\begin{equation}\label{dualvalueintro}
{\text{\Large$\upsilon$}}_0^\Rc \;=\;   \sup_{  \nu\in  \calv} J^\Rc(  \nu).
\end{equation}

The idea behind this formulation is that by choosing a probability $\P^\nu$ we control  the law of the process $I$, without affecting the law of the Brownian motion. In contrast to the original control problem, where we choose a process $\alpha$ and construct the corresponding trajectory $X^\alpha$, here we fix the pair $(I,X)$ solution to \eqref{dynXrandomintro} but we modify its law by choosing $\P^\nu$ and obtaining the corresponding reward. We will see that in great generality the values coincide:
\[
{\text{\Large$\upsilon$}}_0={\text{\Large$\upsilon$}}_0^\Rc.
\]
In order for this result to be true it is required that the random variables $\eta_n$ should have a distribution with full topological support in $A$: this allows the process $I$ to ``visit'' the entire control space and achieve the same value as the control processes $\alpha$ in the original optimization problem.

The required probabilities $\P^\nu$ are constructed in a specific way, namely by a Girsanov change of measure that we now describe. $\calv$ is chosen as the collection of all the random fields 
\[\nu=\nu_t(\omega,a):\Omega\times [0,\infty)\times A\to (0,\infty)
\]
that are bounded and predictable (with respect to the filtration generated by $W$ and the Poisson random measure). Given such $\nu$, one defines the exponential martingale
\[  
\kappa_t^{ \nu} \ 
= \ \exp\left(\int_0^t\int_A (1 -   \nu_s(a))\lambda(da)\,ds
\right)\prod_{0<  S_n\le t}\nu_{  S_n}(  \eta_n),\qquad t\ge 0,
\]
and defines the required probability on $( \Omega, \calf)$ setting
$ \P^{ \nu}(d \omega)=\kappa_T^{ \nu}( \omega)\, \P(d \omega)$. The effect of changing the law of $I$ can be described as follows. 
Denote by $\lambda$ the rate of the exponential variables $  S_n-  S_{n-1} $ and by $\pi(da)$ the law of the variables $\eta_n$; define the (finite) measure $\lambda(da)=\lambda\,\pi(da)$ on $A$. Denote
 $ \mu=\sum_{n\ge 1}\delta_{(  S_n, \eta_n)}$ the Poisson random measure on $(0,\infty)\times A$. Then under $\P$ the measure
 \[
 \mu(dt\,da) -\lambda(da)\,dt
 \]
 is a martingale measure on $(0,T]\times A$, while 
under $\P^\nu$ the measure
 \[
 \mu(dt\,da) -\nu_t(a)\,\lambda(da)\,dt
 \]
 is a martingale measure. With a different terminology, passing from $\P$ to $\P^\nu$ changes the compensator of the random measure from $\lambda(da)\,dt$ to $\nu_t(a)\,\lambda(da)\,dt$.
  By known properties of point processes, the compensator entirely characterizes the law of the random measure and hence of $I$ (while the law of $W$ remains that of a Brownian motion). We may also say that in the randomized problem the  controls are random intensities $\nu$ and that $\calv$ is the set of admissible controls.

\subsubsection{The associated BSDE}

   The equation that we are going to introduce 
will be called the \emph{randomized BSDE}. It is  the following constrained BSDE
on the time interval $[0,T]$:
\begin{equation}\label{BSDEconstrainedintro}
\begin{cases}
\vspace{2mm} \dis Y_t \ = \ g(X_T)  + \int_t^T  f(s,X_s ,I_s) ds + K_T - K_t - \int_t^TZ_s\,dW_s - \int_t^T\!\int_A U_s(a)\,\mu(ds\,da), \\
\dis U_t(a) \ \le \ 0.
\end{cases}
\end{equation}
Here $(X,I)$ the same as in \eqref{dynXrandomintro} and $\mu$ is the Poisson random measure introduced before. As the underlying  filtration we take the one generated by $W$ and  $\mu$. The unknown is a quadruple
$(Y_t,Z_t,U_t(a),K_t)$, $(t\in [0,T]$, $a\in A)$, satisfying appropriate regularity and integrability conditions; in particular, $Y$ and $Z$ are progressive, the process $K$ and the random field $U$ are predictable, $K$ is non-decreasing and $K_0=0$. Since $U$ occurs in the integral with respect to $\mu$, the sign constraint $U_t(a)\le 0$ implies a sign condition on the corresponding jumps: for this reason many authors call this equation a BSDE with constrained jumps. 
For short, we will also call it the constrained BSDE.

In the absence of the constraint $U_t(a)\le 0$ it would be immediate to find a solution to the simpler equation:
\[
\dis Y_t  =  g(X_T)  + \int_t^T  f(s,X_s ,I_s) ds   - \int_t^TZ_s\,dW_s - \int_t^T\int_A U_s(a)\,\mu(ds\,da),\quad t\in [0,T].
\]
Indeed, this equation is linear in the unknown triple $(Y,Z,U)$ and it can be solved by an application of a representation theorem for martingales (with respect to the filtration generated by a Wiener and a Poisson process): the reader may look at the proof of Lemma \ref{lemmapenalized} below, or consult directly  \cite{LiTang94}  Lemma 2.4 for full details. So the equation without constraint has a solution with $K\equiv 0$. The occurrence of the increasing process $K$ is required to deal with the constraint.

We note that the solution to equation 
\eqref{BSDEconstrainedintro} is clearly non unique. Indeed, given a solution $(Y_t,Z_t,U_t(a),K_t)$, another solution is $(Y_t+h(T)-h(t),Z_t,U_t(a),K_t+h(t))$ where $h$ is an arbitrary  non-decreasing deterministic function with $h(0)=0$. 
However, it is possible to prove that the BSDE \eqref{BSDEconstrainedintro} has a unique minimal solution: minimality of a solution $(Y_t,Z_t,U_t(a),K_t)$ means that if $(Y_t',Z_t',U_t'(a),K'_t)$ is another solution then $Y_t\le Y'_t$ for every $t\in[0,T]$. With this definition it is obvious that two minimal solutions have the same first component $Y$, but it can be proved that the others coincide as well, so that uniqueness of the minimal solution holds true.

Existence of a minimal solution  can be obtained by a penalization method: a sequence of uniquely solvable equations is defined depending on a parameter $n$ that later tends to infinity. These equations are constructed in such a way that the first component $Y^n$ of their solutions represents the value of a penalized control problem: in fact, this is exactly the randomized optimal control problem introduced in the previous paragraph \ref{subsub;randomized}, but with the additional condition that the intensity $\nu$ take values in $(0,n]$. From this representation
it follows easily that $Y^n$ are monotonically increasing in $n$ and they are  bounded, so they converge to a limit process $Y$. The convergence of the other components is obtained by direct estimates. 
Full details of the formulation and the proofs can be found in Section 
\ref{Sec:separandom}.

From the representation formulae for the penalized problem one concludes that the minimal solution 
$(Y_t,Z_t,U_t(a),K_t)$ represents the value of the randomized problem:
$Y_0= {\text{\Large$\upsilon$}}_0^\Rc$. From the equality of the values proved above we finally conclude that 
\[Y_0=
{\text{\Large$\upsilon$}}_0.
\]
This is the generalization of the nonlinear Feynman-Kac formula \eqref{FeynmanKac}
that we were looking for.
By considering a control problem starting at time $t$ and a similar BSDE on the time interval $[t,T]$ one obtains a representation formula for the value function $v(t,x)$ in the same way: see Subsection
\ref{subs:randomdynprogr} for more details.
In the end we may conclude that the randomization method provides a BSDE representation for general stochastic optimal control problems and fully non linear Hamilton-Jacobi-Bellman equations.

\begin{Remark}
    \emph{ 
BSDEs of the form \eqref{BSDEconstrainedintro} show similarities with reflected BSDEs: in the reflected case, the component $Y$ is constrained to remain above a given random obstacle and it is necessary to introduce an additional non-decreasing unknown process $K$ to obtain well-posedness. In addition, a minimality condition (the Skorohod condition)  also occurs in the formulation. Finally,  in both cases the solution can be obtained by a penalization method. However, in the present case the constraint is on one 
of the martingale terms and not on the first component. Moreover, the minimality condition is of different kind and no analogue of the Skorohod condition is known in general. \qed
}
\end{Remark}

\begin{Remark}
    \emph{ 
The idea of randomizing controls is not new. In Statistics randomized tests are used with many purposes, for instance to achieve an arbitrary prescribed size for a statistical test
even in presence of discrete distributions. In stochastic optimal control problems randomized strategies have been already employed as well.
Sometimes they have been introduced in a form related to the specific optimization problem, for instance  in the case  of optimal stopping problems. But the more common randomization technique is probably  the concept of 
relaxed controls. When relaxed controls are adopted, at each time the agent does not select a single point in the action space $A$, but rather he/she chooses a probability distribution on $A$, i.e. an element of $\calp(A)$. One assumes that the actual action is selected randomly according to the chosen distribution and independently of the other sources of noise. The dynamics of the controlled system as well as the functional to be optimized must be rewritten accordingly. Using relaxed controls amounts to replacing $A$ with $\calp(A)$ and formulating a new control problem of standard form. The advantage is that  the new action space $\calp(A)$ enjoys better properties; typically the convexity of  $\calp(A)$ allows to prove existence of an optimal control using abstract techniques of optimization. 
}

\emph{ 
From the previous description of the randomization method it should be clear that relaxed controls are a very different technique.
\qed
}
\end{Remark}

\begin{Remark}\label{randomconb}
\emph{ In our exposition of the randomization method we have introduced an auxiliary randomized state process, solution to equation 
\eqref{dynXrandomintro},
where  the control was replaced by a piecewise constant process $I$. This is not the only possibility.  Occasionally, for specific purposes, the process $I$ is taken as an independent Brownian motion (when the control space $A$ is  an Euclidean space) or the image of a Brownian motion under some appropriate map. The associated constrained BSDE needs a corresponding reformulation. This is the point of view, for instance, in \cite{ChoukrounCosso16} or \cite{CoGuaTe2019}, to which we refer the reader for a detailed exposition.
    \qed
}
\end{Remark}

\subsection{Discussion}

Having described the randomization method in general terms we wish to discuss its usefulness in stochastic optimal control and its advantages and disadvantages. The considerations that follow may reflect my personal taste and may be affected by my limited knowledge of some subject, but I will try to justify them with convincing arguments.

One of the good features of the randomization method is its great flexibility. As we will see in the  Section \ref{S:appli}, it can be successfully applied to a large variety of stochastic optimization problems: classical control (with finite or infinite horizon, or even ergodic control), switching control, impulse control, control of point processes, control of piecewise-deterministic Markov processes. It has been applied to systems evolving in infinite-dimensional spaces, or systems showing mean-field effects (control of McKean-Vlasov type). It works equally well for systems that exhibit memory effects, hence in a non-Markovian context, with little additional effort.

This method has contributed to widen applicability of the theory of BSDEs: it shows that rather classical BSDEs (driven by Wiener and Poisson processes, containing a constraint  and understood in the sense of minimal solutions) can represent the value of general optimal control problems and give a probabilistic representation formula for fully non linear PDEs, at least of second order elliptic or parabolic type and for scalar unknown functions. For instance, a BSDE representation has become possible for the value function of a partially observed control problem in the classical general case, a problem that had remained open for a long time: see Subsection \ref{subs:partialobs}.
Since the numerical analysis of BSDEs is also rapidly developing, this opens the possibility of finding new efficient numerical methods in optimal control.

In comparison with the theory of Second Order BSDEs, the randomization method is based on more classical BSDEs and completely avoids the requirements on nondegeneracy of the noise that sometimes occurs in that theory.
In comparison with the theory of $G$-expectation, this method does not require to construct new foundations for probability and stochastic analysis.

On the other hand, the randomization method still has some disadvantages, so that one cannot immediately claim its superiority in comparison with 
Second Order BSDEs or $G$-expectation.

The first limitation is that the randomization method is not very effective in characterizing an optimal control or even to prove its existence. It is designed mainly to represent the value of the optimization problem by means of a BSDE. This disadvantage is common to other successful theories that have been applied to optimal control problems, for instance the theory of viscosity solutions to the Hamilton-Jacobi-Bellman equation. However, in that context there exist some results in this direction, for instance some generalized variants of the verification theorem, that are still missing for the randomization method.

A second inconvenience is the fact that the BSDE that represents the value is rather difficult to handle. Not only does it include an additional source of noise that was not present in the original control problem (the Poisson random measure), but it is also difficult to obtain direct estimates on its solution, the reason being essentially that the solution is not unique (only the minimal solution is uniquely determined). As a consequence, it is difficult to use the BSDE to get additional information on the value function of the original control problem, for example its continuous or regular dependence on parameters. 

One can hope that in the future the randomization method will be further developed. Some  more open problems will be outlined in the last Section.

\section{Preliminaries on marked point processes}
\label{S:Pointprocesses}

In the previous Section the randomization method was presented and applied to a classical optimal control problem without any precise assumption or result. It is our purpose to give a precise formulation of the same control problem,  including all the required assumptions and to present with full details how the randomization method can be used for its study. As we saw, a basic role is played by the introduction of an appropriate point process (a Poisson random measure) and some advanced notions in the theory of point processes are needed. 
The purpose of this Section is to collect all the results we need on the theory of marked point processes and that play a role in the randomization method. For the general theory   we refer the reader to the books
\cite{brandtlast} or \cite{bremaud2} (but see also the more pedagogical introduction \cite{bremaud}), or to the paper \cite{ja}. In fact the latter paper contains all the results we need, even in greater generality.

In paragraph \ref{subs;pointreminders} we only review known results and the expert reader may skip this part and consult it later for  some notation. 
In paragraph
\ref{subs;timechanges} we present some specific  results in point process theory that will be applied in the Sections that follow.

\subsection{Reminders on marked point processes}\label{subs;pointreminders}

  Let  $A$ be a Borel  space, namely   a topological space  homeomorphic to a  Borel subset of a Polish space 
(some authors use the terminology  Lusin space).
We recall the definition of marked point process with values in $A$ (sometimes called multivariate point process). We take a point $\Delta\notin A$ and we add it to $A$ as an isolated point obtaining the space $A\cup \{\Delta\}$.

Let
$(\Omega,\calf)$ be a measurable space with a filtration $\F=(\calf_t)_{t\ge 0}$. A sequence $(T_n,A_n)_{n\ge 1}$ defined on $\Omega$ is called a marked point process if:
\begin{itemize}
\item[\textup{(i)}] for every $n$, $(T_n,A_n)$ takes values in $\{(0,\infty)\times A\}\cup \{(\infty,\Delta)\}$;
\item[\textup{(ii)}] for every $n$  and $\omega$
we have $T_n(\omega) \le T_{n+1}(\omega)$ and
if $T_n(\omega)<\infty$ then $T_n(\omega) <T_{n+1}(\omega)$;
\item[\textup{(iii)}] for every $n$, $T_n $ is an $\F $-stopping time and $A_n $ is $\Fc_{T_n } $-measurable.

\end{itemize}
In the sequel we will only consider processes which are non-explosive, meaning that for every   $\omega$ we have $\lim_{n\rightarrow\infty}T_n(\omega)=\infty$.
Very often we will consider processes satisfying $(T_n(\omega),A_n(\omega))\in (0,\infty)\times A$, or equivalently $T_n(\omega)<\infty$, but this is not assumed everywhere.

To the process $(T_n,A_n)$ we may associate a random measure $\mu$, defined on the Borel subsets of $(0,\infty)\times A$ given by
\[
\mu=\sum_{n\ge1}\delta_{(T_n,A_n)}\,1_{\{T_n<\infty\}},
\]
where $\delta$ is the Dirac measure. 
The counting processes are defined, for every $C\in\calb(A)$, by the formula
\[
N^C_t=\mu((0,t]\times C), \qquad t\ge 0.
\]
The filtration generated by the counting processes is called the natural filtration $\F^\mu=(\calf_t^\mu)_{t\ge0}$ where
\[
\calf_t^\mu=\sigma\Big( \mu((0,s]\times C)\,:\, s\in [0,t],C\in\calb(A)\Big).
\]

The main properties of the natural filtration are summarized in the following proposition.

\begin{Proposition} \label{propnaturalfiltration}

\begin{itemize}
\item[\textup{(i)}] 
  $\F^\mu$ is right-continuous.
  
\item[\textup{(ii)}] 
$\calf_{T_n}^\mu=\sigma(T_1,A_1,\ldots,T_n,A_n)$.
  
\item[\textup{(iii)}] 
Denote by 
$\Pc(\F^\mu) $ the predictable $\sigma$-algebra of $\F^\mu$.
A process $\nu:\Omega\times [0,\infty)\to \R$ is $\F^\mu$-predictable (equivalently, $\Pc(\F^\mu) $-measurable) if and only if it 
 admits the following representation
\begin{align*}
\nu_t(\omega) \ &= \ \nu_t^{(0)} \, 1_{\{0\le t\leq T_1(\omega)\}} \\
&\quad \ + \sum_{n=1}^\infty \nu_t^{(n)}\big( T_1(\omega ), \ldots,T_n(\omega) \big) \, 1_{\{T_n(\omega)<t\leq T_{n+1}(\omega)\}}, \notag
\end{align*}
for all $(\omega, t)\in \Omega\times\R_+$, for some (deterministic) maps $\nu^{(n)}\colon \R_+\times(\R_+ )^n\rightarrow\R$, $n\geq1$, and  $\nu^{(0)}\colon \R_+\rightarrow\R$, which are Borel-measurable.  
\item[\textup{(iv)}] 
More generally, a random field
$\nu:\Omega\times [0,\infty)\times A\to \R$ is $\Pc(\F^\mu)\otimes \calb(A) $-measurable if and only if it has the form
\begin{align*}
\nu_t(\omega, a) \ &= \ \nu_t^{(0)} ( a ) \, 1_{\{0\le t\leq T_1(\omega)\}} \\
&\quad \ + \sum_{n=1}^\infty \nu_t^{(n)}\big( T_1(\omega ),A_1(\omega),\ldots,T_n(\omega),A_n(\omega),a\big) \, 1_{\{T_n(\omega)<t\leq T_{n+1}(\omega)\}}, \notag
\end{align*}
for all $(\omega, t,a)\in \Omega\times\R_+\times A$, for some   maps $\nu^{(n)}\colon \R_+\times(\R_+\times A)^n\times A\rightarrow\R$, $n\geq1$, $($resp. $\nu^{(0)}\colon \R_+\times A\rightarrow\R$$)$, which are $\Bc(\R_+\times (\R_+\times A)^n \ \times A)$-measurable $($resp. $\Bc(\R_+\times A)$-measurable$)$.

\end{itemize}

\end{Proposition}

 \noindent {\bf Proof.} For points $(i)-(ii)$ we refer to
\cite{bremaud2} Theorem A.4.1.
Point $ (iii)$ is in 
\cite{ja}   Lemma 3.3 or 
\cite{bremaud2} Theorem A.4.4.
Point $ (iv)$ follows from 
point $ (iii)$ and a monotone class argument.
\qed

So far no probability was considered. To every probability $\P$ on $(\Omega,\calf)$ there corresponds a random measure $\mu^\P$ on the Borel subsets of $(0,\infty)\times A$, called the compensator (or dual predictable projection) of $\mu$,   whose  properties are stated in the following proposition.

\begin{Proposition} \label{propcompensator}
Let $(T_n,A_n)_{n\ge 1}$ be a marked point process defined in $(\Omega,\calf,\F)$ and let   $\P$ be a probability in
$(\Omega,\calf)$.  Then there exists a transition kernel $\mu^\P(\omega,dt\,da)$ from $(\Omega,\calf)$ to  $((0,\infty)\times A,\calb((0,\infty)\times A))$ with the following properties.

\begin{itemize}
\item[\textup{(i)}] $\mu^\P$ is $\F$-predictable: this means that   the counting processes
$(\mu^\P((0,t]\times C)_{t\ge0}$ are $\F$-predictable, for every choice of $C\in\calb(A)$.
\item[\textup{(ii)}] for every $n$  and $C\in\calb(A)$ the processes
\[
\mu((0,t\wedge T_n]\times C)-
\mu^\P((0,t\wedge T_n]\times C), \qquad t\ge0,
\]
are $\F$-martingales under $\P$; since we are assuming $T_n\to\infty$ the processes
\[
\mu((0,t]\times C)-
\mu^\P((0,t]\times C), \qquad t\ge0,
\]
are local martingales;

\item[\textup{(iii)}] $\mu^\P(\{t\}\times A)\leq 1$ for every $t\ge0$;
\item[\textup{(iv)}] The measure $\mu^\P$, satisfying \emph{(i)-(ii)-(iii)} above, is unique, up to a $\P$-null set.
\item[\textup{(v)}] If two probabilities $\P$ and $\Q$  defined on $(\Omega,\calf)$ admit the same compensator then $\P=\Q$ on the $\sigma$-algebra $\calf_\infty^\mu$. In particular,  the law of $(T_n,A_n)_{n\ge 1}$ is the same under $\P$ and $\Q$ (in other words, the compensator uniquely specifies the law of $(T_n,A_n)_{n\ge 1}$).

\item[\textup{(vi)}] Suppose   that a random field
$H:\Omega\times [0,\infty)\times A\to \R$ is $\Pc(\F^\mu)\otimes \calb(A) $-measurable. If 
\[
\int_{(0,t]}\int_{ A}|H_s(a)|\, \mu^\P(ds\,da)<\infty, \qquad \P-a.s.
\]
for every $t\ge0$,    then the process
\[
M_t:=\int_{(0,t]}\int_{ A}H_s(a)\, \mu(ds\,da)-
\int_{(0,t]}\int_{ A}H_s(a)\, \mu^\P(ds\,da), \qquad t\ge0,
\]
is a local martingale; if we even have
\[
\E \left[\int_{(0,t]}\int_{ A}|H_s(a)|\, \mu^\P(ds\,da)\right]<\infty 
\]
for every $t\ge0$, 
then $M$ is a  martingale.
\item[\textup{(vii)}] If  a random field
$H:\Omega\times [0,\infty)\times A\to [0,\infty)$ is $\Pc(\F^\mu)\otimes \calb(A) $-measurable then
\begin{equation}\label{smoothing}
\E\int_{0}^\infty\int_{ A}H_t(a)\, \mu(dt\,da)=
\E\int_{0}^\infty\int_{ A}H_t(a)\, \mu^\P(dt\,da).    
\end{equation}
Conversely, if a kernel $\bar\mu(\omega, dt\,da)$ is $\F$-predictable (i.e., the associated counting processes are predictable) and \eqref{smoothing} holds for $\bar\mu$ and for every such $H$ then $\bar\mu=\mu^\P$ up to a $\P$-null set.
\end{itemize}
\end{Proposition}

 \noindent {\bf Proof.} For points $(i)-(iv)$ we refer to
\cite{ja}   Theorem 2.1, Proposition 2.3 and 2.4,   or to  \cite{bremaud2} Theorem 5.6.4.
Point $ (v)$ is in 
\cite{ja}   Theorem 3.4.
Point $ (vi)$ can be found in \cite{bremaud2} Theorem 5.1.31 and
Corollary 5.1.32. Point  $(vii)$ is \cite{ja}   Theorem 2.1, or \cite{bremaud2} Theorem 5.1.20 where it is called the ``smoothing formula".
\qed

An important example is the Poisson process, that will be used extensively below: given a finite positive measure $\lambda(da)$ on the Borel subsets of $A$, we say that $(T_n,A_n)_{n\ge1}$, defined on a probability space $(\Omega,\calf,\P)$, is a Poisson process with (finite) intensity $\lambda(da)$ if $(T_n-T_{n-1})_{n\ge1}$ are independent and have common exponential law with parameter $\lambda(A)$ (we set $T_0=0)$, $(A_n)_{n\ge1}$ are independent and have common law $\lambda(da)/\lambda(A)$, the two sequences $(T_n)$ and $(A_n)$ are independent. In this case the compensator with respect to the natural filtration $\F^\mu$ is given by $\mu^\P(dt\,da)=\lambda(da)\,dt$ (a deterministic product measure).

\begin{Remark}\label{largfiltrcomp}\emph{
If $\G=(\calg_t)_{t\ge0}$ is another filtration in $(\Omega,\calf)$ and $\calg_\infty$ and $\calf^\mu_\infty$ are independent under $\P$, it is easy to see that $\mu^\P(dt\,da)=\lambda(da)\,dt$ is also the compensator with respect to the larger filtration $(\calf^\mu_t\vee \calg_t)_{t\ge0}$
    }
\end{Remark}

A converse result also holds, which is known as Watanabe's theorem (see \cite{bremaud2} Theorem 5.7.1):
\begin{Theorem} \label{watanabe}
Let $(T_n,A_n)_{n\ge 1}$ be a marked point process defined in $(\Omega,\calf,\F,\P)$ with compensator $\mu^\P(dt\,da)=\lambda(da)\,dt$, where $\lambda(da)$ is a finite Borel measure on $A$. Then  
$(T_n,A_n)_{n\ge 1}$ is a Poisson process with intensity $\lambda(da)$.
\end{Theorem}

In the sequel we will make use of the canonical   space $(\Omega,\calf)$ of a marked point process; we limit ourselves to the case of a 
 non-explosive process   with finite arrival times:  
  $\Omega$ is the set of sequences $\omega=(t_n,a_n)_{n\geq1}\subset(0,\infty)\times A$ with $t_n<t_{n+1}\to\infty$; 
   the canonical   process 
    $(T_n,A_n)_{n\geq1}$ is defined setting 
   $T_n(\omega)=t_n$, $A_n(\omega)=a_n$;
   the $\sigma$-algebra $\Fc$ is the   smallest $\sigma$-algebra such that all the maps $T_n,A_n$ are measurable.
   If $\F^\mu=(\calf_t^\mu)_{t\ge0}$ denotes the natural filtration then we have $\calf=\calf_\infty^\mu$. We also note that the  $\sigma$-algebras  $\calf^\mu_{T_n}$ and $\calf_\infty^\mu$ are isomorphic, respectively,  to 
   \begin{equation}
       \label{infproduct}
   \Big(
 \calb((0,\infty)\times A))
\Big)^n, \qquad \Big(
 \calb((0,\infty)\times A))
\Big)^\N   , 
   \end{equation} 
 where the latter is  the countably infinite  product of   the indicated Borel $\sigma$-algebra.

The canonical space is often used to obtain existence results. For instance, given an arbitrary finite measure $\lambda(da)$ on $\calb(A)$, it is well known that there exists a unique probability on the canonical   space $(\Omega,\calf)$ such that the canonical process is a Poisson process with   intensity $\lambda(da)$.
   \bigskip

The following result is a version of  the
Girsanov theorem for multivariate point processes.

\begin{Theorem} \label{girsanovpoint}
Let $(T_n,A_n)_{n\ge 1}$ be a marked point process defined in $(\Omega,\calf,\F,\P)$ with compensator $\mu^\P$ satisfying
\[
\mu^\P(\{t\}\times A)=0, \quad
\mu^\P((0,t]\times A)<\infty, \qquad t\ge0.
\]
Let   
$\nu:\Omega\times [0,\infty)\times A\to [0,\infty)$ be a  $\Pc(\F)\otimes \calb(A) $-measurable random field.
\begin{itemize}
    \item[\textup{(i)}] The  exponential process
\begin{eqnarray}  \nonumber
\kappa_t^{\nu}(\omega) \ &=&  \exp\left(\int_0^t\int_A (1 -  \nu_s(\omega,a))\mu^\P(ds\,da)
\right)\prod_{0< T_n(\omega)\le t}\nu_{  T_n(\omega)}(  \omega,A_n(\omega)),\qquad t\ge 0,
%\label{doleans}
\end{eqnarray} 
is a supermartingale with respect to $\P$ and $\F$.
\item[\textup{(ii)}] Suppose that for some (deterministic) $T>0$ we have $\E\,\kappa_T^\nu=1$ and define
a new probability on $(\Omega,\calf)$ setting
$ \P_T^{ \nu}(d\omega)=\kappa_T^{ \nu}( \omega)\, \P(d \omega)$. Then $\kappa^\nu$ is an $\F$-martingale under $\P$ on $[0,T]$ and
the random measure $\nu_t(a)\mu^\P(dt\,da)$ 
is
the $\F$-compensator of $\mu$ under $  \P_T^{ \nu}$ on the set
$(0,T]\times A$, meaning that
for every $n$  and $C\in\calb(A)$ the processes
\begin{equation}\label{compensmartloc}
\mu((0,t\wedge T_n]\times C)-\int_{(0,t\wedge T_n] \times C}\nu_s(a)\, 
\mu(ds\,da), \qquad t\in [0,T],
\end{equation}
are $\F$-martingales under $\P_T^\nu$. 

\end{itemize}

\end{Theorem}

 \noindent {\bf Proof.} For points $(i)$ and $(ii)$ we refer to
\cite{ja} Proposition 4.3 and Theorem 4.5 (and the comment in the  lines that follow it), which even contains a more general statement. For our purposes the slightly less general result in \cite{bremaud2} Theorem 5.5.1 will be enough.
\qed

In the sequel we will often use the previous result in the following special setting.

\begin{Proposition} \label{girsanovpoisson} Let $\lambda(da)$ be a finite Borel measure on $A$ and
let $(T_n,A_n)_{n\ge 1}$ be a marked point process defined in $(\Omega,\calf,\F,\P)$  with   compensator $\mu^\P(dt\,da)=\lambda(da)\,dt$ (hence 
 a Poisson  point process with intensity $\lambda(da)$). 
Let   
$\nu:\Omega\times [0,\infty)\times A\to \R$ be a  $\Pc(\F)\otimes \calb(A) $-measurable random field 
satisfying
$$
0<\nu_t(\omega,a)\le \sup \nu<\infty.
$$
\begin{itemize}
    \item[\textup{(i)}] 
    The exponential process
\begin{eqnarray}  \nonumber
\kappa_t^{\nu}(\omega) \ &=&  \exp\left(\int_0^t\int_A (1 -  \nu_s(\omega,a))\lambda(da)\,ds
\right)\prod_{0< T_n(\omega)\le t}\nu_{  T_n(\omega)}(  \omega,A_n(\omega)),\qquad t\ge 0,%\label{doleans}
\end{eqnarray} 
is a martingale with respect to $\P$ and $\F$ and it satisfies $\E[|\kappa_t^\nu|^2]<\infty$ for every $t\ge0$.
\item[\textup{(ii)}] Given  $T>0$  define
a probability on $(\Omega,\calf)$, equivalent to $\P$, setting
$ \P^{ \nu}_T(d\omega)=\kappa_T^{ \nu}( \omega)\, \P(d \omega)$. Then 
the random measure $\nu_t(a)\,\lambda(da)\,dt$ 
is
the $\F$-compensator of $\mu$ under $  \P^{ \nu}_T$ on  
$(0,T]\times A$, and 
for every  $C\in\calb(A)$ the processes
\begin{equation}\label{compensmartvera}
\mu((0,t]\times C)-\int_{(0,t] \times C}\nu_s(a)\, 
\lambda( da)\,ds, \qquad t\in [0,T],
\end{equation}
are $\F$-martingales under $\P^\nu_T$.

\item[\textup{(iii)}] Suppose in addition that $(\Omega,\calf)$ is the canonical space introduced above and that 
 $\nu$ is  $\Pc(\F^\mu)\otimes \calb(A) $-measurable, where  $\F^\mu$ denotes  the natural filtration. Then
there exists
 a unique probability measure $\P^{\nu}$ on $(\Omega,\calf^\mu_\infty)$ such that $\P^{\nu}(d\omega)=\kappa_t^{\nu}(\omega)\P(d\omega)$
 on each $\sigma$-algebra $\calf^\mu_t$. Moreover   the $\F^\mu$-compensator of $\mu$ under $\P^{\nu}$
 is $\nu_t(a)\lambda(da)\,dt$.

\item[\textup{(iv)}] Under the assumptions of point \textup{(iii)}, recalling the form of   $\nu$    given in Proposition \ref{propnaturalfiltration}-$(iv)$, the law of $(T_n,A_n)_{n\ge1}$ under $\P^\nu$ is described by the following formulae:  $\P^\nu$-a.s., for $t\ge 0$,
\begin{align} \label{leggetunomarcato}
    \P^\nu(T_1>t) &=
    \exp\left(-\int_0^t\int_A
    \nu_s^{(0)} ( a ) \,\lambda(da)\,ds\right),
    \\\label{leggeaunomarcato}
    \P^\nu(A_1\in da\,|\, \sigma(T_1)) (\omega)
    &=
    \frac{\nu_{T_1}^{(0)} ( a ) \,\lambda(da)
    }{
    \int_A\nu_{T_1}^{(0)} ( a ) \,\lambda(da)
    },
    \\\label{leggetduemarcato}
  \P^\nu(T_{2}-T_1>t\,|\, \sigma(T_1,A_1) )(\omega)
  &=  
  \exp\bigg(-\int_{T_1(\omega)}^{T_1(\omega)+t}\int_A
   \nu_s^{(1)}\big( T_1(\omega ),
A_1(\omega),a\big)  \,\lambda(da)\,ds\bigg),
\\\label{leggeaduemarcato}
\P^\nu(A_{2}\in da\,|\, \sigma(T_1,A_1,T_{2}))(\omega) &=
 \frac{
  \nu_{T_{2}(\omega)}^{(1)}\big( T_1(\omega ),
A_1(\omega),a\big)
\,\lambda(da)
    }{
    \int_A
  \nu_{T_{2}(\omega)}^{(1)}\big( T_1(\omega ),
A_1(\omega),a\big)
\,\lambda(da)
    },
    \end{align}
     and in general, for $n\ge1$,    
\begin{align}\label{leggetnmarcato}
    &
 \P^\nu(T_{n+1}-T_n>t\,|\, \sigma(T_1,A_1,\ldots,T_n,A_n) )(\omega)=
 \\
 &\nonumber
\qquad\qquad\qquad    \exp\bigg(-\int_{T_n(\omega)}^{T_n(\omega)+t}\int_A
   \nu_s^{(n)}\big( T_1(\omega ),
A_1(\omega),\ldots,T_n(\omega),A_n(\omega),a\big)  \,\lambda(da)\,ds\bigg),
\\
&\label{leggeanmarcato}
\P^\nu(A_{n+1}\in da\,|\, \sigma(T_1,A_1,\ldots,T_n,A_n,T_{n+1})) (\omega)=\\
&\nonumber
\qquad\qquad\qquad
 \frac{
  \nu_{T_{n+1}(\omega)}^{(n)}\big( T_1(\omega ),
A_1(\omega),\ldots,T_n(\omega),A_n(\omega),a\big)
\,\lambda(da)
    }{
    \int_A
  \nu_{T_{n+1}(\omega)}^{(n)}\big( T_1(\omega ),
A_1(\omega),\ldots,T_n(\omega),A_n(\omega),a\big)
\,\lambda(da)
    }.
\end{align}
\end{itemize}

\end{Proposition}

\begin{Remark}
{\rm  
The formulae in Proposition \ref{girsanovpoisson}-$(iv)$ describe, under $\P^\nu$, the law of 
$T_1$, the conditional law of $A_1$ given $T_1$, the conditional law of $T_2$ given $T_1,A_1$ and so on, and therefore they specify uniquely the law of the entire sequence $(T_n,A_n)$. The fact that this law is entirely determined by the compensator was already  stated in  
 Proposition \ref{propcompensator}-$(v)$
\qed
}
\end{Remark}

 \noindent {\bf Proof.} Using the fact that $\nu$ is bounded and $\lambda(da)$ is finite it is not difficult to show that $\kappa^\nu$ are  $\F$-martingales under $\P^\nu$ and that $\E[|\kappa_t^\nu|^2]<\infty$ for every $t\ge0$ (or one can refer to  \cite{bremaud2} Remark 5.5.2).
In particular we have $\E\,[\kappa_T^\nu]=1$. Since $\nu>0$ we also have $\kappa^\nu>0$ and so $\P^\nu_T$ is equivalent to $\P$.
The indicated form of the compensator follows from 
 Theorem \ref{girsanovpoint} and the fact that the processes in \eqref{compensmartvera} are martingales  follows from the local martingale property of \eqref{compensmartloc} letting $n\to\infty$, using the boundedness of $\nu$ and standard properties of the Poisson process.
 
Now we turn to the
  verification  of  point $(iii)$.
  Using
the boundedness of $\nu$  one first verifies that
$$
\kappa^{\nu}_{t\wedge T_n}(\omega)\le a_n\,e^{b\,T_n(\omega)}
$$
for some constants $a_n,b$,
which implies that $(\kappa^\nu_{t\wedge T_n})_{t\ge 0}$ is a
uniformly integrable martingale with respect to $\P$ and $\F^\mu$.
Then the probabilities $\P_n^{\nu}$ defined on $\calf_{T_n}^\mu$ setting
$\P_n^{\nu}(d\omega)=\kappa_{T_n}^{\nu}(\omega)\,\P(d\omega)$
satisfy the  compatibility condition: $\P_{n+1}^{\nu}=\P_n^{\nu}$ on $\calf_{T_n}^\mu$
for every $n$, as an application of the optional sampling theorem for uniformly integrable martingales. Recalling that $\calf_{T_n}^\mu$ and $\calf_\infty^\mu$ are isomorphic to $\eqref{infproduct}$, by  the Kolmogorov extension theorem
there exists a unique probability $\P^{\nu}$ on $(\Omega,\calf^\mu_\infty)$
such that $\P^{\nu}=\P^{\nu}_n$ on each $\calf^\mu_{T_n}$.

Now let $A\in\calf^\mu_t$ and set $B_n:=A\cap \{T_n\le t<T_{n+1}\}$ for every $n\ge 0$ ($T_0:=0$). By the canonical character of the filtration $\F^\mu$ we have $B_n\in\sigma(T_0,A_0,\ldots,T_n,A_n)=\calf^\mu_{T_n}$. Since $B_n=A\cap \{T_n\le t\}\cap\{T_{n+1}\le t\}^c$, clearly $B\in\calf^\mu_t$ and so $B_n\in \calf^\mu_{T_n}\cap \calf^\mu_t=\calf^\mu_{T_n\wedge t}$. Since $A=\cup_{n\ge 0}B_n$ (disjoint union),
$$
\E[\kappa_{t}^{\nu}1_A]=\sum_{n\ge 0}\E[\kappa_{t}^{\nu}1_{B_n}]=\sum_{n\ge 0}\E[\kappa_{t\wedge T_n}^{\nu}1_{B_n}],
$$
where we used the fact that $\kappa^{\nu}$ is a martingale with respect to $\P$ and $\F^\mu$ and $B_n\in \calf^\mu_{T_n\wedge t}$. Since $(\kappa^\nu_{s\wedge T_n})_{s\ge 0}$ is a uniformly integrable martingale and $B_n\in \calf^\mu_{ t}$ we have $\E[\kappa_{t\wedge T_n}^{\nu}1_{B_n}]=\E[\kappa_{T_n}^{\nu}1_{B_n}]$, and since $B_n\in \calf^\mu_{ T_n}$ we have $\E[\kappa_{T_n}^{\nu}1_{B_n}]=\P^{\nu}(B_n)$ and finally
$$
\E[\kappa_{t}^{\nu}1_A]=\sum_{n\ge 0}\E[\kappa_{t\wedge T_n}^{\nu}1_{B_n}]= \sum_{n\ge 0}\P^{\nu}(B_n)=\P^{\nu}(A).
$$
This proves that $d\P^{\nu}=\kappa_t^{\nu}\,d\P$ on each  $\calf^\mu_t$. The fact that  the $\F^\mu$-compensator of $\mu$ under $\P^{\nu}$ has the desired form follows from the previous points.

Finally, the formulae of point $(iv)$ follow from uniqueness of the $\F^\mu$-compensator and from the  explicit form of the compensator given in \cite{ja} Proposition 3.1 or \cite{bremaud2} Theorem 5.2.2.
\qed

\subsection{Some point process constructions using time changes}
\label{subs;timechanges}

All the results of the previous subsection are already well established. Now we present some statements, that consist in 
time-change constructions, that we believe to have some interest in the theory of point processes.
Although  
time-change techniques are customary in this theory, see for instance 
\cite{bremaud2} Sections 5.4 and 10.3, the following results have been proved in \cite{BCFP16c} in connection with the randomization method.

To make the exposition more readable, and as an introduction to the more complex cases needed later, we start from a result on simple point processes. Let $(T_n)_{n\ge1}$ be a Poisson process with intensity $\lambda>0$ defined on the canonical probability space $(\Omega,\calf,\P)$. This means that $\Omega$ is the set of sequences $\omega=(t_n)_{n\ge1}\subset (0,\infty)$ such that $t_n<t_{n+1}\to\infty$, $T_n(\omega)=t_n$, $\calf$ is generated by the functions $T_n$  and $\P$ is the unique probability on $(\Omega,\calf)$ such that $(T_n-T_{n-1})_{n\ge1}$ (with $T_0:=0$) is an independent family of exponential random variables with parameter $\lambda$. Let $\F^\mu=(\calf^\mu_t)_{t\ge0}$ be the natural filtration, namely the one generated by the counting process $N_t:=\sum_{n\ge1}1_{\{T_n\le t\}}$. Note that $\calf^\mu_\infty=\calf$.
The $\F^\mu$-compensator of $(T_n)$ under $\P$ is the measure $\lambda\,dt$ on the Borel subsets of $(0,\infty)$. 

\begin{Proposition} With the previous notation, 
suppose we are given  an    $\F^\mu$-predictable process
$\nu$ 
satisfying
$$
0<\inf\nu\le \nu_t(\omega)\le \sup \nu<\infty.
$$
Let $\P^\nu$ denote the   probability on $(\Omega,\calf^\mu_\infty)$ under which the $\F^\mu$-compensator of $(T_n)$
is $\lambda\,\nu_t\,dt$, constructed in 
 Proposition  \ref{girsanovpoisson}-(iii). Then there exists a sequence $(T_n^\nu)_{n\ge1}$ of $\calf_\infty^\mu$-measurable random variables with values in $(0,\infty)$ such that
$T_n^\nu<T_{n+1} $, $T_n^\nu\to\infty$ as $n\to\infty$ and
$$
\call_{\P}((T_n^\nu)_{n\ge1})=
\call_{\P^\nu}((T_n)_{n\ge1}),
$$
namely the law of the sequence $(T_n^\nu) $ under $\P$ is the same as the law of the sequence $(T_n) $ under $\P^\nu$.
\end{Proposition}

\begin{Remark}
    {\rm 
    We note that the variables $T^\nu_n$ need not be stopping times for the filtration $\F^\mu$ and so the point process $(T^\nu_n)_{n\ge 1}$ is not a point process for $\F^\mu$ in general.
    }\qed
\end{Remark}

 \noindent {\bf Proof.} 
 Let us write $\nu$ in the form given in  Proposition \ref{propnaturalfiltration}-$(iii)$, where all the functions $\nu^{(n)}$ are bounded above by $\sup \nu$ and bounded below by $\inf \nu$.  By Proposition \ref{girsanovpoisson}-$(iv)$ (or rather its version for simple point processes) we have, for $t\ge 0$,
 \begin{align}
    \P^\nu(T_1>t) &=  \exp\left(-\lambda \int_0^t \nu_s^{(0)}   \,ds\right)
  \label{leggetunosemplice}  \\
  \P^\nu(T_{2}-T_{1}>t\,|\, \sigma(T_1) )
  &=  
  \exp\bigg(-\lambda \int_{T_1 }^{T_1 +t}
   \nu_s^{(1)}( T_1  ) \,ds\bigg),
   \label{leggetduesemplice}
\end{align}
and so on. Now define
\[
T_1^\nu= \inf\{t\ge0\,:\,   \int_0^t \nu_s^{(0)}   \,ds\ge T_1\},
\]
\[
T_2^\nu= \inf\{t\ge T_1^\nu\,:\,  \int_{T_1^\nu}^t \nu_s^{(1)}(T_1^\nu)   \,ds\ge T_2-T_1\},
\]
and in general, for $n\ge 1$,
\[
T_{n+1}^\nu= \inf\{t\ge T_n^\nu\,:\,   \int_{T_n^\nu}^t \nu_s^{(n)}(T_1^\nu,\ldots, T_n^\nu)   \,ds\ge T_{n+1}-T_n\}.
\]

Let us first consider $T_1^\nu$. 
Since $\nu^{(0)}\ge\inf\nu>0$, the function $t\mapsto  \int_0^t\nu_s^{(0)}   \,ds$ is continuous and strictly increasing, so that $T_1^\nu$ is in fact a minimum. Since $T_1\in (0,\infty)$ we deduce that  $T_1^\nu\in (0,\infty)$, $T_1^\nu\cdot\sup\nu\ge T_1$ and
\[
T^\nu_1>t\qquad \Longleftrightarrow\qquad
  \int_0^t \nu_s^{(0)}   \,ds< T_1.
\]
It follows that
\[
\P(
T^\nu_1>t)=\P\left(T_1>
  \int_0^t \nu_s^{(0)}   \,ds\right)
= \exp\left(-\lambda \int_0^t \nu_s^{(0)}   \,ds\right),
\]
where we used the fact that $\P(T_1>r)=\exp(-\lambda r)$ for $r\ge 0$. 
Comparing with 
  \eqref{leggetunosemplice}
we see that  $\P^\nu(T_1>t) =
\P(T_1^\nu >t)$.

Next let us consider
$T_2^\nu$.    
Since $\nu^{(1)}\ge\inf\nu>0$, the function $t\mapsto \int_{T_1^\nu}^t \nu_s^{(1)}(T_1^\nu)     \,ds$ is continuous and strictly increasing, and since $ T_2-T_1\in (0,\infty)$ we see that $T_2^\nu$ is in fact a minimum and  $T_2^\nu\in (T_1^\nu,\infty)$. We also have  $(T_2^\nu-T_1^\nu)\cdot \sup\nu\ge T_2-T_1$  and, recalling the previous inequality $T_1^\nu\cdot\sup\nu\ge T_1$, we deduce that $T_2^\nu\cdot\sup\nu\ge T_2$. Moreover
\[
T^\nu_2 >t\qquad \Longleftrightarrow\qquad
  \int_{T_1^\nu}^{ t} \nu_s^{(1)}(T_1^\nu)   \,ds<T_2- T_1.
\]
It follows that, for $t\ge0$,
\[
\P(
T^\nu_2-T_1^\nu>t\,|\,\sigma(T_1^\nu))=\P\left(T_2-T_1> \int_{T_1^\nu}^{T_1^\nu+t} \nu_s^{(1)}(T_1^\nu)   \,ds\,\bigg|\,\sigma(T_1^\nu)\right).
\]

Recall now that $T_{2}-T_1$ is independent of $T_1$
under $\P$ and note that, by construction, $T_1^\nu$ is $\sigma(T_1)$-measurable. Therefore, $T_{2}-T_1$ is also independent of $T_1^\nu$. As a consequence, 
\begin{align*}
&\P(
T^\nu_2-T_1^\nu>t\,|\, \sigma(T_1^\nu))
=  \P(T_2-T_1>r)\big|_{r=\int_{T_1^\nu}^{T_1^\nu+t} \nu_s^{(1)}(T_1^\nu)   \,ds} 
= 
  \exp\bigg(-\lambda \int_{T_1^\nu }^{T_1^\nu +t}
   \nu_s^{(1)}( T_1^\nu  ) \,ds\bigg),
\end{align*}
where for the last equality we used the formula $\P(T_{2}-T_1>r)=\exp(-\lambda r)$ for $r\ge 0$. 
Comparing with 
  \eqref{leggetduesemplice}
we see that  $\P(
T^\nu_2-T_1^\nu>t\,|\,\sigma(T_1^\nu)) =
\P^\nu(T_2-T_1>t\,|\,\sigma(T_1))$.

Iterating this argument, or by a formal proof by induction, one concludes that, for every $n\ge 1$,
$T_n^\nu\cdot\sup\nu\ge T_n$, which implies that $T_n^\nu\to\infty$, and
\[
\P(
T^\nu_{n+1}-T_n^\nu>t\,|\,\sigma(T_1^\nu,\ldots,T_n^\nu)) =
\P^\nu(T_{n+1}-T_n>t\,|\,\sigma(T_1,\ldots,T_n)),
\]
which proves the final assertion of the proposition.
\qed

Next we move to a more general version of this result, in the case of marked point processes. We first consider the special case  $A=\R$ with a  non atomic measure $\lambda(da)$. Our framework is similar to Proposition \ref{girsanovpoisson}.

\begin{Proposition}
\label{cambiosequenzasur}
Let $\lambda(da)$ be a finite Borel measure on $A=\R$ such that $\lambda(\{a\})=0$ for every $a\in A$. Let $(\Omega,\calf)$ be the canonical space of a marked point process  $(T_n,A_n)_{n\ge 1}$ on $A$ and let $\P$ be the unique probability on $(\Omega,\calf)$ such that $(T_n,A_n)_{n\ge 1}$ is  a Poisson process with intensity $\lambda(da)$.
Denote by $\F^\mu=(\calf^\mu_t)_{t\ge0}$   its natural filtration (in particular, $\calf^\mu_\infty=\calf$) and recall that the $ \F^\mu$-compensator is $\lambda(da)\,dt$.

Let   
$\nu:\Omega\times [0,\infty)\times A\to [0,\infty)$ be a  $\Pc(\F^\mu)\otimes \calb(A) $-measurable random field 
satisfying
$$
0<\inf \nu\le  \nu_t(\omega,a)\le \sup \nu<\infty.
$$
Define the $\F^\mu$-martingale
\begin{eqnarray}  \nonumber
\kappa_t^{\nu} \ &=&  \exp\left(\int_0^t\int_A (1 -  \nu_s( a))\lambda(da)\,ds
\right)\prod_{0< T_n \le t}\nu_{  T_n }(  A_n ),\qquad t\ge 0,
\end{eqnarray} 
and let 
$\P^{\nu}$ be the
  unique  probability measure  on $(\Omega,\calf^\mu_\infty)$ such that $\P^{\nu}(d\omega)=\kappa_t^{\nu}(\omega)\P(d\omega)$
 on each $\sigma$-algebra $\calf^\mu_t$.

 Then there exists a sequence $(T_n^\nu,A_n^\nu)_{n\ge1}$ of $\calf_\infty^\mu$-measurable random variables with the following properties.
  
\begin{itemize}   

\item[\textup{(i)}] For every $n\ge1$,
$(T_n^\nu,A_n^\nu)$ takes
  values in $(0,\infty)\times A$ and 
$T_n^\nu<T_{n+1}^\nu $.
\item[\textup{(ii)}] 
  $T_n^\nu\to\infty$ as $n\to\infty$.

\item[\textup{(iii)}] 
We have
$$
\call_{\P}((T_n^\nu,A_n^\nu)_{n\ge1})=
\call_{\P^\nu}((T_n,A_n)_{n\ge1}),
$$
namely the law of the sequence $(T_n^\nu,A_n^\nu) $ under $\P$ is the same as the law of the sequence $(T_n) $ under $\P^\nu$.
\end{itemize}
\end{Proposition}

 \noindent {\bf Proof.} 
Let us first introduce the cumulative distribution function of the probability $\lambda(da)/\lambda(A)$:
\begin{equation}
    \label{fdrlambda}
    F(b) =    \frac{  \lambda((-\infty,b])}{\lambda(A)},
\qquad b\in A=\R.
\end{equation}
We recall that $F(b) =  \P(A_n\leq b)$ for every $n$. 
Since we assume that $\lambda$ is nonatomic, $F$ is continuous and therefore the random variables $F(A_n)$ have law $U(0,1)$, the uniform distribution on $(0,1)$, under $\P$.

Next let us note that $\nu$ has the form given in  Proposition \ref{propnaturalfiltration}-$(iv)$, where all the functions $\nu^{(n)}$ are bounded above by $\sup \nu$ and bounded below by $\inf \nu$. 
 Recall that 
the law of $(T_n,A_n)_{n\ge1}$ under $\P^\nu$ is described by the formulae in
 Proposition \ref{girsanovpoisson}-$(iv)$. We will construct iteratively $T_1^\nu, A_1^\nu, T_2^\nu, A_2^\nu,\ldots$ and we will check at each step that their laws (or conditional laws) under $\P$ are given by the same formulae: this will prove (iii). At each step we will verify that $T_n^\nu<T_{n+1}^\nu $ and that $T_n^\nu\cdot \sup\nu\ge T_n$,
which implies that  $T_n^\nu\to\infty$ as $n\to\infty$.

First define
\begin{equation}\label{Tnu1}
T^\nu_1  = \inf\{t\ge 0\,:\,\theta^{(1)}_t\ge T_1\}, \qquad {\rm where} \qquad 
\theta^{(1)}_t=
\frac{1}{\lambda(A)} \int_0^t\int_A \nu_s^{(0)}(a) \lambda(da)\,ds.
\end{equation}
Since $\nu^{(0)}\geq\inf\nu>0$, the function
$ t\mapsto \theta^{(1)}_t$
is continuous and strictly increasing. Since $T_1\in (0,\infty)$
we see 
 that the infimum is achieved and
 $T_1^\nu\in (0,\infty)$ is in fact a minimum.  
 Since $\nu^{(0)}\leq\sup\nu$
 we also have
 $T_1^\nu\cdot\sup\nu\ge T_1$.  Finally note that
\[
T^\nu_1>t\qquad \Longleftrightarrow\qquad
   \theta^{(1)}_t
  < T_1.
\]
It follows that, for $t\ge0$,
\[
\P(
T^\nu_1>t)=\P\left(T_1>
  \theta^{(1)}_t\right)
= \exp\left(-\lambda(A) \theta^{(1)}_t\right)
= \exp\left(- \int_0^t\int_A \nu_s^{(0)}(a) \lambda(da)\,ds\right),
\]
where we used the fact that $\P(T_1>r)=\exp(-\lambda(A) r)$ for $r\ge 0$. 
Comparing with 
  \eqref{leggetunomarcato}
we see that  $\P^\nu(T_1>t) =
\P(T_1^\nu >t)$.

Now we define $A_1^\nu$ as follows:
\begin{equation}\label{Anu1}
A^\nu_1  = \inf\{b\in A=\R\,:\,F_b^{(1)} \ge F(A_1)\},
\qquad {\rm where} \qquad 
F_b^{(1)} =   \frac{\int_{-\infty}^b \nu_{T_1^\nu }^{(0)}( a)\,\lambda(da)}{\int_{-\infty}^{+\infty}\nu_{T_1^\nu }^{(0)}( a)\,\lambda(da)}
\end{equation}
and $F$ was introduced in \eqref{fdrlambda}.
Since $\nu^{(0)}\ge\inf\nu>0$ and $\lambda(\{a\})=0$ for any $a\in A$, the function $b\mapsto F_b^{(1)} $ is continuous and strictly increasing. Moreover $\lim_{b\rightarrow-\infty}F_b^{(1)}=0$, $\lim_{b\rightarrow+\infty}F_b^{(1)}=1$ and $F(A_1)$ takes values   in $(0,1)$ a.s.  Therefore
  the infimum defining $A_1^\nu$ is achieved,
 $A_1^\nu$ is in fact a minimum  and 
 $A_1^\nu \leq b$ if and only if $F(A_1)\leq F_b^{(1)}$.
 It follows that
\begin{align*}
\P\big(A_1^\nu \leq b\,\big|\,\sigma(T_1^\nu)\big)  
=   \P\big(F(A_1)\leq F_b^{(1)}\,\big|\,\sigma(T_1^\nu\big).
\end{align*}
Since $A_1$ is independent of $T_1$ under $\P$
and $T_1^\nu$ is $\sigma(T_1)$-measurable by construction, it follows that $A_1$ is also independent of $T_1^\nu$. Moreover, by definition we see that $F_b^{(1)}$ is $\sigma(T_1^\nu)$-measurable. Therefore we have, $\P$-a.s.,
\[
\P\big(F(A_1)\leq F_b^{(1)}\,\big|\,\sigma(T_1^\nu)\big)  =  \P\big(F(A_1)\leq a\big)\big|_{a=F_b^{(1)}}   =   F_b^{(1)},
\]
where we used the fact that $F(A_1)$ has law $U(0,1)$ under $\P$. From the equality of the cumulative distribution functions we deduce
\[
\P\big(A_1^\nu \in da\,\big|\,\sigma(T_1^\nu)\big)  
=
  \frac{  \nu_{T_1^\nu }^{(0)}( a)\,\lambda(da)}{\int_{A}\nu_{T_1^\nu }^{(0)}( a)\,\lambda(da)}.
\]
Comparing with 
  \eqref{leggeaunomarcato}
we conclude that  $\P(A_1^\nu \in da\,|\,\sigma(T_1^\nu))=\P^\nu(A_1\in da\,|\,\sigma(T_1)) $.

Next we define $T_2^\nu$ as follows:
\begin{equation}\label{Tnu2}
T^\nu_2  = \inf\{t\ge T^\nu_1\,:\,\theta^{(2)}_t\ge T_2-T_1\},
\qquad {\rm where} \qquad
\theta^{(2)}_t=
\frac{1}{\lambda(A)} \int_{T_1^\nu}^{t}\int_A \nu_s^{(1)}(T_1^\nu,A_1^\nu,a) \,\lambda(da)\,ds .
\end{equation}
Since $\nu^{(1)}\geq\inf\nu>0$, the function
$ t\mapsto \theta^{(2)}_t$
is continuous and strictly increasing, and since $ T_2-T_1\in (0,\infty)$  the infimum is achieved,  $T_2^\nu$ is  a minimum and $T_2^\nu\in (T_1^\nu,\infty)$. Since $\nu^{(1)}\leq\sup\nu$ we have $(T_2^\nu-T_1^\nu)\cdot \sup\nu\ge T_2-T_1$, and recalling the previous inequality $T_1^\nu\cdot\sup\nu\ge T_1$, we conclude that $T_2^\nu\cdot\sup\nu\ge T_2$. Moreover, the inequality 
$T^\nu_2>t$ holds if and only if 
$  \theta^{(2)}_{t} <T_2- T_1$ and it follows that,  for $t\ge0$,
\[
\P(
T^\nu_2-T_1^\nu>t\,|\,\sigma(T_1^\nu,A_1^\nu))=\P\left(T_2-T_1> \theta^{(2)}_{T_1^\nu+t}\,\bigg|\,\sigma(T_1^\nu,A_1^\nu)\right).
\]
Note that, by construction, $(T_1^\nu,A_1^\nu)$ is $\sigma(T_1,A_1)$-measurable and therefore $\theta^{(2)}_{T_1^\nu+t}$ is also $\sigma(T_1,A_1)$-measurable.
Recall now that $T_{2}-T_1$ is independent of $(T_1,A_1)$
under $\P$. We conclude that  $T_{2}-T_1$ is also independent of $(T_1^\nu,A_1^\nu,\theta^{(2)}_{T_1^\nu+t})$. As a consequence, 
\begin{align*}
&\P(
T^\nu_2-T_1^\nu>t\,|\, \sigma(T_1^\nu,A_1^\nu))
=  \P(T_2-T_1>r)\big|_{r=\theta^{(2)}_{T_1^\nu+t}} 
= 
  \exp\bigg(-\int_{T_1^\nu}^{T_1^\nu+t}\int_A \nu_s^{(1)}(T_1^\nu,A_1^\nu,a) \,\lambda(da)\,ds\bigg)
\end{align*}
where for the last equality we used the formula $\P(T_{2}-T_1>r)=\exp(-\lambda(A) r)$ for $r\ge 0$. 
Comparing with 
  \eqref{leggetduemarcato}
we see that  $\P(
T^\nu_2-T_1^\nu>t\,|\,\sigma(T_1^\nu,A_1^\nu)) =
\P^\nu(T_2-T_1>t\,|\,\sigma(T_1,A_1))$.

Next we define $A_2^\nu$ by the formula
\begin{equation}\label{Anu2}
A^\nu_2  = \inf\{b\in A=\R\,:\,F_b^{(2)} \ge F(A_2)\},
\qquad
{\rm where} \qquad 
F_b^{(2)} =   \frac{\int_{-\infty}^b \nu_{T_2^\nu }^{(1)}(T_1^\nu, A_1^\nu, a)\,\lambda(da)}{\int_{-\infty}^{+\infty}\nu_{T_2^\nu }^{(1)}( T_1^\nu, A_1^\nu, a)\,\lambda(da)}
\end{equation}
and $F$ was introduced in \eqref{fdrlambda}.
Since $\nu^{(1)}\ge\inf\nu>0$ and $\lambda(\{a\})=0$ for any $a\in A$, the function $b\mapsto F_b^{(2)} $ is continuous and strictly increasing. Moreover $\lim_{b\rightarrow-\infty}F_b^{(2)}=0$, $\lim_{b\rightarrow+\infty}F_b^{(2)}=1$ and $F(A_2)$ takes values   in $(0,1)$ a.s.  Therefore
  the infimum defining $A_2^\nu$ is achieved,
 $A_2^\nu$ is in fact a minimum  and 
 $A_2^\nu \leq b$ if and only if $F(A_2)\leq F_b^{(2)}$.
 It follows that
\begin{align*}
\P\big(A_2^\nu \leq b\,\big|\,\sigma(T_1^\nu,A_1^\nu,T_2^\nu)\big)
=   \P\big(F(A_2)\leq F_b^{(2)}\,\big|\,\sigma(T_1^\nu,A_1^\nu,T_2^\nu)\big).
\end{align*}
Since $A_2$ is independent of $(T_1,A_1,T_2)$ under $\P$
and $(T_1^\nu,A_1^\nu,T_2^\nu)$ is $\sigma(T_1,A_1,T_2)$-measurable by construction, it follows that $A_2$ is also independent of $(T_1^\nu,A_1^\nu,T_2^\nu)$. Moreover, by definition we see that $F_b^{(2)}$ is $\sigma(T_1^\nu,A_1^\nu,T_2^\nu )$-measurable. Therefore we have, $\P$-a.s.,
\[
\P\big(F(A_2)\leq F_b^{(2)}\,\big|\,\sigma(T_1^\nu,A_1^\nu,T_2^\nu)\big)  =  \P\big(F(A_2)\leq a\big)\big|_{a=F_b^{(2)}}   =   F_b^{(2)},
\]
where we used the fact that $F(A_2)$ has law $U(0,1)$ under $\P$. From the equality of the cumulative distribution functions we deduce
\[
\P\big(A_2^\nu \in da\,\big|\,\sigma(T_1^\nu,A_1^\nu,T_2^\nu)\big)  
= \frac{  \nu_{T_2^\nu }^{(1)}(T_1^\nu, A_1^\nu, a)\,\lambda(da)}{\int_{A} \nu_{T_2^\nu }^{(1)}( T_1^\nu, A_1^\nu, a)\,\lambda(da)}.
\]
Comparing with 
  \eqref{leggeaduemarcato}
we conclude that  $\P(A_2^\nu \in da\,|\,\sigma(T_1^\nu,A_1^\nu,T_2^\nu))=\P^\nu(A_2\in da\,|\,\sigma(T_1,A_1,T_2)) $.

In general we define, for $n\ge1$,
\begin{align}
    \nonumber
T^\nu_{n+1}  = \inf\{t\ge T^\nu_n\,:\,\theta^{(n+1)}_t\ge T_{n+1}-T_n\},\quad
A^\nu_{n+1}  = \inf\{b\in A=\R\,:\,F_b^{(n+1)} \ge F(A_{n+1})\},
\end{align}
where
\begin{align}
   \nonumber
\theta^{(n+1)}_t&=
\frac{1}{\lambda(A)} \int_{T_n^\nu}^{t}\int_A \nu_s^{(n)}(T_1^\nu,A_1^\nu,\ldots,
T_n^\nu,A_n^\nu,
a) \,\lambda(da)\,ds ,
\\ \nonumber
F_b^{(n+1)} &=   \frac{\int_{-\infty}^b \nu_{T_{n+1}^\nu }^{(n)}(T_1^\nu, A_1^\nu,\ldots,T_n^\nu,A_n^\nu, a)\,\lambda(da)}{\int_{-\infty}^{+\infty}\nu_{T_{n+1}^\nu }^{(n)}( T_1^\nu, A_1^\nu,\ldots,T_n^\nu,A_n^\nu, a)\,\lambda(da)}.
\end{align}
Iterating the previous arguments, or writing down a formal proof by induction, one proves that $T_n^\nu<T_{n+1}^\nu $ and that $T_n^\nu\cdot \sup\nu\ge T_n$,
which implies that  $T_n^\nu\to\infty$ as $n\to\infty$, as well as
\[
\P(
T^\nu_{n+1}-T_n^\nu>t\,|\,\sigma(T_1^\nu,A_1^\nu,\ldots,T_n^\nu,A_n^\nu))
=
\int_{T_n^\nu}^{T_n^\nu+t}\int_A \nu_s^{(n)}(T_1^\nu,A_1^\nu,\ldots,
T_n^\nu,A_n^\nu,
a) \,\lambda(da)\,ds ,
\]
 \[
\P\big(A_{n+1}^\nu \in da\,\big|\,\sigma(T_1^\nu,A_1^\nu,\ldots,T_n^\nu,A_n^\nu,T_{n+1}^\nu)\big)  
= \frac{  \nu_{T_{n+1}^\nu }^{(n)}(T_1^\nu, A_1^\nu,\ldots,T_n^\nu,A_n^\nu, a)\,\lambda(da)}{\int_{A} \nu_{T_{n+1}^\nu }^{(n)}( T_1^\nu, A_1^\nu,\ldots,T_n^\nu,A_n^\nu, a)\,\lambda(da)}.
\]
Comparing with
\eqref{leggetnmarcato},
\eqref{leggeanmarcato} we see that, for $n\ge1$ and $t\ge0$,   
\[  
\P(
T^\nu_{n+1}-T_n^\nu>t\,|\,\sigma(T_1^\nu,A_1^\nu,\ldots,T_n^\nu,A_n^\nu))=
 \P^\nu(T_{n+1}-T_n>t\,|\, \sigma(T_1,A_1,\ldots,T_n,A_n) ),
\]
\[
\P\big(A_{n+1}^\nu \in da\,\big|\,\sigma(T_1^\nu,A_1^\nu,\ldots,T_n^\nu,A_n^\nu,T_{n+1}^\nu)\big)=
\P^\nu(A_{n+1}\in da\,|\, \sigma(T_1,A_1,\ldots,T_n,A_n,T_{n+1}),
\]
which completes the proof of the proposition.
\qed

We finally consider the more general case when the space $A$ is an arbitrary Borel space and $\lambda(da)$ is an arbitrary finite measure on $\calb(A)$. We will reduce to the previous situation by means of the following lemma.
 
\begin{Lemma}
    \label{lambdaliftato}
    Let $\lambda(da)$ be a finite  measure on $(A,\calb(A))$. Then there exists a finite measure $\lambda' (dr)$ on $(\R,\calb(\R))$ and a Borel map $\pi:\R\to A$ such that $\lambda=\lambda' \circ\pi^{-1}$  and
    $\lambda' (\{r\})=0$ for every $r\in \R$ 
    \end{Lemma}

 \noindent {\bf Proof.} 
Recall that $A$ is assumed to be a Borel space and it is  known  that  any such  space is either finite or countable (with the discrete topology) or isomorphic, as a measurable space, to the real line (or equivalently
to the half line $(0,\infty)$): see e.g. \cite{BertsekasShreve78}, Corollary 7.16.1.
Let us denote by $A_c$ the subset of $A$ consisting of all points
$a\in A$ such that $\lambda(\{a\})>0$, and let
$A_{nc}=A\backslash A_c$. Since $\lambda$ is finite, the set
$A_c$ is either empty or at most countable, and it follows in particular
that both $A_c$ and $A_{nc}$ are also Borel spaces.  

Suppose initially that $A_{c}\neq\emptyset$ and $A_{nc}\neq\emptyset$.
Choose an arbitrary nonatomic finite measure on $(\R,\calb(\R))$  (say, the standard Gaussian measure, denoted by $\gamma$).
For every
$j\in A_c$ choose a (nontrivial) interval   $\Ic_j\subset (-\infty,0]$
in such a way that $\{\Ic_j,\,j\in A_c\}$ is a partition of $(-\infty,0]$.
Moreover, there exists a
  bijection $\pi_1:(0,\infty)\to A_{nc}$ such that $\pi_1$ and its inverse are both
  Borel measurable.
Denote by  $\lambda' $  the unique positive measure  on $(\R,\Bc(\R))$ such that
\begin{eqnarray*}
\lambda' (B) \ = \ \lambda(\{j\})\gamma(B)/\gamma(\Ic_j),
& \text{for every }B\subset \Ic_j,\,  B\in\calb(\R),\,j\in A_c,
\nonumber\\
\lambda' (B)=\lambda (\pi_1(B))&
 \text{for every }B\subset (0,\infty),\, B\in\calb(\R).
\end{eqnarray*}
Then $\lambda' $ is a finite measure
 satisfying $\lambda'(\Ic_j)=\lambda(\{j\})$ for every $j\in A_c$
and $\lambda'(\{r\})=0$ for every $r\in\R$.
We also define $\pi\colon\R\rightarrow A$ given by
\begin{equation}\label{defproiezdue}
\pi(r) \ = \
\begin{cases}
j, \qquad & \text{if }
r\in \Ic_j  \text{ for some } j\in A_c, \\
\pi_1(r),
& \text{if }
r\in (0,\infty),
\end{cases}
\end{equation}
and it is immediate to check that  $\lambda$ $=$ $\lambda'  \circ \pi^{-1}$.

In the remaining cases we slightly modify the previous construction in an obvious way: when $A_{nc}=\emptyset$ we cover the whole $\R$ with intervals $\Ic_j$, $j\in A=A_c$, and we do not introduce $\pi_1$; when $A_{c}=\emptyset$ we identify $A=A_{nc}$ with $\R$ by a Borel isomorphism $\pi$.
\qed

We finally arrive at the statement we need for the randomization method.

\begin{Proposition}
\label{cambiosequenzasua}
Let $\lambda(da)$ be a finite Borel measure on the Borel space $A$. Let $\pi:\R\to A$ and $\lambda' (dr)$ denote
the map and the measure on $\R$ constructed in Lemma 
\ref{lambdaliftato}. Let $(\Omega,\calf)$ be the canonical space for a marked point process   $(T_n,R_n)_{n\ge 1}$ on $\R$ and let $\P$ be the unique probability on $(\Omega,\calf)$ such that  $(T_n,R_n)_{n\ge 1}$ is a Poisson  process on $\R$ with intensity $\lambda' (dr)$.
Define 
\[A_n=\pi(R_n)
\]
and 
denote by $\F^\mu$   the natural filtration generated by  $(T_n,A_n)_{n\ge 1}$.

Let   
$\nu:\Omega\times [0,\infty)\times A\to \R$ be a  $\Pc(\F^\mu)\otimes \calb(A) $-measurable random field 
satisfying
$$
0<\inf\nu\le \nu_t(\omega,a)\le \sup \nu<\infty
$$
and define 
\begin{eqnarray}  \label{kappanuaux}
\kappa_t^{\nu} \ &=&  \exp\left(\int_0^t\int_A (1 -  \nu_s( a))\lambda(da)\,ds
\right)\prod_{0< T_n \le t}\nu_{  T_n }(  A_n ),\qquad t\ge 0.
\end{eqnarray} 
Then
\begin{itemize}   

\item[\textup{(a)}]
$\kappa^{\nu}$
is an $\F^\mu$-martingale.
\item[\textup{(b)}]
There exists a   
  unique   probability measure $\P^{\nu}$ on $(\Omega,\calf^\mu_\infty)$   such that $\P^{\nu}(d\omega)=\kappa_t^{\nu}(\omega)\P(d\omega)$
 on each $\sigma$-algebra $\calf^\mu_t$.
\item[\textup{(c)}]
The $\F^\mu$-compensator of $(T_n,A_n)_{n\ge 1}$ under $\P^{\nu}$ is $\nu_t(a)\,\lambda(da)\,dt$.
\end{itemize}

Moreover,  there exists a sequence $(T_n^\nu,A_n^\nu)_{n\ge1}$ of $\calf$-measurable random variables with the following properties.
\begin{itemize}   

\item[\textup{(i)}] For every $n\ge1$,
$(T_n^\nu,A_n^\nu)$ takes
  values in $(0,\infty)\times A$ and 
$T_n^\nu<T_{n+1}^\nu $.
\item[\textup{(ii)}] 
  $T_n^\nu\to\infty$ as $n\to\infty$.

\item[\textup{(iii)}] 
We have
$$
\call_{\P}((T_n^\nu,A_n^\nu)_{n\ge1})=
\call_{\P^\nu}((T_n,A_n)_{n\ge1}),
$$
namely the law of the sequence $(T_n^\nu,A_n^\nu) $ under $\P$ is the same as the law of the sequence $(T_n,A_n) $ under $\P^\nu$.
\end{itemize}
\end{Proposition}

 \noindent {\bf Proof.} 
We first note that points (a)-(b)-(c) do not follow from Proposition
\ref{cambiosequenzasur} because  $(\Omega,\calf)$ is not the canonical space of the process $(T_n,A_n)_{n\ge 1}$.

Recalling that $\lambda=\lambda' \circ\pi^{-1}$ it is immediate to see that $(T_n,A_n)_{n\ge1}$ is a Poisson process in $A$ with intensity $\lambda(da) $,  defined in $(\Omega,\calf,\P)$. We may consider
the  random measures $\mu=\sum_{n\ge1}\delta_{(T_n,A_n)}$ on the Borel subsets of $(0,\infty)\times A$
and $\mu'=\sum_{n\ge1}\delta_{(T_n,R_n)}$
on the Borel subsets of $(0,\infty)\times \R$, as well as  the corresponding natural filtrations
$\F^\mu=(\calf_t^\mu)_{t\ge0}$ and
$\F^{\mu'}=(\calf_t^{\mu'})_{t\ge0}$, 
where
\[
\calf_t^\mu=\sigma\Big( \mu((0,s]\times C)\,:\, s\in [0,t],C\in\calb(A)\Big),
\quad
\calf_t^{\mu'}=\sigma\Big( \mu'((0,s]\times B)\,:\, s\in [0,t],B\in\calb(\R)\Big).
\]
Clearly, $
\calf_t^\mu\subset 
\calf_t^{\mu'}$ and $
\calf_\infty^\mu\subset 
\calf_\infty^{\mu'}=\calf$.
Now define 
\begin{equation}
    \label{fieldaux}
\nu_t'(\omega,r)=\nu_t(\omega,\pi(r)),
\qquad 
\omega\in \Omega,\,t\ge 0,\, r\in \R.
\end{equation}
A simple monotone class argument, or the representation in 
 Proposition \ref{propnaturalfiltration}-(iv), shows that $\nu'$ is $\Pc(\F^{\mu'})\otimes \calb(\R) $-measurable. Then
we can apply 
Proposition \ref{cambiosequenzasur} to the process $(T_n,R_n)$ and the random field $\nu'$:  setting
\begin{align*}\kappa_t^{\nu'} = \exp\left(\int_0^t\int_\R (1 -  \nu_s'( r))\lambda'(dr)\,ds
\right)\prod_{0< T_n \le t}\nu'_{  T_n }(  R_n ), \qquad t\ge0,
\end{align*}
we deduce that 
 $\kappa^{\nu'}$ is an
$\F^{\mu'}$-martingale under $\P$
and there exists a unique probability measure
$\P^{\nu'}$   
     on $(\Omega,\calf^{\mu'}_\infty)=(\Omega,\calf)$
such that 
 $\P^{\nu'}(d\omega)=\kappa_t^{\nu'}(\omega)\P(d\omega)$
 on each $\sigma$-algebra $\calf^{\mu'}_t$; moreover, the $\F^{\mu'}$-compensator of $\mu'$ under $\P^{\nu'}$ is $\nu'_t(r)\,\lambda'(dr)\,dt$. 
Next we note that, since $\lambda=\lambda' \circ\pi^{-1}$, we have
\[
\int_A (1 -  \nu_s( a))\lambda(da)=
\int_\R (1 -  \nu_s( \pi(r)))\lambda'(dr)=
\int_\R (1 -  \nu_s'( r))\lambda'(dr)
\]
and since $\nu_{  T_n }(  A_n )=\nu_{  T_n }(  \pi(R_n) )=\nu'_{  T_n }(  A_n )$ we conclude that in fact the process $\kappa^{\nu'}$ and the process 
$\kappa^{\nu}$ defined in \eqref{kappanuaux} are exactly the same. 
Define $\P^\nu$ to be the restriction of $\P^{\nu'}$ to $\calf^\mu_\infty$.
Then clearly 
$\P^{\nu}(d\omega)=\kappa_t^{\nu}(\omega)\P(d\omega)$
 on each $\sigma$-algebra $\calf^\mu_t$. Since 
 $\kappa^\nu$ is $\F^\mu$-adapted, it is an
$\F^\mu$-martingale under 
$\P$. Next we prove that the $\F^\mu$-compensator of $(T_n,A_n)_{n\ge 1}$ under $\P^{\nu}$ is $\nu_t(a)\,\lambda(da)\,dt$.
We note that this is an $\F^\mu$-predictable measure, so by 
Proposition \ref{propcompensator}-(vii) it is enough to show that 
for any   random field 
$H:\Omega\times [0,\infty)\times A\to [0,\infty)$ which is $\Pc(\F^\mu)\otimes \calb(A) $-measurable we have
\begin{equation}\label{eqcompsuar}
    \E^\nu\int_{0}^\infty\int_{ A}H_t(a)\, \mu(dt\,da)=
\E^\nu\int_{0}^\infty\int_{ A}H_t(a)\, \nu_t(a)\,\lambda(da)\,dt.
\end{equation}
This equality is the same if we take the expectation $
\E^{\nu'}$.
Recalling the description of 
$\Pc(\F^\mu)\otimes \calb(A) $-measurable random fields given in
Proposition \ref{propnaturalfiltration}-(iv), it is enough to consider $H$ of the form
\begin{align*}
H_t( a)  =
 h_t \big( T_1,A_1,\ldots,T_n,A_n,a\big) \, 1_{\{T_n<t\leq T_{n+1}\}},  
\end{align*}
for some $n\ge0$ ($T_0:=0)$ and
some   map  Borel measurable $h\colon \R_+\times(\R_+\times A)^n\times A\rightarrow\R_+$.
Then
\begin{align*}
\E^{\nu'}\int_{0}^\infty\int_{ A}H_t(a)\, \mu(dt\,da)&=
\E^{\nu'}h_{T_{n+1}} \big( T_1,A_1,\ldots,T_n,A_n,A_{n+1}\big)
\\&=
\E^{\nu'}h_{T_{n+1}} \big( T_1,\pi(R_1),\ldots,T_n,\pi(R_n),\pi(R_{n+1})\big)
\\&=
\E^{\nu'}\int_{0}^\infty\int_{ \R}
h_t \big( T_1,\pi(R_1),\ldots,T_n,\pi(R_n),\pi(r)\big) \, 1_{\{T_n<t\leq T_{n+1}\}}
\, \mu'(dt\,dr).
\end{align*}
Since the integrand is a
$\Pc(\F^{\mu'})\otimes \calb(\R) $-measurable random field and since the 
$\F^{\mu'}$-compensator of $\mu'$ under $\P^{\nu'}$ is $\nu'_t(r)\,\lambda'(dr)\,dt=\nu_t(\pi(r))\,\lambda'(dr)\,dt$ 
we get
\begin{align*}
&\E^{\nu'}\int_{0}^\infty\int_{ A}H_t(a)\, \mu(dt\,da)
\\&\qquad =
\E^{\nu'}\int_{0}^\infty\int_{ \R}
h_t \big( T_1,\pi(R_1),\ldots,T_n,\pi(R_n),\pi(r)\big) \, 1_{\{T_n<t\leq T_{n+1}\}}
\,\nu_t(\pi(r))\,\lambda'(dr)\,dt
\\&\qquad =
\E^{\nu'}\int_{0}^\infty\int_{ A}
h_t \big( T_1,A_1 ,\ldots,T_n,A_n,a\big) \, 1_{\{T_n<t\leq T_{n+1}\}}
\,\nu_t(a)\,\lambda(da)\,dt
\\&\qquad =
\E^{\nu'}\int_{0}^\infty\int_{ A}
H_t(a)
\,\nu_t(a)\,\lambda(da)\,dt,
\end{align*}
where the second equality holds because $\lambda=\lambda'\circ \pi^{-1}$. 
We have now verified \eqref{eqcompsuar}
and we have proved (a)-(b)-(c) of the Proposition.
 
 We finish the proof by constructing the required sequence $(T_n^\nu,A_n^\nu)_{n\ge1}$.
 Applying 
Proposition
\ref{cambiosequenzasur}-(i)-(ii)-(iii) to the process $(T_n,R_n)_{n\ge1}$, we know  that there exists a sequence 
$(T_n^{\nu'},R_n^{\nu'})_{n\ge1}$ of functions defined in $\Omega$, measurable with respect to $\calf_\infty^{\mu'}=\calf$,  taking
  values in $(0,\infty)\times \R$, such that and 
$T_n^{\nu'}<T_{n+1}^{\nu'}\to\infty$ and
\begin{align}
    \label{eqleggiantn}
    \call_{\P}((T_n^{\nu'},R_n^{\nu'})_{n\ge1})=
\call_{\P^{\nu'}}((T_n,R_n)_{n\ge1}).
\end{align}
We define the   sequence $(T_n^\nu,A_n^\nu)_{n\ge1}$ setting
\begin{equation}\label{anproiettati} 
T_n^\nu=T_n^{\nu'}, \qquad
A_n^\nu=\pi(R_n^{\nu'}).    
\end{equation}
Recalling that $A_n=\pi(R_n)$ it follows from \eqref{eqleggiantn} that
\begin{align}\nonumber
    %\label{eqleggiantn}
    \call_{\P}((T_n^{\nu},A_n^{\nu})_{n\ge1})=
\call_{\P^{\nu'}}((T_n,A_n)_{n\ge1}).
\end{align}
Since $(T_n^{\nu},A_n^{\nu})$ are $\calf^\mu_\infty$-measurable we have
\begin{align}\nonumber
    \call_{\P}((T_n^{\nu},A_n^{\nu})_{n\ge1})=
\call_{\P^{\nu}}((T_n,A_n)_{n\ge1})
\end{align}
and the proof is finished.
\qed

\begin{Remark}
{\rm
Note that, in contrast to Proposition
\ref{cambiosequenzasur}, the random variables $T_n^\nu$ and $A_n^\nu$ are only measurable with respect to $\calf=\calf_\infty^{\mu'}$, not $\calf_\infty^{\mu}$. The marked point process $(T_n^\nu,A_n^\nu)_{\ge1}$ in general is not a point process with respect neither to $\F^{\mu}$ nor to   $\F^{\mu'}$.
}\qed
\end{Remark}

\section{A classical optimal control problem and its randomized version}
\label{S:Formulation}

\subsection{Basic notation and assumptions}
\label{SubS:Notation}

In the following we will consider controlled stochastic equations  of the form
\begin{eqnarray}  \label{dynX1}
dX_t^\alpha & = & b(t, X_t^\alpha, \alpha_t)\,dt +
\sigma( t,X_t^\alpha, \alpha_t)\,dW_t,
\end{eqnarray} 
for $t\in [0,T]$, where $T>0$ is a fixed deterministic and finite
terminal time, and reward functional
\begin{eqnarray*}
J(\alpha) & = & \E\Big[\int_0^Tf(t,X_t^\alpha,\alpha_t)\,dt+g(X^\alpha_T)\Big].
\end{eqnarray*}
The initial condition in \eqref{dynX1} is $X_0^\alpha=x_0$, a given point in $\R^n$.  The controlled process $X^\alpha$ takes values in $\R^n$ while $W$ is a  Wiener process in $\R^{d}$. 
The control process, denoted by $\alpha$, takes values in a set $A$ of control actions, that will be assumed to be a Borel space.
  We  recall that a  Borel  space   $A$ is a topological space  homeomorphic to a  Borel subset of a Polish space.
When needed, $A$ will be  endowed with its Borel $\sigma$-algebra  $\calb(A)$.

Let us state our standing  assumptions on our data.

\vspace{3mm}

\noindent {\bf (A1)}
\begin{itemize}
\item [(i)]  $A$ is a  Borel   space.
\item [(ii)]
The functions $b,\sigma,f$ are defined on $[0,T]\times \R^n\times  A$ and take 
values in $\R^n$, $\R^{n\times d }$, $\R$ respectively; they  are assumed
to be Borel measurable.
\item [(iii)]
 The function $g:\R^n\to \R$ is continuous. For fixed $x\in \R^n$ and $t\in [0,T]$
the functions 
  $$
a\mapsto b(t,x,a),\qquad a\mapsto\sigma(t,x,a)
\qquad a\mapsto f (t,x,a) 
 $$
 are continuous on $A$.
\item [(iv)]
 There exist nonnegative constants  $L$ and $r$ such that
\begin{eqnarray} 
|b(t,x,a) - b(t,x',a)| + |\sigma(t,x,a)-\sigma(t,x',a)|
& \leq & L   |x-x'|, \label{lipbsig} \\
|b(t,0,a)| + |\sigma(t,0,a)| & \leq & L,  \label{borbsig}
\\
|f(t,x,a)| + |g(x)| & \leq & L \big(1 + |x|^r \big), \label{PolGrowth_f_g}
\end{eqnarray} 
for all $(t,x,x',a)$ $\in$ $[0,T]\times\R^n \times \R^n \times A$.
\end{itemize}

\vspace{1mm}

\begin{Remark}
{\rm No non-degeneracy assumption on the diffusion coefficient $\sigma$ is imposed.
\qed
}
\end{Remark}

To implement the randomization method  we need   to introduce two additional objects $\lambda, a_0$ satisfying
the following conditions, which are assumed to hold from now on, unless specified otherwise.

\vspace{3mm}

\noindent {\bf (A2)}
\begin{itemize}
\item [(i)]  $\lambda$ is a finite positive measure on  $(A,\calb(A))$
with full topological support.
\item [(ii)] $a_0$ is a fixed, deterministic point in $A$.
\end{itemize}

\vspace{3mm}

We note that such a measure always exists: for instance, it may be constructed as a (possibly infinite) convex combination of Dirac's delta measures located at points which form a dense set in $A$, but in many cases other constructions are easily devised.  
We anticipate that  $\lambda$ will play the role of an intensity measure for a Poisson process and
$a_0$ will be the starting point of some auxiliary random
process introduced later.
Notice that the optimal control problem \eqref{primalvalue} that we are going to introduce in the next paragraph does not depend on
$\lambda,a_0$.
In this sense, {\bf (A2)} is not a restriction 
and we have the choice to fix $a_0\in A$ and
 $\lambda$ satisfying the previous assumption.

\subsection{Formulation of the   control problem}
\label{Primal}

We assume that
$A,b,\sigma,f,g$
 are given and
satisfy the assumptions {\bf (A1)}. We formulate a control problem fixing a setting
$(\Omega,\calf,\P, \F, W)$, where
 $(\Omega,\calf,\P)$ is a complete probability space with
 a right-continuous and $\P$-complete filtration $\F=(\calf_t)_{t\ge 0}$ and
$W$ is an $\R^{d}$-valued
standard Wiener process with respect to $\F$ and $\P$.

Let us denote $\F^W=(\calf^W_t)_{t\ge 0}$ the right-continuous and $\P$-complete filtration generated by $W$.
An admissible control process is any $\F^W$-progressive process $\alpha$ with values in $A$.
The set of admissible control processes is denoted by $\Ac^W$. The controlled equation has the form
\begin{equation}\label{stateeq}
    dX_t^\alpha \ = \ b(t, X_t^\alpha, \alpha_t)\,dt +
\sigma(t, X_t^\alpha, \alpha_t)\,dW_t
\end{equation}
on the interval $[0,T]$ with initial condition $X_0^\alpha=x_0\in\R^n$,
and the reward functional is
\begin{equation}\label{gaineq}
J(\alpha) \ = \ \E\left[\int_0^Tf(t,X_t^\alpha,\alpha_t)\,dt+g(X^\alpha_T)\right].
\end{equation}

Since we assume that {\bf (A1)} holds,  by standard results (see e.g. \cite{rowi} Theorem V. 11.2,
 or   \cite{jacod_book}  Theorem 14.23),
there exists a unique $\F$-adapted strong solution $X^\alpha$ $=$ $(X_t^\alpha)_{0\leq t \leq T}$   to \eqref{stateeq}
with continuous trajectories and such that  
\begin{equation}
    \label{boundlponx}
        \E\,\Big[\sup_{t\in [0,T]}|X_t^\alpha|^p\Big]   \le C\,(1+|x_0|^p) ,
\end{equation}
for every $ p\in [1,\infty)$ and 
for some constant $C$ depending only on $p,T$, and $L$  as defined in {\bf (A1)}. 
The stochastic optimal control problem   consists in maximizing $J(\alpha)$ over all $\alpha\in\Ac^W$:
\begin{equation}\label{primalvalue}
{\text{\Large$\upsilon$}}_0 = \sup_{\alpha\in\Ac^W} J(\alpha).
\end{equation}

\begin{Remark}{\rm 
    In the exposition that follows the filtration $\F$ does not play any special role. If only a complete space $(\Omega,\calf,\P)$ with a Brownian motion $W$ are given, one may take $\F=\F^W$ to obtain a setting as described before.
}
\end{Remark}

\subsection {Formulation of the randomized control problem}
\label{randomizedformulation}

We   assume that
$A,b,\sigma,f,g,\lambda,a_0$
 are given and
satisfy the assumptions {\bf (A1)} and {\bf (A2)}.
We formulate a different optimal control problem, that we call the \emph{randomized}  control problem, associated with the control problem of subsection
\ref{Primal}.
The randomized control problem is formulated fixing a setting
$(\hat \Omega, \hat \calf,\hat \P, \hat W, \hat \mu )$, where
$(\hat \Omega, \hat \calf,\hat \P)$
is an arbitrary complete probability space
with independent random elements  $\hat W$, $\hat \mu$.
The  process
$\hat   W$
is a standard Wiener process in $\R^{d}$ under $\hat \P$.
$\hat \mu$ is a Poisson random measure on  $A$ with intensity $\lambda(da)$
  under $\hat \P$;
thus, $\hat \mu$ is a sum
 of Dirac measures of the form $\hat\mu=\sum_{n\ge 1}\delta_{(\hat S_n,\hat \eta_n)}$,
 where $(\hat \eta_n)_{n\ge 1}$ is an independent sequence of $A$-valued random variables, all having distribution $\lambda(da)/\lambda(A)$,
 and
 $(\hat S_n)_{n\ge 1}$ is a  sequence of random variables
with values in $(0,\infty)$ such that $(\hat S_n-\hat S_{n-1})_{n\ge1}$ ($\hat S_0:=0$) is an independent sequence of exponentially distributed random variables, each with parameter $\lambda(A)$; $(\hat S_n)$,  $(\hat \eta_n)$ and $W$ are all independent.
We also define the $A$-valued, piecewise constant process associated to $\hat\mu$ as follows:
\begin{equation}
\label{I}
\hat I_t \ = \ \sum_{n\ge 0}\hat \eta_n\,1_{[\hat S_n,\hat S_{n+1})}(t), \qquad t\ge 0,
\end{equation}
where we use the convention that $\hat S_0=0$ and $\hat I_0=a_0$, the
  point in  assumption {\bf (A2)}-(ii).
Notice that the formal sum in \eqref{I} makes sense
even if there is no addition operation defined in $A$.

Let $\hat X$ be the solution to the equation
\begin{eqnarray}  \label{dynXrandom}
d\hat X_t &=&  b(t, \hat X_t,\hat I_t)\,dt + \sigma(t,\hat X_t,\hat I_t)\,d\hat W_t,
\end{eqnarray} 
for $t\in [0,T]$, starting from $\hat X_0$ $=$ $x_0\in\R^n$ (the same starting point as in  subsection
\ref{Primal}). We define the filtration
$\F^{\hat W,\hat \mu}=(\calf^{\hat W,\hat \mu}_t)_{t\ge 0}$
setting
\begin{eqnarray} 
\calf^{\hat W,\hat \mu}_t&=&\sigma (\hat W_s,
\hat \mu((0,s]\times C)\,:\, s\in [0,t],\, C\in\calb(A))
\vee \caln,
 \label{expandedfiltration}
\end{eqnarray} 
where $\caln$ denotes the family of $\hat \P$-null sets of $\hat \calf$.
We denote by $\calp(\F^{\hat W,\hat \mu})$ the corresponding predictable
$\sigma$-algebra. We recall that, by
Remark \ref{largfiltrcomp}, the $\F^{\hat W,\hat \mu}$-compensator of $\hat\mu$ is 
$\lambda(da)\,dt$.

Under {\bf (A1)} it is well-known (see e.g. Theorem 14.23 in \cite{jacod_book}) that there exists a unique $\F^{\hat W,\hat \mu}$-adapted strong solution
$\hat X$ $=$ $(\hat X_t)_{0\leq t \leq T}$   to \eqref{dynXrandom}, satisfying
$\hat X_0=\hat x_0$, with continuous trajectories and such that  
\begin{eqnarray} \label{EstimateX}
    \hat \E\,\Big[\sup_{t\in [0,T]}|\hat X_t|^p\Big]  & \le & C\,(1+|x_0|^p),
\end{eqnarray} 
for every $ p\in [2,\infty)$ and 
for some constant $C$ depending only on $p,T$, and $L$  as defined in {\bf (A1)}. 

We can now define the randomized optimal control problem via a change of probability measure of Girsanov type.  
We define the set 
of admissible controls
\begin{align*}
     \hat \calv=\{
     \hat\nu:\Omega\times [0,\infty)\times A\to \R,\;   \Pc(\F^{\hat W,\hat \mu})\otimes \calb(A){\rm  -measurable},\;
0<\hat\nu_t(\omega,a)\le \sup \hat\nu<\infty\}.
\end{align*}
According to Proposition \ref{girsanovpoisson}-(i),
  the   exponential process
\begin{eqnarray}   
\kappa_t^{\hat\nu} \ 
&=& \ \exp\left(\int_0^t\int_A (1 - \hat \nu_s(a))\lambda(da)\,ds
\right)\prod_{0<\hat S_n\le t}\nu_{\hat S_n}(\hat \eta_n),\qquad t\ge 0,\label{doleans}
\end{eqnarray} 
is a martingale with respect to $\hat\P$ and $\F^{\hat W,\hat \mu}$,
and we can define a new probability on $(\hat\Omega,\hat\calf)$ setting
$\hat\P^{\hat\nu}(d\hat\omega)=\kappa_T^{\hat\nu}(\hat\omega)\,\hat\P(d\hat\omega)$. \footnote{In  Proposition \ref{girsanovpoisson} this was denoted $\hat\P^{\hat\nu}_T$. Here we use a shorter notation.}
From
  Proposition \ref{girsanovpoisson}-(ii), it follows that under $\hat \P^{\hat\nu}$
the $\F^{\hat W,\hat \mu}$-compensator of $\hat\mu$ on the set
$[0,T]\times A$ is the random measure $\hat\nu_t(a)\lambda(da)dt$. 
We finally introduce the gain functional of the randomized control problem
\begin{eqnarray}  \label{defJrandomized}
J^\Rc(\hat\nu) &=&  \hat\E^{\hat\nu}
\left[\int_0^Tf(t,\hat X_t,\hat I_t)\,dt+g(\hat X_T)\right],
\end{eqnarray} 
where $\hat\E^{\hat\nu}$ is the expectation under $\hat\P^{\hat\nu}$.
The randomized stochastic optimal control problem consists in maximizing
$ J^\Rc(\hat\nu)$ over all $\hat\nu\in\hat\calv$. Its value
is defined as
\begin{equation}\label{dualvalue}
{\text{\Large$\upsilon$}}_0^\Rc \;=\;   \sup_{\hat \nu\in\hat \calv} J^\Rc(\hat \nu).
\end{equation}

To show that $ J^\Rc(\hat\nu)$ and ${\text{\Large$\upsilon$}}_0^\Rc$ are well defined we need 
the following lemma.

\begin{Lemma}\label{Wrestawiener}
With the previous notations and assumptions, $\hat W$ remains an $\F^{\hat W,\hat\mu}$-Brownian motion under $\hat\P^{\hat\nu}$.     
\end{Lemma}

\noindent
{\bf Proof.}  Since $\kappa_T^{\hat\nu}>0$, the probabilities
$\hat\P$ and $\hat\P^{\hat\nu}$ are equivalent, so the quadratic variation of $\hat W$ computed under $\hat\P$ and $\hat\P^{\hat\nu}$
is the same, and equals $\langle \hat W\rangle_t=t$. By Lévy's theorem, see e.g. 
\cite{CohenElliottbook} Theorem 14.4.1,  it is enough to show
that $\hat W$ is a $(\hat\P^{\hat\nu},\F^{\hat W,\hat\mu})$-local martingale, which is equivalent
to the fact that $\kappa^{\hat\nu}\hat W$ is a $(\hat\P,\F^{\hat W,\hat\mu})$-local martingale. In turn, this follows from a general fact (see see e.g. 
\cite{CohenElliottbook} Theorem 10.2.21): since $\kappa^{\hat\nu}$ is a $(\hat\P,\F^{\hat W,\hat\mu})$-martingale
of finite variation, it is purely discontinuous  and
therefore   orthogonal (under $\hat\P$)
to $\hat W$; thus, their product $\kappa^{\hat\nu}\hat W$ is a $(\hat\P,\F^{\hat W,\hat\mu})$-local martingale.
\qed

By this Lemma,   the same passages that led to \eqref{EstimateX} now give the following more general estimate: 
\begin{eqnarray} \label{EstimateX_nu}
\sup_{\hat\nu\in\Vc}\,\hat\E^{\hat\nu}\,\Big[\sup_{t\in [0,T]}|\hat X_t|^p\Big]  &\le & C\,(1+|x_0|^p),
\end{eqnarray} 
where $\hat\E^{\hat\nu}$ denotes the expectation with respect to $\hat\P^{\hat\nu}$ and $C$ is the same constant as in \eqref{EstimateX}.

\begin{Remark}
\emph{
Later on we need a conditional form of the inequality \eqref{EstimateX_nu}. Considering the stochastic equation on an interval $[t,T]\subset [0,T]$ with initial condition $X_t$, the same estimates that led to 
\eqref{EstimateX_nu} now give
\begin{eqnarray} \label{EstimateX_nu_cond}
\esssup_{\hat\nu\in\Vc}\,\hat\E^{\hat\nu}\,\Big[\sup_{s\in [t,T]}|\hat X_s|^p\,\Big|\, \calf_t^{\hat W,\hat \mu}\Big]  &\le & C\,(1+|\hat X_t|^p),
\end{eqnarray} 
for every $ p\in [2,\infty)$ and 
for   $C$ depending only on $p,T$, and $L$  as defined in {\bf (A1)}. 
}
\end{Remark} 

\begin{Remark}\label{infzero}\emph{
Let us define
$\hat\Vc_{\inf\,>\,0}=\{\hat\nu\in \hat\Vc\,:\, \inf_{\hat \Omega\times[0,T]\times A}\hat\nu>0\}$.
Then
\begin{equation}\label{eqinfzero}
{\text{\Large$\upsilon$}}_0^\Rc \;= \;
\sup_{\hat\nu\in\hat\Vc_{\inf\,>\,0}}J^\Rc(\hat\nu).
\end{equation}
Indeed, given $\hat\nu\in\hat\Vc$ and $\epsilon>0$, define $\hat\nu^\epsilon=
\hat\nu\vee \epsilon\in \hat\Vc_{\inf\,>\,0}$
and write the gain \eqref{defJrandomized} in the form
 $$
J^\Rc(\hat\nu^\epsilon) =  \hat\E
\left[\kappa_T^{\hat\nu^\epsilon}\left( \int_0^Tf(t,\hat X_t,\hat I_t)\,dt+g(\hat X_T)\right)\right].
 $$
 It is easy to see that $J^\Rc(\hat\nu^\epsilon) \to J^\Rc(\hat\nu) $ as $\epsilon\to 0$,
 which implies
$${\text{\Large$\upsilon$}}_0^\Rc=\sup_{\hat\nu\in\hat\Vc}J^\Rc(\hat\nu)\le \sup_{\hat\nu\in\hat\Vc_{\inf\,>\,0}}J^\Rc(\hat\nu).
 $$
 The other inequality being obvious, we obtain  \eqref{eqinfzero}.
}
\qed
\end{Remark}

\vspace{2mm}

\begin{Remark} \label{remchoice}
{\rm
We end this section noting that a randomized control problem can be constructed starting from the initial control problem
with partial observation. Indeed, let $(\Omega,\calf,\P, \F, W)$
be the setting for the  stochastic optimal control problem formulated in subsection \ref{Primal}.
Suppose that  $(\Omega',\calf',\P')$ is another  probability space where a  Poisson random
 measure $\mu$ with intensity $\lambda$ is defined (for instance by a classical result,
 see \cite{Zabczyk96} Theorem 2.3.1, we may take $\Omega'=[0,1]$,
 $\calf'$ the corresponding Borel sets and $\P'$ the Lebesgue measure).
Then we define $\bar \Omega=\Omega\times \Omega'$, we denote by $\bar \calf$ the completion of
$\calf\otimes \calf'$ with respect to $\P\otimes \P'$ and by
$\bar \P$ the extension of $\P\otimes \P'$ to $\bar \calf$.
The Brownian motion $W$ in $\Omega$ and the random measure
$\mu$ in $\Omega'$ have obvious extensions to $\bar \Omega$, that will be denoted by
the same symbols. Clearly, $(\bar \Omega, \bar \calf,\bar \P,     W,   \mu )$
is a setting for a  randomized control problem as formulated before, that we call
\emph{product extension} of the  setting $(\Omega,\Fc,\P, W )$ for the initial control problem \eqref{primalvalue}.

We note that the initial formulation of a randomized
setting $(\hat \Omega, \hat \calf,\hat \P,  \hat W, \hat \mu)$ was more general, since  it was not required that
$\hat \Omega$ should be a product space $\Omega\times \Omega'$ and, even if it were the case, it was not required that the process $\hat W$
should depend only on $\omega\in\Omega$ while the random measure $\hat\mu$ should depend only on  $\omega'\in\Omega'$.

The aim of the following paragraph (subsection \ref{Sub-valuerandpro}) is to prove that the value ${\text{\Large$\upsilon$}}_0^\Rc$
does
not depend on the precise construction of the \emph{randomized setting}. In
particular, it justifies in the proofs of section 
\ref{proofthm}
 to consider the product
extension (see subsection \ref{Sub-proofineq} for instance).
}
\qed
\end{Remark}

\subsection {The value of the randomized control problem}\label{Sub-valuerandpro}

In this section it is our purpose to show that the value ${\text{\Large$\upsilon$}}_0^\Rc$
of the randomized control problem
defined in \eqref{dualvalue}
does not depend on the specific
setting $(\hat \Omega, \hat \calf,\hat \P,  \hat W, \hat \mu  )$,
so that
 it is just a functional of the (deterministic) elements
$A,b,\sigma,f,g, \lambda, a_0$ and the starting point $x_0$.
 Later on, in
 Theorem
\ref{MainThm},  we will prove
that in fact ${\text{\Large$\upsilon$}}_0^\Rc$
does not depend on the choice of $\lambda$ and $a_0$ either.

So let now $(\tilde \Omega, \tilde \calf,\tilde \P,  \tilde W, \tilde \mu )$
be another setting for the randomized control problem,
as in Section \ref{randomizedformulation}, and let $\F^{\tilde W,\tilde \mu}$,
    $\tilde X$, $\tilde I$, $\tilde \calv$ be defined in analogy
 with what was done before. So, for any   $\tilde\nu\in \tilde \calv$, we also define $\kappa^{\tilde\nu}$ and the probability
$d\tilde\P^{\tilde\nu} =\kappa_T^{\tilde\nu} \,d\tilde\P$ as well as the gain and the value
\begin{eqnarray*}
\tilde J^\Rc(\tilde\nu) \; = \;   \tilde\E^{\tilde\nu}
\left[\int_0^Tf(t,\tilde X_t,\tilde I_t)\,dt+g(\tilde X_T)\right],
\qquad
\tilde{\text{\Large$\upsilon$}}_0^\Rc \; = \;    \sup_{\tilde \nu\in\tilde \calv} \tilde J^\Rc(\tilde\nu).
\end{eqnarray*}
We recall that the  gain functional and value for the setting
$(\hat \Omega, \hat \calf,\hat \P,  \hat W, \hat \mu )$
was defined in \eqref{defJrandomized} and \eqref{dualvalue} and
denoted   by $J^\Rc$ and ${\text{\Large$\upsilon$}}_0^\Rc$ rather than $\hat J^\Rc$ and $\hat{\text{\Large$\upsilon$}}_0^\Rc$, to simplify
the notation.

\begin{Proposition} \label{indepofthesetting}
With the previous notation, we have ${\text{\Large$\upsilon$}}_0^\Rc$ $=$ $\tilde{\text{\Large$\upsilon$}}_0^\Rc$.
In other words, ${\text{\Large$\upsilon$}}_0^\Rc$ only depends on the objects $A,b,\sigma,f,g, \lambda, a_0$ appearing
in the assumptions {\bf (A1)} and {\bf (A2)}, as well as on the starting point $x_0$.
\end{Proposition}

 \noindent {\bf Proof.} It is enough to prove that ${\text{\Large$\upsilon$}}_0^\Rc\le \tilde{\text{\Large$\upsilon$}}_0^\Rc$,
 since the opposite inequality is established by the same arguments. Writing the gain $J^\Rc(\hat\nu)$ defined in
 \eqref{defJrandomized}
 in the form
 $$
J^\Rc(\hat\nu) =  \hat\E
\left[\kappa_T^{\hat\nu}\left( \int_0^Tf(t,\hat X_t,\hat I_t)\,dt+g(\hat X_T)\right)\right],
 $$
 recalling the definition \eqref{doleans} of the process $\kappa^{\hat\nu}$
 and noting that the
process $\hat I$ is completely determined by $\hat\mu$,
we see that  $J^\Rc(\hat\nu)$ only depends on the
(joint) law of $(\hat X,\hat \mu,\hat \nu)$ under $\hat\P$.
Since, however,  $\hat X$ is the solution to equation \eqref{dynXrandom} with initial
condition  $\hat X_0$ $=$ $ x_0$,
 it is easy to check that under  our assumptions
the  law of $(\hat X,\hat \mu,\hat \nu)$  only depends on the law of
$(\hat W,\hat \mu, \hat\nu)$ under $\hat\P$.  Similarly, $\tilde J^\Rc(\tilde\nu)$ only depends on the law of $(\tilde W,\tilde \mu,\tilde \nu)$
 under $\tilde\P$.

 Next we claim that, given $\hat\nu\in\hat \calv$ there exists
 $\tilde\nu\in\tilde \calv$ such that the law of
 $(\hat W,\hat \mu,\hat \nu)$ under $\hat \P$ is the same as
 the law of
 $(\tilde W,\tilde\mu,\tilde \nu)$ under $\tilde\P$.
 Assuming the claim for a moment, it follows from the previous
 discussion that for this choice of $\tilde\nu$ we have
 \begin{eqnarray*}
 J^\Rc(\hat\nu) \; = \; \tilde J^\Rc(\tilde\nu) &\le&  \tilde{\text{\Large$\upsilon$}}_0^\Rc,
 \end{eqnarray*}
 and taking the supremum over $\hat\nu\in\hat \calv$ we deduce that ${\text{\Large$\upsilon$}}_0^\Rc$ $\le$ $\tilde{\text{\Large$\upsilon$}}_0^\Rc$, which proves the result.

It only remains to prove the claim. By a monotone class argument
we may suppose that $\hat\nu_t(a)=k(a)\,\phi_t\,\psi_t$,
where $k$ is a $\calb(A)$-measurable, $\phi$ is $\F^{\hat W}$-predictable
and $\psi$ is $\F^{\hat \mu}$-predictable
(where these filtrations are the ones generated by $\hat W$ and $\hat \mu$ respectively).
We may further suppose that $\phi_t=1_{(t_0,t_1]}(t)\phi_0(\hat W_{s_1},\ldots,\hat W_{s_h})$,
for an integer $h$ and deterministic times $0\le s_1\le\ldots s_h\le t_0 <t_1$
and a Borel function $\phi_0$ on $\R^h$,
since this class of processes generates the predictable $\sigma$-algebra
of $\F^{\hat W}$, and that $\psi_t=1_{(\hat S_n,\hat S_{n+1}]}(t)\psi_0(\hat S_1,\ldots, \hat S_n, \hat \eta_1,\ldots, \hat \eta_n,t)$,
for an integer $n\ge 1$ and
 a Borel function $\psi_0$ on $\R^{2n+1}$,
since this class of processes generates the predictable $\sigma$-algebra
of $\F^{\hat \mu}$ (see \cite{ja}, Lemma (3.3)). It is immediate to verify that
the required process $\tilde\nu$ can be defined setting
$$
\tilde\nu_t(a)=k(a)\,
1_{(t_0,t_1]}(t)\phi_0(\tilde W_{s_1},\ldots, \tilde W_{s_h}) \,
1_{(\tilde S_n,\tilde S_{n+1}]}(t)\psi_0(\tilde S_1,\ldots,
\tilde S_n, \tilde \eta_1,\ldots, \tilde \eta_n,t),
$$
where $(\tilde S_n,\tilde \eta_n)_{n\ge 1}$
are associated to the measure $\tilde \mu$, i.e.
$\tilde\mu=\sum_{n\ge 1}\delta_{(\tilde S_n,\tilde \eta_n)}$.
 \qed

\subsection{Equivalence of the original and the randomized control pro\-blem}
We can now state one of the main results.

\begin{Theorem}
\label{MainThm}
Assume that {\bf (A1)} and {\bf (A2)} are satisfied. 
Then the values of the partially observed control problem
and of the randomized control problem are equal:
\begin{eqnarray}  \label{equivrandom}
{\text{\Large$\upsilon$}}_0 &=& {\text{\Large$\upsilon$}}_0^\Rc,
\end{eqnarray} 
where ${\text{\Large$\upsilon$}}_0$ and ${\text{\Large$\upsilon$}}_0^\Rc$ are defined by  \eqref{primalvalue} and \eqref{dualvalue} respectively.
This common value only depends on   the  starting point  $x_0$ and on the objects $A,b,\sigma,f,g$ appearing in assumption {\bf (A1)}.
\end{Theorem}

The last sentence follows immediately from
Proposition
\ref{indepofthesetting}, from the equality ${\text{\Large$\upsilon$}}_0$ $=$  ${\text{\Large$\upsilon$}}_0^\Rc$
and from the obvious fact that ${\text{\Large$\upsilon$}}_0$ cannot depend on $\lambda, a_0$ introduced in assumption {\bf (A2)}.
The proof of the equality \eqref{equivrandom} is contained in the next section.

\begin{Remark}
    {\rm 
The requirement that $\lambda$ has full support will not be used in the proof of the inequality ${\text{\Large$\upsilon$}}_0^\Rc$ $\le$
${\text{\Large$\upsilon$}}_0$. This inequality holds just assuming that $\lambda(da)$ is an arbitrary finite Borel measure on $A$. }
\qed
\end{Remark}

\section{Proof of Theorem \ref{MainThm}}
\label{proofthm}

The proof  is split into two parts, corresponding to the inequalities ${\text{\Large$\upsilon$}}_0^\Rc$ $\le$ ${\text{\Large$\upsilon$}}_0$ and ${\text{\Large$\upsilon$}}_0$ $\le$ ${\text{\Large$\upsilon$}}_0^\Rc$.
In the sequel, {\bf (A1)} and {\bf (A2)} are always assumed to hold. 

Before starting with the rigorous proof, let us have a look at the main points.
\begin{itemize}
\item ${\text{\Large$\upsilon$}}_0^\Rc$ $\le$ ${\text{\Large$\upsilon$}}_0$. First, in  Lemma \ref{L:AWmu}, we prove that the value of the primal problem ${\text{\Large$\upsilon$}}_0$ does not change if we reformulate it on the enlarged probability space where the randomized problem lives, taking the supremum over $\Ac^{W,\mu'}$, which is the set of controls $\bar\alpha$ progressively measurable with respect to a larger filtration, denoted $\F^{W,\mu'_\infty}$ and associated to  $W$ and the Poisson random measure $\mu'$. Second, in Proposition \ref{P:Ineq_I}  (based on the key Lemma \ref{L:Replication_I}), we prove that for every $\nu\in\Vc_{\inf>0}$ there exists $\bar\alpha^\nu\in\Ac^{W,\mu'}$ such that $\mathscr L_{\P^\nu}(W,I)=\mathscr L_{\bar\P}(W,\bar\alpha^\nu)$ and consequently $J^\Rc(\nu)=  J(\bar\alpha^\nu)$. It follows that
\[
{\text{\Large$\upsilon$}}_0^\Rc \ = \ \sup_{\nu\in\calv_{\inf>0}}J^\Rc(\nu) \ \le \ \sup_{\bar\alpha\in\Ac^{W,\mu'} }J(\bar\alpha^\nu).
\]
Since ${\text{\Large$\upsilon$}}_0$ $=$ $\sup_{\bar\alpha\in\Ac^{W,\mu'}} \bar J(\bar\alpha)$ by Lemma \ref{L:AWmu},   we   obtain the inequality ${\text{\Large$\upsilon$}}_0^\Rc$ $\le$ ${\text{\Large$\upsilon$}}_0$.
%Finally, we notice that the set $\Vc_{\inf>0}$ (rather than $\Vc$) is required only for the proof of this inequality (${\text{\Large$\upsilon$}}_0^\Rc$ $\le$ ${\text{\Large$\upsilon$}}_0$), and in particular for the proof of Lemma
%\ref{L:Replication_I} (see for instance after formula \eqref{vartheta_1}).

\item ${\text{\Large$\upsilon$}}_0$ $\le$ ${\text{\Large$\upsilon$}}_0^\Rc$. The proof of this inequality is based on a ``density'' result (which corresponds to the key Proposition \ref{extensionapproximation}). Roughly speaking, we prove that the class $\{\bar\alpha^\nu\colon\nu\in\Vc_{\inf>0}\}$ is dense in $\Ac^{W,\mu'}$, with respect to an appropriate metric $\tilde\rho$ in the space of control processes   (the same metric used in Lemma 3.2.6 in \cite{80Krylov}). Then, the inequality ${\text{\Large$\upsilon$}}_0$ $\le$ ${\text{\Large$\upsilon$}}_0^\Rc$ follows from the stability Lemma \ref{contrhotilde}, which states that, under Assumption {\bf (A1)}, the gain functional is continuous with respect to the metric $\tilde\rho$.
\end{itemize}

\subsection{Proof of the inequality ${\text{\Large$\upsilon$}}_0^\Rc$ $\le$ ${\text{\Large$\upsilon$}}_0$}
\label{Sub-proofineq}

The requirement that $\lambda$ has full support will not be used in the proof of the inequality ${\text{\Large$\upsilon$}}_0^\Rc$ $\le$
${\text{\Large$\upsilon$}}_0$ and will not be assumed to hold in this paragraph.

 Let $(\Omega,\calf,\P, \F, W)$
be a setting for the stochastic optimal control problem  formulated in Subsection \ref{Primal}. We construct a setting for a randomized control problem in the form of a product
 extension as described in Remark \ref{remchoice}.

Let $\lambda$ be a finite Borel measure on $A$. 
As in  Lemma \ref{lambdaliftato}, we construct a finite measure $\lambda' (dr)$ on $(\R,\calb(\R))$ and a Borel map $\pi:\R\to A$ such that $\lambda=\lambda' \circ\pi^{-1}$  and
    $\lambda' (\{r\})=0$ for every $r\in \R$ 
Now  let  $(\Omega',\calf',\P')$ denote the canonical probability space of a
 non-explosive Poisson point process on $(0,\infty)\times \R$ with intensity $\lambda'(dr)$.
 Thus, $\Omega'$ is the set of sequences $\omega'=(t_n,r_n)_{n\geq1}\subset(0,\infty)\times \R$ with $t_n<t_{n+1}\to\infty$,
 $(T_n,R_n)_{n\geq1}$ is the canonical marked point process  (i.e.  $T_n(\omega')=t_n$, $R_n(\omega')=r_n$),  and
 $\mu'$ $=$ $\sum_{n\ge 1}\delta_{(T_n,R_n)}$ is the corresponding random measure.
 Let $\Fc'$ denote the   smallest $\sigma$-algebra such that all the maps $T_n,R_n$ are measurable, and   $\P'$  the unique probability on
$\Fc'$ such that  $\mu'$ is a Poisson random measure with intensity $\lambda'$.

Define 
\[A_n=\pi(R_n)
\]
and consider the marked point process $(T_n,A_n)_{n\ge1}$ in $A$.
Recalling that $\lambda=\lambda' \circ\pi^{-1}$ it is immediate to see that $(T_n,A_n)_{n\ge1}$ is a Poisson process in $A$ with intensity $\lambda(da) $,  defined in $(\Omega',\calf',\P')$. We   consider
the  random measures 
\[
\mu=\sum_{n\ge1}\delta_{(T_n,A_n)},
\qquad
\mu'=\sum_{n\ge1}\delta_{(T_n,R_n)}
\]
on the Borel subsets of $(0,\infty)\times A$
and  $(0,\infty)\times \R$ respectively, as well as  the corresponding natural filtrations
$\F^\mu=(\calf_t^\mu)_{t\ge0}$ and
$\F^{\mu'}=(\calf_t^{\mu'})_{t\ge0}$, 
where
\[
\calf_t^\mu=\sigma\Big( \mu((0,s]\times C)\,:\, s\in [0,t],C\in\calb(A)\Big),
\quad
\calf_t^{\mu'}=\sigma\Big( \mu((0,s]\times B)\,:\, s\in [0,t],B\in\calb(\R)\Big).
\]
Note that $
\calf_t^\mu\subset 
\calf_t^{\mu'}$ and $
\calf_\infty^\mu\subset 
\calf_\infty^{\mu'}=\calf'$, the latter equality due to the canonical form of the space $\Omega'$.

Then we define $\bar \Omega=\Omega\times \Omega'$, we denote by $\bar \calf$ the completion of
$\calf\otimes \calf'$ with respect to
$\P\otimes \P'$ and by
$\bar \P$ the extension of $\P\otimes \P'$ to $\bar \calf$.
The Brownian motion $ W$ in $\Omega$ and the random measures
$\mu,\mu'$ in $\Omega'$ have obvious extensions to $\bar \Omega$, that will be denoted by
the same symbols. Then $(\bar \Omega, \bar \calf,\bar \P,     W,   \mu)$
is a setting for a  randomized control problem as formulated  in section \ref{randomizedformulation}. Recall that $\F^{ W}$ denotes the $\P$-completed filtration in $( \Omega,   \calf)$ generated by the Wiener
process $W$. All filtrations $\F^{ W}$,  $\F^\mu$,   $\F^{\mu'}$  can also be lifted to  filtrations in $(\bar \Omega, \bar \calf)$, and $\bar\P$-completed.
In the sequel it should be clear from the context whether they are considered as filtrations in $(\bar \Omega, \bar \calf)$ or in their original spaces.
We will compare the value of the randomized optimal control problem in this setting with the value of the original control problem.

\vspace{2mm}

Recall that in the original optimal control problem the class $\cala^W$ of admissible controls consists of all $\F^W$-progressive processes. As a preliminary step  to proving the inequality ${\text{\Large$\upsilon$}}_0^\Rc$ $\le$ ${\text{\Large$\upsilon$}}_0$,  in Lemma \ref{L:AWmu}  we  formulate 
a new equivalent control problem (i.e., with the same value) on the enlarged space
$(\bar\Omega,\bar\Fc)$. We show that the value ${\text{\Large$\upsilon$}}_0$ remains the same even if 
we allow the controls to be progressive with respect to 
 another, larger  filtration
$\F^{W,\mu'_\infty}=(\calf^{W,\mu'_\infty}_t)_{t\ge 0}$
in $(\bar \Omega, \bar \calf)$ defined setting
$$
\calf^{W,\mu'_\infty}_t \ = \ \Fc_t^W\vee\Fc'\vee\caln
$$
where $\caln$ denotes the family of $\bar \P$-null sets of $\bar \calf$ and $\Fc'$ denotes the $\sigma$-algebra extended to $(\bar \Omega, \bar \calf)$, namely
$\{\Omega\times B\, :\, B\in \Fc'\}$.
 This fact has an intuitive meaning: for an $\F^{W,\mu'_\infty}$-progressive process, at each time $t$ the control action is chosen based on knowledge of $\calf^W_t$ and $\calf'$: but $\calf'$ is generated by an independent Poisson process and its knowledge should be irrelevant for optimization. 
 
\begin{Lemma}\label{L:AWmu}
We have ${\text{\Large$\upsilon$}}_0$ $=$ $\sup_{\bar\alpha\in\Ac^{W,\mu'}} \bar J(\bar\alpha)$, where
\[
\bar J(\bar\alpha) \ = \ \bar\E\left[\int_0^Tf(t,X_t^{\bar\alpha},\bar\alpha_t)\,dt+g(X^{\bar\alpha}_T)\right],
\]
and $\Ac^{W,\mu'}$ is the set of all $\F^{W,\mu'_\infty}$-progressive processes $\bar\alpha$ with values in $A$, defined on $(\bar\Omega,\bar\calf)$.  $X^{\bar\alpha}$ $=$ $(X_t^{\bar\alpha})_{  t\in[0,  T]}$ is the  solution to \eqref{stateeq} $($with $\bar\alpha$ in place of $\alpha$$)$
satisfying $X_0^{\bar\alpha}=x_0$, which exists and is unique in the class of continuous processes adapted to the filtration
$\F^{W,\mu'_\infty}$.
\end{Lemma}

\noindent
\textbf{Proof.}
The   inequality ${\text{\Large$\upsilon$}}_0$ $\le$ $\sup_{\bar\alpha\in\Ac^{W,\mu'}} \bar J(\bar\alpha)$ is immediate, since every control
$\alpha\in\Ac^W$ also lies in $\Ac^{W,\mu'}$ and $J(\alpha)=\bar J(\alpha)$, whence $J(\alpha)\leq\sup_{\bar\alpha\in\Ac^{W,\mu'}}
\bar J(\bar\alpha)$ and so ${\text{\Large$\upsilon$}}_0$ $=$  $\sup_{\alpha\in\Ac^W} J(\alpha)\leq\sup_{\bar\alpha\in\Ac^{W,\mu'}} \bar J(\bar\alpha)$.

Let us prove the opposite inequality. Fix $\tilde\alpha\in\Ac^{W,\mu'}$ and consider the (uncompleted) filtration
$\F'':=( \Fc_t^W\vee\Fc')_{t\ge0}$. Then we can find  an  $A$-valued $\F''$-progressive process $\bar   \alpha$
such that  $\bar\alpha=\tilde\alpha$   $\bar\P(d\omega\,d\omega')dt$-almost surely, so that  in particular $\bar J(\bar\alpha)=\bar J(\tilde\alpha)$.
It is easy to verify that, for every $\omega'\in\Omega'$,  the process $\alpha^{\omega'}$, defined by
$\alpha_t^{\omega'}(\omega):=\bar\alpha_t(\omega,\omega')$, is $\F^W$-progressive. Consider now the controlled equation on $[0,T]$
\begin{align}
\label{SDEomega'}
X_t \ &= \ x_0 + \int_0^t b(s,X_s, \alpha_s^{\omega'})\,ds + \int_0^t \sigma(s,X_s, \alpha_s^{\omega'})\,dW_s \\
&= \ x_0 + \int_0^t b(s,X_s, \bar\alpha_s(\cdot,\omega'))\,ds + \int_0^t \sigma(s,X_s, \bar\alpha_s(\cdot,\omega'))\,dW_s. \notag
\end{align}
From the first line of \eqref{SDEomega'} we see that, under Assumption {\bf (A1)}, for every $\omega'$ there exists a unique (up to indistinguishability) continuous  process $X^{\alpha^{\omega'}}=(X_t^{\alpha^{\omega'}})_{t\in [0,T]}$   solution to \eqref{SDEomega'}, adapted to the filtration
$\F^W$.
 On the other hand, from the second line of \eqref{SDEomega'}, it follows that the process $X^{\bar\alpha}(\cdot,\omega')=(X_t^{\bar\alpha}(\cdot,\omega'))_{t\in [0,T]}$ solves the above equation. From the pathwise uniqueness of strong solutions to equation \eqref{SDEomega'}, it follows that $X_t^{\alpha^{\omega'}}(\omega)=X_t^{\bar\alpha}(\omega,\omega')$, for all $t\in[0,T]$, $\P(d\omega)$-a.s. By the Fubini theorem
\[
\bar J(\tilde\alpha) =
\bar J(\bar\alpha)  = \int_{\Omega'} \E\bigg[\int_0^Tf(t,X_t^{\alpha^{\omega'}},\alpha_t^{\omega'})\,dt+g(X_T^{\alpha^{\omega'}})
\bigg] \, \P'(d\omega').
\]
Since the inner expectation equals the gain $J(\alpha^{\omega'})$,
it cannot exceed ${\text{\Large$\upsilon$}}_0$ and it follows that
$\bar J(\tilde\alpha)\le {\text{\Large$\upsilon$}}_0$.  The claim follows from the  arbitrariness of $\tilde\alpha$.
\qed

\vspace{2mm}

Next we formulate the randomized control problem in the setting 
  $(\bar \Omega, \bar \calf,\bar \P,     W,   \mu)$
that we have constructed. The formulation is the same as in  Subsection \ref{randomizedformulation}, but with a different notation. Thus, as in \eqref{I}  
we define the $A$-valued process
\[
 I_t \ = \ \sum_{n\ge 0}A_n\,1_{[ T_n, T_{n+1})}(t), \qquad t\ge 0,
 \]
where we use the convention that $ T_0=0$ and $ I_0=a_0$ (the
  point in  assumption {\bf (A2)}-(ii)) and 
    the solution $  X$ to the equation
\eqref{dynXrandom}, namely
\begin{eqnarray}  \nonumber
d X_t &=&  b(t,  X_t, I_t)\,dt + \sigma(t, X_t, I_t)\,d W_t,
\quad t\in [0,T],
\qquad 
X_0=x_0.
\end{eqnarray} 
We   define the filtration
$\F^{ W,  \mu}=(\calf^{ W,  \mu}_t)_{t\ge 0}$ in $(\bar \Omega, \bar \calf)$  by
$$
\calf^{  W,  \mu}_t  =
\Fc_t^W\vee\calf^\mu_t\vee\caln,
$$
 ($\caln$ denotes the family of $\bar \P$-null sets of $\bar \calf$),   we   introduce the classes $  \calv, \Vc_{\inf\,>\,0}$ and, for any
 admissible control $ \nu \in \calv$, the   process
\begin{eqnarray}  \nonumber
\kappa_t^{\nu} (\omega,\omega')\ 
&=& \ \exp\left(\int_0^t\int_A (1 -  \nu_s(\omega,\omega',a))\lambda(da)\,ds
\right)\prod_{0< T_n(\omega')\le t}\nu_{ T_n(\omega')}(\omega,\omega',A_n),\qquad t\ge 0, 
\end{eqnarray} 
which 
is a martingale with respect to $\bar\P$ and $\F^{ W, \mu}$,
 the probability
 $\P^{ \nu}(d\omega\, d\omega' )=\kappa_T^{\nu}(\omega,\omega')\,\bar\P(d \omega\,d\omega')$ on $(\bar\Omega,\bar\calf)$, and finally the gain  
\begin{eqnarray}  
J^\Rc(\nu) &=&  \E^{\nu}
\left[\int_0^Tf(t, X_t, I_t)\,dt+g( X_T)\right],
\end{eqnarray} 
as in \eqref{defJrandomized}.

\vspace{3mm}

The next result provides a decomposition of any element $\nu$ $\in$ $\Vc$; in fact, the result holds for any  $\calp(\F^{W,\mu})\otimes \calb(A)$-measurable random field.

\begin{Lemma}\label{L:nu_pred}
 Let $\nu$ be a $\calp(\F^{W,\mu})\otimes \calb(A)$-measurable random field.
 \newline
\textup{(i)} There exists a $\bar\P$-null set $\bar N\in\Nc$ such that $\nu$ admits the following representation
\begin{align*}
\nu_t(\omega,\omega',a) \ &= \ \nu_t^{(0)}\big(\omega,a\big) \, 1_{\{0\le t\leq T_1(\omega')\}} \\
&\quad \ + \sum_{n=1}^\infty \nu_t^{(n)}\big(\omega,(T_1(\omega'),A_1(\omega')),\ldots,(T_n(\omega'),A_n(\omega')),a\big) \, 1_{\{T_n(\omega')<t\leq T_{n+1}(\omega')\}}, \notag
\end{align*}
for all $(\omega,\omega',t,a)\in\bar\Omega\times\R_+\times A$, $(\omega,\omega')\notin\bar N$, for some maps $\nu^{(n)}\colon\Omega\times\R_+\times(\R_+\times A)^n\times A\rightarrow(0,\infty)$, $n\geq1$, $($resp. $\nu^{(0)}\colon\Omega\times\R_+\times A\rightarrow(0,\infty)$$)$, which are $\Pc(\F^W)\otimes\Bc((\R_+\times A)^n)\otimes\Bc(A)$-measurable $($resp. $\Pc(\F^W)\otimes\Bc(A)$-measurable$)$. 

\vspace{2mm}

\noindent\textup{(ii)} For every $\omega\in\Omega$ define a map $\nu^\omega=\nu_t^\omega(\omega',a): \Omega'\times \R_+\times A\to \R$ by
\begin{align}\label{representation_nu^omega}
\nu_t^\omega(\omega',a) \ &= \ \nu_t^{(0)}\big(\omega,a\big) \, 1_{\{0\le t\leq T_1(\omega')\}} \\
&\quad \ + \sum_{n=1}^\infty \nu_t^{(n)}\big(\omega,(T_1(\omega'),A_1(\omega')),\ldots,(T_n(\omega'),A_n(\omega')),a\big) \, 1_{\{T_n(\omega')<t\leq T_{n+1}(\omega')\}}, \notag
\end{align}
where, clearly, $\nu_\cdot^{(n)}(\omega,\cdot)$ $($resp. $\nu_\cdot^{(0)}(\omega,\cdot)$$)$ is $\Bc(\R_+)\otimes\Bc((\R_+\times A)^n)\otimes\Bc(A)$-measurable $($resp. $\Bc(\R_+)\otimes\Bc(A)$-measurable$)$.

Denote  $\caln'$ the family of $\P'$-null sets of $\Fc'$. 
Then there exists $\tilde N\in\Fc$, with $\P(\tilde N)=0$, such that
 for every $\omega\notin\tilde N$ there exists $N'_\omega\in\Nc'$ such that 
\[
\nu_t^\omega(\omega',a) \ = \ \nu_t(\omega,\omega',a), 
\] 
for all $(\omega',t,a)\in\Omega'\times\R_+\times A$, $\omega'\notin N_\omega'$.

\end{Lemma}

\noindent {\bf Proof.} Consider the (uncompleted) filtration $(\Fc_t^W\vee\calf^\mu_t)_{t\ge0}$.  Up to a $\bar\P$-null set $\bar N$, the random field is measurable with respect to 
$\calp
((\Fc_t^W\vee\calf^\mu_t)_{t\ge0})\otimes \calb(A)$. 
Recalling the form of $\calp(\F^\mu)\otimes \calb(A)$-measurable random fields given in Proposition \ref{propnaturalfiltration}, point (i) can be proved by a 
monotone class argument. Point (ii) follows immediately.
\qed

\vspace{3mm}
 
Given $\nu\in\Vc_{\inf\,>\,0}$,  the maps $\nu^{(n)}$ in 
 \eqref{representation_nu^omega}    can be chosen to satisfy $0<\inf\nu\le \nu^{(n)}\leq \sup\nu<\infty$, so the process  $\nu^\omega$ introduced in the previous Lemma 
  also satisfies 
$0<\inf\nu\le \nu^\omega\leq \sup\nu<\infty$.

For any $\omega\in\Omega$ define a process $\kappa^{\nu^\omega}:\Omega'\times [0,\infty)\to\R$, setting
\begin{eqnarray*} 
\kappa_t^{\nu^\omega}(\omega') \ 
= \ \exp\left(\int_0^t\int_A (1 - \nu_s^\omega(\omega',a))\lambda(da)\,ds
\right)\prod_{0<  S_n(\omega')\le t}\nu^\omega_{  S_n(\omega')}(\omega', \eta_n(\omega')),\quad \omega'\in\Omega',\,t\ge 0.
\end{eqnarray*}
Let $\tilde N\in\Fc$ and $N'_\omega\in\calf'$ be the  $\P$-null sets given by   Lemma \ref{L:nu_pred}-(ii). 
Then, for $\omega\notin\tilde N$, we have $\kappa_t^{\nu^\omega}(\omega')=\kappa_t^\nu(\omega,\omega')$, for all $ t\ge0$
and 
$ \omega' \in N'_\omega$. For later use we record the following fact:
\begin{align}\label{kappanuomega}
    \omega\notin\tilde N\qquad\Longrightarrow\qquad
\kappa_T^{\nu^\omega}=\kappa_T^\nu(\omega,\cdot),\,\P'-a.s.
\end{align}

Now fix   $\omega\in\Omega$. The representation 
 \eqref{representation_nu^omega}  shows that the process $\nu^\omega$ is $\calp(\F^\mu)\otimes \calb(A)$-measurable, so we can apply 
Proposition  
\ref{cambiosequenzasua} on the canonical probability space  $(\Omega',\Fc',\P')$ introduced  above.
By that proposition, $(\kappa_t^{\nu^\omega})_{t\ge 0}$ is a martingale
with respect to $\P'$ and $\F^\mu$
and there exists
 a unique probability measure $\P^{\nu^\omega}$ on $(\Omega',\calf^\mu_\infty)$ such that $\P^{\nu^\omega}(d\omega')=\kappa_t^{\nu^\omega}(\omega')\P'(d\omega')$
 on each $\sigma$-algebra $\calf^\mu_t$. Moreover,  the $\F^\mu$-compensator of $
\mu(\omega')=\sum_{n\ge1}\delta_{(T_n(\omega'),A_n(\omega'))},$ under $\P^{\nu^\omega}$
 is given by 
 $\nu^\omega_t(\omega',a)\,\lambda(da)\,dt$.

\vspace{2mm}

We can now state the following key  result  from  which the required conclusion  of this Subsection
follows readily (see Proposition \ref{P:Ineq_I}). 

\begin{Lemma}
\label{L:Replication_I} Given $\nu\in\Vc_{\inf\,>\,0}$, there exist a  sequence $(T_n^\nu,A_n^\nu)_{n\geq1}$ on $(\bar\Omega,\bar\Fc,\bar\P)$ such that:
\begin{itemize}
\item[\textup{(i)}] for every $n\geq1$, $(T_n^\nu,A_n^\nu)$ takes values in $(0,\infty)\times A$ and $T_n^\nu<T_{n+1}^\nu$;
\item[\textup{(ii)}] for every $n\geq1$, $T_n^\nu$ is an $\F^{W,\mu'_\infty}$-stopping time and $A_n^\nu$ is $\Fc_{T_n^\nu}^{W,\mu'_\infty}$-measurable;
\item[\textup{(iii)}] $\lim_{n\rightarrow\infty}T_n^\nu=\infty$;
\item[\textup{(iv)}] for every $\omega\in \Omega$, we have
\[
\mathscr L_{\P'}\big((T_n^\nu(\omega,\cdot),A_n^\nu(\omega,\cdot))_{n\geq1}\big) \ = \ \mathscr L_{\P^{\nu^\omega}}\big((T_n,A_n)_{n\geq1}\big).
\]
\end{itemize}
Finally, let $\bar\alpha_t^\nu=a_01_{[0,T_1^\nu)}+\sum_{n=1}^\infty A_n^\nu 1_{[T_n^\nu,T_{n+1}^\nu)}(t)$ be the step process associated with $(T_n^\nu,A_n^\nu)_{n\geq1}$. Then, $\bar\alpha^\nu\in\Ac^{W,\mu'}$ and $\mathscr L_{\P'}(\bar\alpha^\nu(\omega,\cdot))=\mathscr L_{\P^{\nu^\omega}}(I)$, $\omega\in\Omega$.
\end{Lemma}
\textbf{Proof.} Below we use the shortened notation $\bar\omega=(\omega,\omega')$.
Suppose that we have already constructed a multivariate point process $(T_n^\nu,A_n^\nu)_{n\geq1}$ satisfying points (i)-(ii)-(iii)-(iv) of the Theorem. Then, by (ii) it follows that $\bar\alpha^\nu$ is c\`adl\`ag and $\F^{W,\mu'_\infty}$-adapted, hence progressive. Moreover, by (iii), for every $(\bar\omega,t)\in\bar\Omega\times[0,T]$ the series $\sum_{n=1}^\infty A_n^\nu(\bar\omega) 1_{[T_n^\nu(\bar\omega),T_{n+1}^\nu(\bar\omega))}(t)$ is a finite sum, and thus $\bar\alpha^\nu$
$\in$   $\Ac^{W,\mu'}$. Furthermore, by (iv) we see that $\mathscr L_{\P'}(\bar\alpha^\nu(\omega,\cdot))=\mathscr L_{\P^{\nu^\omega}}(I)$, $\omega\in\Omega$.

Let us now construct $(T_n^\nu,A_n^\nu)_{n\geq1}$ satisfying points (i)-(ii)-(iii)-(iv).
For the given $\nu\in\Vc_{\inf\,>\,0}$, 
consider again, for 
$\omega\in\Omega$,
the random fields $\nu^\omega$
introduced in
 Lemma \ref{L:nu_pred}-(ii). By 
Proposition
\ref{cambiosequenzasua},  applied on the canonical space $(\Omega',\calf',\P')$,
we conclude that
there exists   a sequence $(T_n^{\nu^\omega},A_n^{\nu^\omega})_{n\ge1}$ of functions defined on $\Omega'$,  
taking values in $(0,\infty)\times A$,
measurable with respect to $\calf_\infty^{\mu'}=\calf'$,
such that
$T_n^{\nu^\omega} <T_{n+1}^{\nu^\omega}\to\infty$ and
\[
\mathscr L_{\P'}\big((T_n^{\nu^\omega},A_n^{\nu^\omega})_{n\geq1}\big) \ = \ \mathscr L_{\P^{\nu^\omega}}\big((T_n,A_n)_{n\geq1}\big).
\]
We define, for $\omega\in\Omega$, $\omega'\in\Omega'$, $n\ge1$:
\[
T_n^\nu(\omega,\omega')= 
T_n^{\nu^\omega}(\omega'),
\quad
 A_n^\nu(\omega,\omega')=
A_n^{\nu^\omega}(\omega').
\]
Clearly,  the constructed sequence 
$(T_n^\nu(\bar\omega),A_n^\nu(\bar\omega))_{n\ge1}$ satisfies  (i)-(iii)-(iv). To finish the proof it remains to prove that it
satisfies (ii) as well, namely that
$T_n^\nu$ is an $\F^{W,\mu'_\infty}$-stopping time and $A_n^\nu$ is $\Fc_{T_n^\nu}^{W,\mu'_\infty}$-measurable, when they are thought as functions of $\bar\omega=(\omega,\omega')$. 
We will verify this by inspection of the explicit formulae used for the construction of $T_n^\nu$ and $A_n^\nu$, introduced in the proofs of 
  Propositions
\ref{cambiosequenzasur}
and 
\ref{cambiosequenzasua}, that now will be written  with the present notation.

First we recall that 
$T^\nu_1(\bar\omega)$ was defined 
in \eqref{Tnu1}. With the present notation and recalling the random field introduced in \eqref{fieldaux}, formula  \eqref{Tnu1} becomes
\[
T^\nu_1(\bar\omega)  = \inf\{t\ge 0\,:\,\theta^{(1)}_t(\omega)\ge T_1(\omega')\}, \qquad {\rm where} \qquad 
\theta^{(1)}_t(\omega)=
\frac{1}{\lambda'(\R)} \int_0^t\int_\R \nu_s^{(0)}(\omega,\pi(r)) \lambda'(dr)\,ds.
\]

 Let $\bar E_{T_1}:=\{(\bar\omega,t)\in\bar\Omega\times\R_+\colon \theta_t^{(1)}(\omega)=T_1(\omega')\}$. Since
the process $(\bar\omega,t)\mapsto (\theta_t^{(1)}(\omega),T_1(\omega'))$
is $\F^{W,\mu_\infty}$-adapted and continuous,
 $\bar E_{T_1}$ is an $\F^{W,\mu_\infty}$-optional set (in fact, predictable). Since $T_1^\nu(\bar\omega)=\inf\{t\in\R_+\colon(\bar\omega,t)\in\bar E_{T_1}\}$
  is the d\'ebut of $\bar E_{T_1}$ it follows that $T_1^\nu$ is an $\F^{W,\mu_\infty}$-stopping time, by a standard result in the general theory of stochastic processes (see e.g.
   \cite{jacod_book} Theorem 1.14).

Now let us consider $A_1^\nu$. It was defined in 
 \eqref{anproiettati} as
$A_1^\nu=\pi(R_1^{\nu'})$, where $R_1^{\nu'}$ was defined in
\eqref{Anu1}. 
With the present notation and recalling the random field introduced in \eqref{fieldaux} we have
\[
R_1^{\nu'}  = \inf\{b\in \R\,:\,F_b^{(1)}(\bar\omega) \ge F'(R_1(\omega'))\},
\qquad {\rm where} \qquad 
F_b^{(1)}(\bar\omega) \ := \ \frac{\int_{-\infty}^b \nu_{T_1^\nu(\bar\omega)}^{(0)}(\omega,\pi(r))\,\lambda'(dr)}{\int_{-\infty}^{+\infty}\nu_{T_1^\nu(\bar\omega)}^{(0)}(\omega,\pi(r))\,\lambda'(dr)}
\]
and $F'$ is now the cumulative distribution function of $\lambda'$, namely $F'(r)=\lambda'((-\infty,r])/\lambda'(\R)$. 
We note that the process
$$
(\bar\omega,t)\mapsto
\frac{
\int_{-\infty}^b \nu_{t}^{(0)}(\omega,\pi(r))\,\lambda'(dr)
}{
\int_{-\infty}^{+\infty}\nu_{t}^{(0)}(\omega,\pi(r))\,\lambda'(dr)
}
$$
is predictable with respect to $\F^W$, hence
it is also $\F^{W,\mu_\infty}$-progressive.
Substituting $t$ with $T_1^\nu(\bar\omega)$ we conclude that
$F_b^{(1)}$ is $(\Fc_{T_1^\nu}^{W,\mu_\infty})$-measurable.
Since $R_1$ is clearly $\calf'$-measurable and
$\calf'\subset \Fc_{0}^{W,\mu_\infty}\subset \Fc_{T_1^\nu}^{W,\mu_\infty}$,
$A_1$  is also $(\Fc_{T_1^\nu}^{W,\mu_\infty})$-measurable.
Recalling   $b\mapsto F_b^{(1)}(\bar\omega)$ was continuous and strictly increasing, 
it is easy to conclude that $R_1^{\nu'}$ is $(\Fc_{T_1^\nu}^{W,\mu_\infty})$-measurable.
This implies that $A_1^\nu$ is also $\Fc_{T_1^\nu}^{W,\mu_\infty}$-measurable.

Next we consider
$T^\nu_2(\bar\omega)$. It was defined 
in \eqref{Tnu2}, that now becomes, recalling \eqref{fieldaux},
\[
T^\nu_2 (\bar\omega) = \inf\{t\ge T^\nu_1(\bar\omega)\,:\,\theta^{(2)}_t(\bar\omega)\ge T_2(\omega')-T_1(\omega')\},
\]
where
\[
\theta^{(2)}_t(\bar\omega)=
\frac{1}{\lambda'(\R)} \int_{T_1^\nu(\bar\omega)}^{t}\int_A \nu_s^{(1)}(\omega,T_1^\nu(\bar\omega),A_1^\nu(\bar\omega),\pi(r)) \,\lambda'(dr)\,ds .
\]

Taking into account the results already proved for $T_1^\nu(\bar\omega)$ and $A_1^\nu(\bar\omega)$, since $T_{2}^\nu$ is the d\'ebut of $\bar E_{T_{2}}:=\{(\bar\omega,t)\in\bar\Omega\times\R_+\colon \theta_t^{(2)}(\bar\omega)=T_{2}(\omega')-T_{1}(\omega')\}$,   it is an $\F^{W,\mu_\infty}$-stopping time.

Now let us consider $A_{2}^\nu$. It was defined in 
 \eqref{anproiettati} as
$A_2^\nu=\pi(R_2^{\nu'})$, where $R_2^{\nu'}$ 
was given in \eqref{Anu2} by the formula
\[
A^\nu_2  = \inf\{b\in \R\,:\,F_b^{(2)}(\bar\omega) \ge F'(R_2(\omega'))\},
\]
  where
\[ 
F_b^{(2)}(\bar\omega) =   \frac{\int_{-\infty}^b \nu_{T_2^\nu(\bar\omega) }^{(1)}(\omega,T_1^\nu(\bar\omega), A_1^\nu(\bar\omega), \pi(r))\,\lambda'(dr)}{\int_{-\infty}^{+\infty}\nu_{T_2^\nu(\bar\omega) }^{(1)}(\omega, T_1^\nu(\bar\omega), A_1^\nu(\bar\omega), \pi(r))\,\lambda'(dr)}
\]
and $F'$ is the same as before.
Recalling that the function $b\mapsto F_b^{(2)}(\bar\omega) $ is continuous and strictly increasing and arguing as before we conclude that 
$R_2^{\nu'}$, and hence 
 $A_{2}^\nu$, is
 $\Fc_{T_{2}^\nu}^{W,\mu_\infty}$-measurable.

Iterating these arguments, or using a formal inductive proof, we conclude the
 proof of the proposition.
\qed

\vspace{3mm}

We can now prove the main result of this subsection and conclude the proof of the inequality ${\text{\Large$\upsilon$}}_0^\Rc$ $\le$  ${\text{\Large$\upsilon$}}_0$.

\begin{Proposition}
\label{P:Ineq_I}
For every $\nu\in\Vc_{\inf\,>\,0}$ there exists $\bar\alpha^\nu\in\Ac^{W,\mu'}$ such that
\begin{equation}\label{LawB_I=LawB_alpha}
\mathscr L_{\P^\nu}(W,I) \ = \ \mathscr L_{\bar\P}(W,\bar\alpha^\nu).
\end{equation}
In particular,   $W$ is a standard Wiener process under $\P^\nu$ and we have
\begin{equation}\label{LawX_I=LawX_alpha}
\mathscr L_{\P^\nu}(X,I)  =  \mathscr L_{\bar\P}(X^{\bar\alpha^\nu},\bar\alpha^\nu),
\qquad J^\Rc(\nu)=\bar J(\bar\alpha^\nu).
\end{equation}
Finally, ${\text{\Large$\upsilon$}}_0^\Rc$ $:=$ $\sup_{\nu\in\Vc}J^\Rc(\nu)$ $\leq$ $\sup_{\bar\alpha\in\Ac^{W,\mu'}}\bar J(\bar\alpha)$  $=$
$\sup_{\alpha\in\Ac^W}J(\alpha)$ $=:$ ${\text{\Large$\upsilon$}}_0$.
\end{Proposition}
\textbf{Proof.}
Suppose that \eqref{LawB_I=LawB_alpha} holds. Then, from equation \eqref{stateeq} and Assumption {\bf (A1)} it is well-known that the first
equality in \eqref{LawX_I=LawX_alpha} holds as well, and this   implies the second equality.
 From the arbitrariness of $\nu\in\Vc_{\inf\,>\,0}$, we deduce that $\sup_{\nu\in\Vc_{\inf\,>\,0}}J^\Rc(\nu)\leq\sup_{\bar\alpha\in\Ac^{W,\mu'}}\bar J(\bar\alpha)$. Since  $\sup_{\nu\in\Vc_{\inf\,>\,0}}J^\Rc(\nu)=\sup_{\nu\in\Vc}J^\Rc(\nu)$
 by \eqref{eqinfzero}, we conclude that $\sup_{\nu\in\Vc}J^\Rc(\nu)\leq\sup_{\bar\alpha\in\Ac^{W,\mu'}}\bar J(\bar\alpha) = \sup_{\alpha\in\Ac^W}J(\alpha)$, where the last equality follows from Lemma \ref{L:AWmu}.

Let us now prove \eqref{LawB_I=LawB_alpha}. Fix $\nu\in\Vc_{\inf\,>\,0}$ and consider the process $\bar\alpha^\nu$ given by Lemma \ref{L:Replication_I}. In order to prove \eqref{LawB_I=LawB_alpha}, we have to show that
\begin{equation}\label{psi_phi}
\bar\E[\kappa_T^\nu \psi(W)\phi(I)] \ = \ \bar\E[ \psi(W)\phi(\bar\alpha^\nu)],
\end{equation} 
for any real valued  bounded Borel function $\psi$ defined on the space of continuous trajectories  $[0,T]\to \R^{d}$ (endowed with the supremum norm) and any real valued bounded Borel  function $\phi$ on the space  of c\`adl\`ag paths from $[0,T]$ to $A$ (endowed with the supremum norm). By the Fubini theorem,  \eqref{psi_phi} can be rewritten as
\begin{align*}
&\int_\Omega  \psi(W(\omega)) \bigg(\int_{\Omega'} \kappa_T^\nu(\omega,\omega') \phi(I(\omega')) \P'(d\omega')\bigg) \P(d\omega) \\
&= \ \int_\Omega  \psi(W(\omega)) \bigg(\int_{\Omega'} \phi(\bar\alpha^\nu(\omega,\omega')) \P'(d\omega')\bigg) \P(d\omega).
\end{align*}
Let $\tilde N\in\Fc$ be as in Lemma \ref{L:nu_pred}. Then it is enough to prove that $\E'[\kappa_T^\nu(\omega,\cdot)\phi(I)]=\E'[\phi(\bar\alpha^\nu(\omega,\cdot))]$ whenever $\omega\notin\tilde N$. Recalling 
\eqref{kappanuomega}, 
this is equivalent to
\[
\E' [\kappa_T^{\nu^\omega}\phi(I)] \ = \ \E'[\phi(\bar\alpha^\nu(\omega,\cdot))], \qquad  \omega\notin\tilde N.
\]
By  the definition of $\P^{\nu^\omega}$ this is the same as:
\[
\E^{\nu^\omega}[\phi(I)] \ = \ \E'[\phi(\bar\alpha^\nu(\omega,\cdot))], \qquad  \omega\notin\tilde N,
\]
and this follows immediately from the last statement of Lemma \ref{L:Replication_I}.
\qed

\subsection{Proof of the inequality ${\text{\Large$\upsilon$}}_0$ $\le$ ${\text{\Large$\upsilon$}}_0^\Rc$}
\label{secVleqVR}

This proof uses arguments introduced in  Proposition 4.1 in \cite{FP15} and Section 4.2. in \cite{BCFP16c}.
It is based on a general point process construction reported in the Appendix.

We still suppose that assumptions
{\bf (A1)} and {\bf (A2)} hold true. 
Let us consider the optimal control problem  in Section  \ref{Primal}, formulated in a  setting $(\Omega,\calf,\P, \F, W)$. 

In the proof we need to consider the continuous dependence of the reward functional on the control process. To this end we will introduce a metric in the space of   control processes. This can be done in various ways. For instance, in 
\cite{YongZhou99} Lemma 6.4
  a metric on controls was used in connection with the stochastic maximum principle
in the presence of state constraints.
More sophisticated constructions are used in the framework of so-called relaxed controls, see for instance
\cite{YongZhou99} Section 5.3. 
  Here we will use
the metric   introduced in
 \cite{80Krylov} and presented below.

 Recall that the space $A$ of control actions was assumed to be a Borel space. Let $\rho$ denote a metric inducing its topology. Without loss of generality we will assume that $\rho(x,y)<1$ for all $x,y\in A$: if this is not the case we replace $\rho$ by the metric $(x,y)\mapsto \rho(x,y)/(1+\rho(x,y))$, which is equivalent, in the sense that it induces the same topology. 
Next define, for any pair $\alpha^1,\alpha^2:\Omega\times [0,T]\to A$ of measurable processes (namely, measurable with respect to $\calf\otimes\calb([0,T])$ and $\calb(A)$):
\begin{eqnarray*}
\tilde \rho(\alpha^1,\alpha^2) &=& \E\Big[\int_0^T\rho(\alpha^1_t,\alpha^2_t)\,dt \Big].
\end{eqnarray*}
Note that a sequence $\alpha^n$ converges  to a limit $\alpha$ with respect to this metric if and only if $\alpha^n\to\alpha$
in $dt\otimes d\P$-measure, i.e. if and only if for  any  $\epsilon>0$ we have
\begin{eqnarray*}
\lim_{n\to\infty}(dt\otimes d\P)(\{(t,\omega)\in[0,T]\times\Omega\,:\, \rho(\alpha^n_t(\omega),\alpha_t(\omega))>\epsilon\}) &=& 0,
\end{eqnarray*}
Recall that the space of controls was denoted  $\Ac^W$ and consists of $A$-valued processes progressive for the completed filtration  $\F^W$ generated by $W$. The state equation was introduced in \eqref{stateeq} 
and the gain functional $J(\alpha)$ was defined in
 \eqref{gaineq}.
The following continuity result  of the gain functional with respect to the control  was proved in \cite{80Krylov}.

\begin{Lemma}\label{Jcontinuous}
The map $\alpha\mapsto 
J(\alpha) $ is continuous with respect to the metric $\tilde \rho$.
\end{Lemma}
{\bf Proof.} Suppose $\alpha^n,\alpha\in \Ac^W$ and  $\tilde\rho(\alpha^n,\alpha)\to 0$ or equivalently $\alpha^n\to\alpha$
in $dt\otimes d\P$-measure. We will prove that $J (\alpha^n)\to J(\alpha)$. Denote $X^n, X$ the corresponding
trajectories. Then, starting from the state equation \eqref{stateeq}, using usual arguments involving the Burkholder-Davis-Gundy inequalities
and the Gronwall lemma, for every $p\in [{2},\infty)$ we arrive at
\begin{eqnarray*}
\E \Big[\sup_{t\in[0,T]}|X^n_t-X_t|^p \Big] &\le&  C \, \E \Big[ \int_0^T \Big\{ | b(t,X_t,\alpha^n_t)-b(t,X_t,\alpha_t)|^p
\\&&
+ \,| \sigma(t,X_t,\alpha^n_t)-\sigma(t,X_t,\alpha_t)|^p \Big\} \,dt \Big] ,
\end{eqnarray*}
for a suitable constant $C$, independent of $n$.  Recalling the bound \eqref{boundlponx} on the solutions $X$ and $X_n$, using the linear growth of the coefficients $b$ and $\sigma$ that follows from
{\bf (A1)}-\eqref{lipbsig} and {\bf (A1)}-\eqref{borbsig},
 by the dominated convergence theorem  we first conclude that $\E\Big[\sup_{t\in[0,T]}|X^n_t-X_t|^p\Big]$ $\to 0$ as $n\to\infty$.
Next we have
\begin{eqnarray*}
|J( \alpha^n)-J( \alpha)| \le
\E \,\Big[ \int_0^T|f(t,X^n_t,\alpha^n_t)- f(t,X_t,\alpha_t)|\,dt \Big]
+ \E\,|g(X^n_T)-g(X_T)|.
\end{eqnarray*}
To finish the proof we show that the right-hand side tends to zero.
Suppose, on the contrary,   that there exist $\eta>0$ and a subsequence (denoted $(X^{n'},\alpha^{n'})$) such that
\begin{equation} \label{contradictiondominated}
\E \Big[\int_0^T|f(t,X_t^{n'},\alpha^{n'}_t)- f(t,X_t,\alpha_t)|\,dt \Big] \; \ge \;  \eta,
\end{equation}
for every $n'$. Passing to a sub-subsequence, still denoted by the same symbol,  we can assume that
\begin{eqnarray*}
\sup_{t\in[0,T]}|X^{n'}_t-X_t|\to 0, \;\;\;d\P-a.s.,\quad
\rho(\alpha^{n'}_t,\alpha_t)\to 0,\;\;\; dt\otimes d\P-a.e.
\end{eqnarray*}
as $n'\to\infty$, and by the   continuity properties of $f$ assumed in {\bf (A1)} it follows that
$f(t,X_t^{n'},\alpha^{n'}_t)\to f(t,X_t,\alpha_t)$,  $dt\otimes d\P$-a.e.
Next we extract a further subsequence $(n'_j)$ such that
\begin{eqnarray*}
\left(\E \Big[ \sup_{t\in[0,T]}|X^{n'_j}_t-X_t|^p \Big] \right)^{1/p} &\le&  2^{-j},
\end{eqnarray*}
so that the random variable $\bar X:=\sum_j \sup_{t\in[0,T]}|X^{n'_j}_t-X_t|$ satisfies $\E|\bar X|^p<\infty$ as well as
$|X^{n'_j}_t|\le |X_t|+ |\bar X|$ for every $t$ and $j$. Recalling the polynomial growth condition \eqref{PolGrowth_f_g} of index $r$ imposed on $f$   in {\bf (A1)}, we obtain
\begin{eqnarray*}
|f(t,X_t^{n'_j},\alpha^{n'_j}_t)- f(t,X_t,\alpha_t)| & \le &  C(1+\sup_{t\in[0,T]}|X^{n'_j}_t|^r
+ \sup_{t\in[0,T]}|X_t|^r) \\
&\le &  C(1+ |\bar X|^r+ \sup_{t\in[0,T]}|X_t|^r)
\end{eqnarray*}
for a suitable constant $C$, and choosing $p$ $=$ $\max(r,2)$,  we conclude that the right-hand side is integrable, which gives a contradiction with
\eqref{contradictiondominated} by the dominated convergence theorem. This shows that
$\E \big[\int_0^T|f(t,X_t^n,\alpha^n_t)- f(t,X_t,\alpha_t)|\,dt\big]$ $\to 0$, and in a similar way one shows that $\E|g(X^n_T)- g(X_T)|\to 0$.
\qed

Let us denote by $\cala^W_0$ the set of processes $\alpha\in\cala^W$ of the form $$\alpha_t=\sum_{n= 0}^{N-1} \alpha_{n}1_{ [t_n,t_{n+1})}(t),
$$
where  $0=t_0<t_1<\ldots t_N=T$ is a deterministic subdivision of $[0,T]$,
$\alpha_0,\ldots,\alpha_{N-1}$ are $A$-valued
random variables that take only a finite number of values, and each $\alpha_n$ is $\calf^W_{t_n}$-measurable.
The following result is 
 Lemma 3.2.6 in \cite{80Krylov};
the proof is not long but it will not be reported here.

\begin{Lemma}\label{densitykrylov}
The set $\cala_0^W$  is dense in the set of admissible controls $\cala^W$ with respect to the norm $\tilde\rho$.
\end{Lemma}

Together with 
Lemma \ref{Jcontinuous}
this proves that, instead of 
\eqref{primalvalue}, 
the value  ${\text{\Large$\upsilon$}}_0$ of the control problem can be defined by taking the supremum of $J(\alpha)$ over all $\alpha\in\Ac^W_0$.
This will be used in the proof of
Proposition \ref{extensionapproximation} in the Appendix, which is a crucial step in the proof that follows.

After these preliminaries we can conclude the proof of the desired inequality.

\vspace{2mm}

\noindent {\bf Proof of the inequality ${\text{\Large$\upsilon$}}_0$ $\le$ ${\text{\Large$\upsilon$}}_0^\Rc$.}

Suppose we are given a setting $(\Omega,\calf,\P, \F, W)$ for the original optimal control problem, satisfying the conditions in Section \ref{Primal}, and  consider the controlled equation \eqref{stateeq} and the gain
\eqref{gaineq}. We fix an $\F^W$-progressive process $\alpha$ with values in $A$.  We will show how to construct a sequence of settings
$(\hat \Omega, \hat \calf,\hat \P_k,  \hat W, \hat \mu_k )_k$
for the randomized control problem of Section \ref{randomizedformulation},
and a sequence $(\hat\nu^k)_k$ of corresponding admissible controls
(both sequences depending on $\alpha$),
  such that for the corresponding gains, defined by \eqref{defJrandomized}, we have:
\begin{eqnarray}  \label{desiredineq}
J^\Rc(\hat\nu^k) & \rightarrow &  J(\alpha), \;\;\; \mbox{ as } \; k \rightarrow \infty.
\end{eqnarray} 
Admitting for a moment that this has been done, the proof
is easily concluded by the following arguments. By \eqref{desiredineq}, we can find,
 for any $\varepsilon $ $>$ $0$, some $k$ such that $J^\Rc(\hat\nu^k)$ $>$ $J(\alpha)-\varepsilon $.
Since $J^\Rc(\hat\nu^k)$ is a gain of a randomized control problem, it can not exceed the value ${\text{\Large$\upsilon$}}_0^\Rc$ defined in \eqref{dualvalue} which,
by Proposition \ref{indepofthesetting}, does not depend on $\varepsilon$ nor on $\alpha$.
It follows that
\begin{eqnarray*}
   {\text{\Large$\upsilon$}}_0^\Rc &>&  J(\alpha)-\varepsilon
\end{eqnarray*}
and by the arbitrariness of $\varepsilon $ and $\alpha$, we obtain the required inequality ${\text{\Large$\upsilon$}}_0^\Rc$ $\ge$ ${\text{\Large$\upsilon$}}_0$.

In order to construct the sequences $(\hat \Omega, \hat \calf,\hat \P_k,  \hat W, \hat \mu_k )_k$ and   $(\hat \nu^k)_k$
satisfying \eqref{desiredineq}, as a first step
 we apply Proposition \ref{extensionapproximation} in the Appendix  to the
probability  space $(\Omega,\calf,\P)$ of the original control problem and to the filtration $\G:=\F^W$. The fact that $\lambda(da)$ has full topological support, contained in assumption {\bf (A2)}, is used there.
In that Proposition a suitable probability space $(\Omega',\calf',\P')$ is introduced and the product space $(\hat\Omega,
\hat \calf,\Q)$ is constructed:
\begin{eqnarray*}
\hat\Omega \; = \; \Omega\times \Omega',
\qquad
\hat \calf \; = \; \calf\otimes \calf',
\qquad
\Q \; = \; \P\otimes \P'.
\end{eqnarray*}
Then   
the processes $\alpha$ and $W$, which are random processes in  
$(\Omega,\calf)$, can
be extended to random processes 
$\hat\alpha$ and $ \hat W$ in $(\hat\Omega,
\hat \calf)$
in a natural way setting
\[ \hat\alpha_t(\hat\omega)=\alpha_t(\omega),\quad \hat W_t(\hat\omega)=W_t(\omega),\qquad
{\rm where} \qquad
 \hat\omega=(\omega,\omega'),
\]
and the extension $ \hat W$ of $W$  remains a Wiener process
under $\Q$. 
The filtration $\F^W$ can also be canonically extended to the  filtration 
 $\hat\calg_t=\{A\times \Omega'\,:\, A\in\calg_t\}$ in $(\hat\Omega,
\hat \calf)$.
It is not difficult to verify that
this filtration  coincides with the filtration $\F^{\hat W}$ generated
by $\hat W$.

The metric $\tilde\rho$ can also be extended to
  any pair  $\beta^1,\beta^2:\hat\Omega\times [0,T]\to A$ of   measurable
  processes  setting
\begin{eqnarray*}
\tilde \rho(\beta^1,\beta^2) &=& \E^\Q \Big[\int_0^T\rho(\beta^1_t,\beta^2_t)\,dt \Big],
\end{eqnarray*}
where  $\E^\Q$ denotes the expectation under  $\Q$.
Note that we use the same symbol $\tilde\rho$ to denote  the extended metric as well.

Then, by Proposition \ref{extensionapproximation}, for any integer $k\ge1$
 there exists a marked point process $(\hat S_n^k,\hat \eta^k_n)_{n\ge 1}$ defined in
$(\hat\Omega,
\hat \calf,\Q)$
satisfying the following conditions.
\begin{enumerate}
\item
Setting
$
\hat S_0^k=0$, $\hat \eta_0^k=a_0$, $
\hat I_t^k=\sum_{n\ge 0}\hat  \eta^k_{n}1_{ [\hat S^k_n,\hat S^k_{n+1})}(t)$,
we have
$    \tilde\rho (\hat I^k, \hat \alpha)<1/k$.
\newline\noindent (We recall that $a_0$ is the point in 
assumption  {\bf (A2)})

\item Denote
$\hat \mu_k=\sum_{n\ge1}\delta_{(\hat S_n^k,\hat \eta^k_n)}$ the random measure
associated to $(\hat S_n^k,\hat \eta_n^k)_{n\ge 1}$ and
$\F^{\hat\mu_k}=(\calf_t^{\hat\mu_k})_{t\geq 0}$
 the natural filtration of
$\hat\mu_k$;
then
the  compensator of    $\hat \mu_k$ under $\Q$
with respect to the filtration $ (\calf^{\hat W}_t\vee\calf_t^{\hat\mu_k})_{t\geq 0}$
is absolutely continuous
 with respect to $\lambda(da)\,dt$ and
it can be written in the form
$$
 \hat\nu_t^k(\hat\omega,a)\, \lambda(da)\,dt
$$
 for  some  nonnegative $\calp(\F^{\hat W,\hat \mu})
 \otimes \calb(A)$-measurable  function $\hat\nu^k$ satisfying
\[
0< \inf_{\hat\Omega\times [0,T]\times A}\hat\nu^k\le 
  \sup_{\hat\Omega\times [0,T]\times A}\hat\nu^k<\infty.
  \]
\end{enumerate}
We note that $\hat\mu_k$ (and so also $\hat I^k$ and $\hat\nu^k$)
depend on $\alpha$ as well, but we do not make it explicit in the notation.
Let us now consider the completion of the probability
space $(\hat\Omega,\hat\calf,\Q)$, that will be denoted
by the same symbol for simplicity of notation, and let
$\caln$ denote the family of $\Q$-null sets of the completion.
Then the filtration
$(\calf^{\hat W}_t\vee\calf_t^{\hat\mu_k}\vee \caln)_{t\geq 0}$
coincides with  the filtration previously denoted by
$\F^{\hat W,\hat \mu_k}=(\calf^{\hat W,\hat \mu_k}_t)_{t\ge 0}$
(compare with formula  \eqref{expandedfiltration} in
  section \ref{randomizedformulation}). It is easy to see that
$ \hat\nu_t^k(\hat\omega,a)\, \lambda(da)\,dt$ is the compensator
of $\hat\mu_k$  with respect to $\F^{\hat W,\hat \mu_k}$ and the extended
probability $\Q$ as well.

We recall that our aim is to construct 
 the sequences $(\hat \Omega, \hat \calf,\hat \P_k,  \hat W, \hat \mu_k )_k$ and   $(\hat \nu^k)_k$
satisfying \eqref{desiredineq}. So far we have constructed 
 $(\hat \Omega, \hat \calf,\Q,  \hat W, \hat \mu_k )_k$ and   $(\hat \nu^k)_k$. We will
change the probability $\Q$ and define an equivalent probability measure $\hat \P_k$ by an application of the Girsanov theorem for point processes, 
 Theorem  \ref{girsanovpoint}.

The random field 
$(\hat\nu^k)^{-1}$ 
is $\calp(\F^{\hat W,\hat\mu_k})\otimes\calb(A)$-measurable. Consider
 the corresponding exponential process
\begin{equation}\label{defgirsanovdensity}
 M_t^k \; := \; \exp\Big(\int_0^t\int_A(1-\hat\nu^k_s(a)^{-1})\,
 \hat\nu_s^k(a) \lambda(da)\,ds\Big) \prod_{\hat S_n^k\le t} \hat\nu^k_{\hat S^k_n}
 (\hat\eta^k_n)^{-1},\qquad t\in [0,T],
\end{equation}
as defined in   
 Theorem  \ref{girsanovpoint}-(i)
(to match the notation of that theorem,  $\mu$ is now replaced by 
$\hat\mu^k$, the random field
$\nu$ is now  
$(\hat\nu^k)^{-1}$ and the compensator $\mu^\P(dt\,da)$ is now 
$\hat\nu_t^k(a) \lambda(da)\,dt$).
By point (i) of
 Theorem  \ref{girsanovpoint}, 
$M^k$ is a   positive supermartingale (with respect to $\F^{\hat W,\hat\mu_k}$ and $\Q$) on $[0,T]$. Since, however, we have $\inf\hat\nu^k>0$, 
the random field $(\hat\nu^k)^{-1}$ and bounded and it is easy to prove that $\E^\Q[M_T^k]=1$ 
(or one can refer to  \cite{bremaud2} Remark 5.5.2 for a careful verification).
According to   point (ii) of
 Theorem  \ref{girsanovpoint},
$M^k$
is in fact a strictly positive martingale and we can
define an equivalent probability $\hat\P_k$ on the space $(\hat\Omega,\hat\calf)$ setting
 $\hat\P_k(d\hat\omega)$ $=$ $M_T^k(\hat\omega)\Q(d\hat\omega)$. The expectation under
 $\hat\P_k$ will be denoted  $\hat\E_k$. 
Now we make the following claims:
\begin{enumerate}
  \item [(i)] On the set $(0,T]\times A$
 the random measure 
 $\hat\mu_k$      
   has $(\hat\P_k,\F^{\hat W,\hat\mu_k})$-compensator $\lambda(da)\,dt$; moreover, 
it is a Poisson random measure with intensity $\lambda(da)$.

  This follows from  point (ii) of
 Theorem  \ref{girsanovpoint},
  which guarantees that
  under the new probability $\hat\P_k$   the   compensator of $\hat\mu_k$   is given by
  $
  (\hat\nu^k)^{-1}\cdot \hat\nu^k \,\lambda(da)\,dt $ $=$ $ \lambda(da)\,dt.
  $ The Poisson character follows from
Watanabe's Theorem \ref{watanabe}.

   \item [(ii)] $\hat W$ is a $(\hat\P_k,\F^{\hat W,\hat\mu_k})$-Wiener process.

The proof is as follows. Since  the probabilities
$\hat\P_k$ and $\Q$ are equivalent,  the quadratic variation of $\hat W$ computed under $\hat\P_k$ and $\Q$
is the same, and equals $\langle \hat W\rangle_t=t$. So it is enough to show
that $\hat W$ is a $(\hat\P_k,\F^{\hat W,\hat\mu_k})$-local martingale, which is equivalent
to the fact that $M\hat W$ is a $(\Q,\F^{\hat W,\hat\mu_k})$-local martingale. In turn, this follows from a general fact: since $M$ is a $(\Q,\F^{\hat W,\hat\mu_k})$-martingale
of finite variation, it is purely discontinuous  and
therefore   orthogonal (under $\Q$)
to $\hat W$; thus, their product $M\hat W$ is a $(\Q,\F^{\hat W,\hat\mu_k})$-local martingale.

  \item  [(iii)] $\hat W$ and $\hat \mu_k$ are independent under $\hat\P_k$.

To prove this claim it is enough to show that, for any Borel subset $B\subset A$, the counting process
 $  N^B_t  :=  \hat\mu_k((0,t]\times B)  $
is  independent from
  $\hat W$  under $\hat\P_k$.  From  claims (i) and (ii) it  follows that
  $N^B$ is a Poisson process and $\hat W$ is a Wiener process, both with
  respect to $\F^{\hat W,\hat\mu_k}$ and $\hat\P_k$. By a general result,
  see e.g. Theorem 11.43 in \cite{HeWangYanbook},  to check the independence it
  is enough to note that their right bracket $[N^B,\hat W]$ is null,
  which is obvious, since $\hat W$ is continuous and $N^B$ has
  no continuous part.
\end{enumerate}

From claims (i), (ii), (iii) we conclude that $(\hat \Omega, \hat \calf,\hat \P_k,  \hat W, \hat \mu_k )$
is indeed a setting  for a randomized
control problem as in Section \ref{randomizedformulation}.
Since $\hat\nu^k$ is a bounded,  strictly positive and
$\calp(\F^{W,\hat\mu})\otimes\calb(A)$-measurable
random field
  it  belongs to the class $\hat \calv^k$ of admissible controls for the
  randomized control problem and we now proceed
  to evaluating its gain $J^\Rc(\hat\nu^k)$ and to comparing it with $J(\alpha)$.
Our  aim is to show that,
as a consequence of the fact that
$    \tilde\rho (\hat I^k, \hat \alpha)<1/k$,
we have $J^\Rc(\hat\nu^k)   \to J(\alpha)$ as $k\to\infty$.

 We introduce the exponential process $\kappa^{\hat \nu^k}$
 corresponding to $\hat\nu^k$
  (compare formula
\eqref{doleans}):
\begin{equation}\label{doleansbis}
\kappa_t^{\hat\nu^k} =
 \exp\left(\int_0^t\int_A (1 - \hat \nu^k_s(a))\lambda(da)\,ds
\right)\prod_{\hat S_n^k\le t}\nu^k_{\hat S^k_n}(\hat \eta^k_n),\qquad t\in [0,T],
\end{equation}
we define
 the probability
$d\hat\P_k^{\hat\nu^k} =\kappa_T^{\hat\nu^k} d\hat\P_k $ and we obtain
the  gain
$$
J^\Rc(\hat\nu^k) =  \hat\E^{\hat\nu^k}
\left[\int_0^Tf(t,\hat X_t^k,\hat I^{k}_T)\,dt+g(\hat  X^k_T)\right],
$$
where      $\hat X^k$ is the
solution to
the equation
\begin{equation}\label{controlledhatkappa}
d\hat X_t ^k=  b(t, \hat X_t^k, \hat I^{k}_t)\,dt
+ \sigma (t,\hat X^k_t, \hat I^{k}_t)\,d\hat W_t,
\qquad
\hat X_0^k =\hat x_0.
\end{equation}
  However comparing \eqref{defgirsanovdensity}  and
 \eqref{doleansbis}  shows that
  $\kappa^{\hat \nu^k}_T\,M_T^k\equiv 1$, so that
 the  Girsanov transformation $\hat\P_k\mapsto \hat\P_k^{\hat\nu^k}$
is the inverse to the transformation $\Q\mapsto \hat\P_k$ made above,
and   changes back the probability $\hat\P_k$ into
$\Q$ considered above. Therefore we have
$\hat\P_k^{\hat\nu^k}=\Q$ and also
\begin{eqnarray}  \label{JR}
J^\Rc(\hat\nu^k) &=&  \E^\Q \left[ \int_0^Tf(t,\hat X_t^k,\hat I^{k}_t)\,dt+g(\hat X^k_T)\right].
\end{eqnarray} 
On the other hand, the gain $J(\alpha)$ of the initial control problem was defined in \eqref{gaineq} in terms
of the solution $X^\alpha$ to the controlled equation \eqref{stateeq}.
Denoting $\hat X^\alpha$ the extension
of $X^\alpha$
to $\hat\Omega$, it is easy to verify
that it is the solution to
\begin{equation}\label{controlledhat}
        d\hat X_t^\alpha \ = \ b(t,\hat  X_t^\alpha,\hat  \alpha_t)\,dt +
\sigma(t,\hat  X_t^\alpha, \hat \alpha_t)\,d\hat W_t, \qquad \hat X_0^\alpha=\hat x_0,
\end{equation}
and that
\begin{eqnarray}  \label{Ja}
J(\alpha) & = & \E^\Q\left[\int_0^Tf (t,\hat X_t^\alpha,\hat \alpha_t)\,dt+
g(\hat X^\alpha_T)\right].
\end{eqnarray} 
Equations \eqref{controlledhat} and \eqref{controlledhatkappa} are
 considered in the same probability space $(\hat\Omega,
\hat \calf,\Q)$.  In \eqref{controlledhatkappa}  we find a
solution adapted to the filtration
$\G^k:=\F^{  \hat W,\hat \mu_k}$
(defined as in
\eqref{expandedfiltration}) and in
\eqref{controlledhat}
 we find a
solution adapted to the filtration
$\G^0:=\F^{ \hat W}$ generated by   $\hat W$.

\vspace{2mm}

In order to conclude that $J^\Rc(\hat\nu^k)   \to J(\alpha)$ we need the following stability lemma, whose proof is entirely analogous to the proof of Lemma \ref{Jcontinuous} and therefore will be omitted.

\begin{Lemma}
\label{contrhotilde}
Given a probability space $(\hat\Omega,
\hat \calf,\Q)$ with filtrations
$\G^k=(\calg_t^k)_{t\ge 0}$ ($k\ge 0$)
consider
the equations
$$
d  Y_t ^k=  b(t, Y^k_t, \gamma^{k}_t)\,dt
+ \sigma(t,Y_t^k, \gamma_t^{k})\,d\beta_t,
\qquad
\hat Y_0^k =y_0\in\R^n,
$$
where $\beta$  is a Wiener process with respect to each $\G^k$,
$\E^\Q[|y_0|^p]<\infty$
and $\gamma^k$ is $\G^k$-progressive for every $k$.
If $\tilde\rho (\gamma^k, \gamma^0)\to 0$
as $k\to\infty$,
then
$$
\E^\Q\sup_{t\in [0,T]}|Y^k_t-Y^0_t|^p\to 0,
\quad
\E^\Q \Big[ \int_0^Tf(t,Y^k_t,\gamma^{k}_t)\,dt+g(Y^k_T)\Big]
\to
\E^\Q \Big[ \int_0^Tf (t,Y_t^0,\gamma_t^{0})\,dt+g(Y^0_T)\Big].
$$
\end{Lemma}

Applying the Lemma to $\beta=\hat W$, $Y^k=\hat X^k$,
$\gamma^k=\hat I^k$
 (for $k\ge 1$) and $Y^0=\hat X^\alpha$, $\gamma^0=\hat \alpha$,
and recalling that $\tilde\rho (\hat I^k, \hat \alpha)<1/k\to 0$ we conclude by \eqref{JR}, \eqref{Ja} that   $J^\Rc(\hat\nu^k)   \to J(\alpha)$ as $k\to\infty$.
Therefore relation  \eqref{desiredineq} is satisfied for this choice of the sequence $(\hat\nu^k)_k$ and for the corresponding  settings $(\hat \Omega, \hat \calf,\hat \P_k, \hat W, \hat \mu_k)_k$.
This ends the proof of the inequality ${\text{\Large$\upsilon$}}_0$ $\le$ ${\text{\Large$\upsilon$}}_0^\Rc$.
\qed

\section{Representing the value of a classical control problem by a BSDE with constrained jumps}
\label{Sec:separandom}

In this section the assumptions {\bf (A1)} and {\bf (A2)}
are assumed to hold. We will 
consider the   randomized optimal control problem 
 as described in Section
\ref{randomizedformulation}.
So we choose
a setting $(  \Omega,   \calf, \P,      W,   \mu )$ as described there and let
 $\F^{W,\mu}$ denote the corresponding 
 completed filtration generated by $W$ and $\mu$. 
Also recall the definition of 
  the randomized process $X$ as the solution to equation \eqref{dynXrandom} and the
gain functional
\begin{eqnarray}  \label{separ}
{\text{\Large$\upsilon$}}_0^\Rc &=& \sup_{\nu\in\Vc} \,\E^\nu \Big[ \int_0^T  f(t,X_t, I_t)\, dt + g(X_T) \Big].
\end{eqnarray}

Our purpose is to show that one can introduce a suitable backward SDE and prove that 
${\text{\Large$\upsilon$}}_0^\Rc$ can be represented in terms of its solution. 
In view of 
  Theorem \ref{MainThm}, this will give a BSDE representation of the value 
${\text{\Large$\upsilon$}}_0$  
  of the original control problem as well.
  
   The equation that we are going to introduce 
will be called the \emph{randomized BSDE}. It is  the following constrained BSDE
on the time interval $[0,T]$:
\begin{equation}\label{BSDEconstrained}
\begin{cases}
\vspace{2mm} \dis Y_t \ = \ g(X_T)  + \int_t^T  f(s,X_s ,I_s) ds + K_T - K_t - \int_t^TZ_s\,dW_s - \int_t^T\!\int_A U_s(a)\,\mu(ds\,da), \\
\dis U_t(a) \ \le \ 0.
\end{cases}
\end{equation}

We look for a (minimal) solution to \eqref{BSDEconstrained}  in the sense of the following definition.

\begin{Definition}\label{BSDEdef}
A quadruple $(Y_t,Z_t,U_t(a),K_t)$ $($$t\in [0,T]$, $a\in A$$)$
is called a solution to the BSDE  \eqref{BSDEconstrained} if
\begin{enumerate}
  \item $Y$ $\in$  $\Sc^2(\F^{W,\mu})$, the space of real-valued c\`adl\`ag  $\F^{W,\mu}$-adapted processes  satisfying  
\[
\|Y\|_{\Sc^2}^2:=\E\bigg[\sup_{ t\in [0, T]}|Y_t|^2\bigg]<\infty.
  \]
  \item $Z$ $\in$ $L_W^2(\F^{W,\mu})$, the space of     $\F^{W,\mu}$-progressive  processes
  with values in $\R^d$
  satisfying 
\[
 \|Z\|_{L_W^2}^2 := \E\bigg[\int_0^T|Z_t|^2\, dt\bigg]  < \infty.
\]  

  \item $U$ $\in$ $L_{\lambda}^2(\F^{W,\mu})$, the space of real-valued $\calp(\F^{W,\mu})\otimes \calb(A)$-measurable processes  satisfying 
  \begin{align}
\label{normpoissonlambda}
 \|U\|^2_{L_{\lambda}^2} := \E\bigg[\int_0^T\int_A|U_t(a)|^2\,\lambda(da)\,dt\bigg] < \infty.
  \end{align}

    \item $K$ $\in$ $\Kc^2(\F^{W,\mu})$, the subset of $\Sc^2(\F^{W,\mu})$ consisting of  $\F^{W,\mu}$-predictable nondecreasing process with $K_0=0$.
  \item $\P$-a.s., the equality in \eqref{BSDEconstrained}
holds for every $t\in[0,T]$,
  and the constraint $U_t(a)\le 0$ is understood to hold
  $dt\otimes d\P\otimes d\lambda$-almost everywhere.
\end{enumerate}
A minimal solution $(Y,Z,U,K)$ is a solution to \eqref{BSDEconstrained} such that for any other solution $(Y',Z'$, $U',K')$ we have,
$\P$-a.s., $Y_t\le Y'_t$ for all $t\in [0,T]$.
\end{Definition}

We note that  the space $\Sc^2(\F^{W,\mu})$   becomes a Banach space, with the indicated norm, provided
  we identify elements which are indistinguishable.
Similarly, $L_W^2(\F^{W,\mu})$ (respectively, $L_{\lambda}^2(\F^{W,\mu})$) becomes a Hilbert spaces, with the natural inner product, provided we identify elements that are equal  $dt\otimes d\P$-a.e.  (respectively, 
  $dt\otimes d\P\otimes d\lambda$-a.e.).
These identifications will be assumed to hold in the sequel.
  
We now state the main results of this section.

\begin{Theorem}
\label{Thm:RandomizedFormula}
There exists a unique minimal solution $(Y,Z,U,K)$ $\in$ $\Sc^2(\F^{W,\mu})\times L_W^2(\F^{W,\mu})\times L_{\lambda}^2(\F^{W,\mu})\times\Kc^2(\F^{W,\mu})$ to the randomized equation \eqref{BSDEconstrained}. 
\end{Theorem}

\begin{Theorem}
\label{Thm:RandomizedFormulabis}
For the minimal solution $(Y,Z,U,K)$   to the randomized equation \eqref{BSDEconstrained}  we have the following representation formula: for every $t\in [0,T]$,
\begin{eqnarray} 
\label{RandomizedFormula}
Y_t &=& \esssup_{\nu\in\calv} \E^\nu\bigg[ \int_t^T   f (s,X_s,I_s )\,ds + g(X_T) \,\bigg|\, \calf_t^{W,\mu}\bigg].
\end{eqnarray} 
In particular
\begin{align}\label{RandomizedFormuladet}
    Y_0=\sup_{\nu\in\calv} J^\Rc(\nu)=: {\text{\Large$\upsilon$}}_0^\Rc,
\end{align}
the value of the randomized control problem.
\end{Theorem}

\begin{Remark}
{\rm
Combining Theorems \ref{MainThm} and \ref{Thm:RandomizedFormulabis} we deduce the BSDE representation for the original control problem
\begin{equation}\label{Feynman-Kac}
Y_0 \ = \ \sup_{\alpha\in\Ac^W} J(\alpha)=: {\text{\Large$\upsilon$}}_0.
\end{equation}
Recalling that $Y$ is a component of the minimal solution, we obtain a sort of duality relation:
\[
\sup_{\nu\in\calv} J^\Rc(\nu)=
\inf\{Y_0'\,\colon\, (Y',Z',U',K') \;{\rm solution \;to}\; \eqref{BSDEconstrained}\}.
\]
\qed
}
\end{Remark}

\begin{Remark}
{\rm
The assumption that $\lambda(da)$ has full support, contained in {\bf (A2)}, is not needed in Theorems 
\ref{Thm:RandomizedFormula} and
\ref{Thm:RandomizedFormulabis}. It is however required to arrive at the conclusion \eqref{Feynman-Kac}.
\qed
}
\end{Remark}

BSDEs of this form were first considered in \cite{KMPZ10}, where the generator $f$ was allowed to depend on $Z$ and $U$ as well. In \cite{KP12} some more general case was addressed. In these papers a deep result by Peng (\cite{Peng99})
on monotonic limits of BSDEs had to be used, as well as its extension \cite{Essaky08} to equations with jumps. In \cite{FP15}, for our special case needed for optimal control, some simplification were introduced: below we follow this reference, adding more details to the proof.

\vspace{3mm}

We proceed to the proof of 
 Theorem
\ref{Thm:RandomizedFormula}. The uniqueness part is stated in the following proposition.

\begin{Proposition}
\label{Propuniq:RandomizedFormula}
There exists at most one  minimal solution $(Y,Z,U,K)$   to the randomized equation \eqref{BSDEconstrained}. 
\end{Proposition}

\noindent {\bf Proof.}  Suppose that
$(Y',Z',U',K')$ is another minimal solution. We have $Y=Y'$ up to indistinguishability, by the definition of minimality. Subtracting the corresponding equations and rearranging terms we obtain
\begin{align}
    \label{partiprev}
    \int_0^t(Z_s-Z'_s)\,dW_s=
  K_t -K'_t-  \int_0^t\int_A (U_s(a)-U_s'(a))\,\mu(ds\,da).
\end{align}
Computing the quadratic variations and noting that the right-hand side has finite variation paths we obtain 
$\int_0^t|Z_s-Z'_s|^2\,ds=0$, $\P$-a.s. and so $Z=Z'$, $ d\P\otimes dt$-a.e. It follows that
\[
K_t - K'_t =
    \int_0^t\int_A (U_s(a)-U_s'(a))\,\mu(ds\,da)=
    \sum_{T_n\le t} \Big(U_{T_n}(A_n)-U'_{T_n}(A_n)\Big),
\]
where $(T_n,A_n)$ is the marked point process associated to $\mu$.
The left-hand side is a predictable process, so it has no totally inaccessible  jumps times (see e.g. 
\cite{JS03} Chapter I Proposition 2.24) while the right-hand side is a pure jump process with totally inaccessible jump times. It follows that we have $U_{T_n}(A_n)=U'_{T_n}(A_n)$, $\P$-a.s. for all $n$, and so $K=K'$. Finally,
\[
\E \int_0^T\int_A |U_s(a)-U_s'(a)|\,\lambda(da)\,ds=
 \E\int_0^T\int_A |U_s(a)-U_s'(a)|\,\mu(ds\,da)=0
\]
where the first equality follows from by predictability of the integrand and from
Proposition \ref{propcompensator}-(vii).
We conclude that $U=U'$, $dt\otimes d\P\otimes\lambda(da)$-almost everywhere as required.
\qed

The existence part in the proof of 
 Theorem
\ref{Thm:RandomizedFormula} needs a series of Lemmas. We first
 introduce a family of penalized BSDEs associated to \eqref{BSDEconstrained}, parametrized by  an integer $n$ $\geq$ $1$:
\begin{eqnarray} \label{bsdepenalized}
Y_t^n &=&  g(X_T) + \int_t^T f(s,X_s,I_s)\, ds  + n \int_t^T \int_A U_s^n(a)^+ \lambda(da) \,ds \\
 & &  -  \int_t^T Z^n_s\,dW_s  - \int_t^T\int_AU^n_s(a)\,\mu (ds\,da), \;\;\;   t\in [0,T],  \nonumber
\end{eqnarray} 
where $u^+$ $=$ $\max(u,0)$.

\begin{Lemma}\label{lemmapenalized}
There exists a unique solution
$(Y^n,Z^n,U^n)\in  \Sc^2(\F^{W,\mu})\times L_W^2(\F^{W,\mu})\times L_{\lambda}^2(\F^{W,\mu})$
to the penalized BSDE \eqref{bsdepenalized}.   
\end{Lemma}

\noindent {\bf Proof.}
We only sketch the proof and we refer the reader to  \cite{LiTang94}  Lemma 2.4 or to the paper \cite{BarBuckPar97}
for full details.

Let us first consider the linear equation
\begin{eqnarray} %\label{bsdepenalized}
Y_t &=&  g(X) + \int_t^T f(s,X_s,I_s)\, ds 
%\\& & 
-  \int_t^T Z_s\,dW_s  - \int_t^T\int_AU_s(a)\,\mu (ds\,da), \;\;\;   t\in [0,T],   \nonumber
\end{eqnarray} 
with unknown $(Y,Z,U)$. 
Let us define 
\begin{eqnarray} %\label{bsdepenalized}
Y_t &=&  \E\left[ g(X) + \int_t^T f(s,X_s,I_s)\, ds \,\bigg| \, \calf^{W,\mu}_t\right]
=M_t-\int_0^t f(s,X_s,I_s)\, ds,\nonumber
\end{eqnarray} 
where
\[
M_t= 
\E\left[ g(X) + \int_0^T f(s,X_s,I_s)\, ds \,\bigg| \, \calf^{W,\mu}_t\right],\;\;\;   t\in [0,T].
\]
It is easy to check that under our assumptions the process $M$ is a square integrable martingale with respect to $\F^{W,\mu}$. By a martingale representation theorem for the filtration generated by a Wiener process and an independent Poisson random measure (see \cite{LiTang94}  Lemma 2.4, or \cite{ja} Theorem 5.4 for a more general case), there exist 
$(Z,U)\in   L_W^2(\F^{W,\mu})\times L_{\lambda}^2(\F^{W,\mu})$ such that
\[
M_t=   \int_0^t Z_s\,dW_s  + \int_0^t\int_AU_s(a)\,\mu (ds\,da),
\]
and  $(Y,Z,U)$ is the unique required solution.

In the general case the proof is finished by a standard Picard iteration technique, due to the Lipschitz character of the additional term 
$  n \int_t^T \int_A U_s(a)^+ \lambda(da)\, ds$.
\qed

We provide an explicit representation of the solution to the penalized BSDE in terms of a family of auxiliary randomized control problems. For every integer $n\ge1$,
let $\calv^n $ denote the subset of elements $\nu_t(\omega,a)$ in $\calv $ valued in $(0,n]$.

\begin{Lemma}\label{bsdepenlizaedandcontrol}
We have for all $n$ $\geq$ $1$ and $t \in [0,T]$,
\begin{equation} \label{expliYn}
Y_t^n  =  \mathop{\rm ess\, sup}_{\nu\in\calv^n}  \E^\nu \Big[ \int_t^T f (s,X_s,I_s) \,ds +  g(X_T) \big| \calf_t \Big],  \qquad \P-a.s.
\end{equation}
\end{Lemma}
{\bf Proof.} Fix $n$ $\geq$ $1$ and consider $(Y^n,Z^n,U^n)$   solution to \eqref{bsdepenalized}.
We first prove that 
the processes
\begin{eqnarray*}
\int_0^t Z^n_s dW_s, & & \int_0^t \int_A U_s^n(a) \mu (ds\,da)-   \int_0^t \int_A U_s^n(a)\,\nu_s(a)\, \lambda (da)\,ds
,\qquad t\in [0,T]
\end{eqnarray*}
are $\P^\nu$-martingales. 

Indeed, $W$   remains a Wiener process under $\P^\nu$, so $\int_0^\cdot Z^n_s dW_s$ is a $\P^\nu$-local martingale; to show that it is a martingale, by the Burkholder-Davis-Gundy inequality it is enough to note that
\[
\E^\nu\left[ \left(\int_0^T|Z^n_s|^2\,ds\right)^\frac{1}{2}\right]=
\E\left[\kappa_T^\nu \left(\int_0^T|Z^n_s|^2\,ds\right)^\frac{1}{2}\right]
\le
(\E\,|\kappa_T^\nu|^2)^\frac{1}{2}
\left(\E \int_0^T|Z^n_s|^2\,ds\right)^\frac{1}{2} <\infty,
\]
since 
$Z^n$ $\in$ $L_W^2(\F^{W,\mu}) $ and
$\E\,|\kappa_T^\nu|^2<\infty$
by 
 Proposition \ref{girsanovpoisson}-(i).
Concerning the other process we note that
\[
\E\left[\int_0^t \int_A |U_s^n(a)|\,\nu_t(a)\, \lambda (da)\,ds\right]<\infty,
\]
which follows from $U^n$ $\in$ $L_{\lambda}^2(\F^{W,\mu})$ and the boundedness of $\nu$, and we conclude recalling 
Proposition \ref{propcompensator}-(vi).

 From the established martingale property, by taking the conditional expectation $\E^\nu$ given $\calf_t$ in  \eqref{bsdepenalized}, we obtain, for every $t \in [0,T]$:
\begin{eqnarray*}
 Y_t^n &=& \E^\nu \Big[ \int_t^T f(s,X_s,I_s)\,ds + g(X_T) \big| \calf_t\Big]  \\
 & & + \; \E^\nu\Big[ \int_t^T\int_A [n  U^n_s(a)^ + -  \nu_s(a) U^n_s(a) ] \lambda(da) \,ds   \big| \calf_t \Big], \;\;\; \P-a.s.
\end{eqnarray*}
From the elementary numerical inequality: $nu^+-\nu u$ $\geq$ $0$ for all $u$ $\in$ $\R$, $\nu$ $\in$ $[0,n]$, we deduce that
\begin{eqnarray*}
 Y_t^n & \geq & \mathop{\rm ess\, sup}_{\nu\in\calv^n}  \E^\nu \Big[ \int_t^T f (s,X_s,I_s)\,ds + g(X_T) \big| \calf_t\Big],   \;\;\;\P-a.s.
\end{eqnarray*}
To prove the opposite inequality it would be enough to  take $\nu_s(a)$ $=$ $n 1_{\{U^n_s(a)\ge 0\}}$; however,   this random field does not belong to $\calv^n$
because of the requirement of strict positivity of $\nu$ in the definition of $\calv$ and $\calv^n$. Therefore, we proceed taking an approximation:
for  $\epsilon$ $\in$ $(0,1)$,  define
$$\nu^\epsilon_s(a)= n1_{\{U^n_s(a)\ge 0\}} + \epsilon 1_{\{-1<U^n_s(a)< 0\}}
 - \epsilon U^n_s(a)^{-1} 1_{\{ U^n_s(a)\le-1\}}.
 $$
Then  $\nu^\epsilon$ $\in$  $\calv^n $  and we have
\begin{eqnarray*}
n U^n_s(a)^+  -  \nu_s^\epsilon(a) U^n_s(a) &\le & \epsilon, \;\;\;   s\in [0,  T],
\end{eqnarray*}
so that
\begin{eqnarray*}
Y_t^n &\le& \E^{\nu^\epsilon} \Big[ \int_t^T f (s,X_s,I_s)\,ds + g(X_T) \big| \calf_t\Big] + \epsilon T\lambda(A)  \\
 &\le& \mathop{\rm ess\, sup}_{\nu\in\calv^n }  \E^\nu \Big[ \int_t^Tf (s,X_s,I_s) \,ds + g(X_T) \big| \calf_t\Big]+\epsilon T\lambda(A),
\end{eqnarray*}
which is enough to conclude the proof. 
\qed

\vspace{3mm}

As a consequence of this explicit representation of  the penalized BSDE, we  obtain the follo\-wing uniform estimate on the sequence $(Y^n)_n$:

\begin{Lemma}\label{estimatesonbsdepenlizaed}
The sequence $(Y^n)_n$ is monotonically increasing in $n$: $\P$-a.s.,
\begin{equation}
    \label{ynmonotono}
    Y^n_t\le Y^{n+1}_t, \qquad t\in [0,T].
\end{equation}
Morover we have
\begin{eqnarray} \label{ynbounded}
\sup_{t\in [0,T]}|Y_t^n| &\le & C  \Big(1+\sup_{t\in [0,T]}|X_t|^{\max(r,2)}\Big),  \qquad \P-a.s.
\\\label{ynboundedmean}
\E\,\bigg[\sup_{t\in[0,T]}|Y_t^n|^2\bigg]&\le&  C\,( 1 +  |x_0 |^{\max(r,2)}  ).
\end{eqnarray} 
 for some constant $C$ depending only on  $T$ and on the constants $r$ and $L$ as defined in {\bf (A1)}. 
 \end{Lemma}
 
 \noindent
{\bf Proof.} Monotonicity follows from the formula for $Y^n$ presented in Lemma  \ref{bsdepenlizaedandcontrol}, since $\calv^n $ $\subset$ $\calv^{n+1} $.
Then the inequality $Y^n_t$ $\le$ $Y^{n+1}_t$ holds $\P$-a.s. for all $t\in[0,T]$ since these processes are c\`adl\`ag.

   Recall the polynomial growth condition \eqref{PolGrowth_f_g} on $f$ and $g$ in {\bf (A1)} and
  denote by $C$ a generic constant  depending only on  $T,r,L$, whose precise value may possibly change at each occurrence.
It follows from  Lemma \ref{bsdepenlizaedandcontrol} that for every $t \in [0,T]$,
\begin{eqnarray*}
|Y_t^n| &\le&  C \mathop{\rm ess\, sup}_{\nu\in\calv^n }  \E^\nu \Big[1+ \sup_{s\in [t,T]}|X_s|^r  \Big| \calf_t\Big], \qquad  \P-a.s.
\end{eqnarray*}
Set $p=\max(r,2)$ and recall the conditional estimate \eqref{EstimateX_nu_cond}
that now gives 
\begin{eqnarray*}
 \E^\nu \Big[ \sup_{s\in [t,T]}|X_s|^p \Big| \calf_t\Big] &\le & C \Big( 1+ |X_t|^p \Big),
\qquad\P-a.s.
\end{eqnarray*}
 where the constant $C$ can be chosen to be the same for every $\nu\in \calv $.
It follows that $|Y_t^n| \le C\big(1+ |X_t|^p\big)$,  $\P$-a.s.  and the inequality \eqref{ynbounded}  follows immediately. The inequality \eqref{ynboundedmean} is now a consequence of estimate
 \eqref{EstimateX}. 
\qed

 In the following lemma 
similar uniform estimates are also proved for  $(Z^n,U^n,K^n)_n$:
 
\begin{Lemma} \label{boundsuznun}
Consider the  solution
$(Y^n,Z^n,U^n)$
to \eqref{bsdepenalized} and 
let $B$ denote a constant such that
$\E\,[\sup_{s\in[0,T]}|Y_s^n|^2]\le B$ for all $n$.  Then
\begin{eqnarray} \nonumber
\|Y^n\|^2_{\Sc^2(\F^{W,\mu})} + \|Z^n\|^2_{L_W^2(\F^{W,\mu})}
+\|U^n\|^2_{L_{\lambda}^2(\F^{W,\mu})} +\|K^n\|^2_{\Sc^2(\F^{W,\mu})}
\\\qquad \le C\left(B+\E\,|g(X_T)|^2+\E\int_0^T|f(t,X_t,I_t)|^2\,dt\right),
\end{eqnarray} 
where $K_t^n= n \int_0^t \int_A U_s^n(a)^+ \lambda(da) \,ds$, $t\in [0,T]$.
\end{Lemma}

\noindent {\bf Proof.} We note that by
the previous Lemma  
we may choose   $B=  C\,( 1 +  |x_0 |^{\max(r,2)}  )$ and the bound on  $
\|Y^n\|^2_{\Sc^2(\F^{W,\mu})}$ is already established. 
Then   we write
\eqref{bsdepenalized} in the form
\begin{eqnarray} \label{penalizzconkappa}
Y_t^n &=&  g(X_T) + \int_t^T f(s,X_s,I_s)\, ds  +K_T^n-K_t^n 
\\
 & &  -  \int_t^T Z^n_s\,dW_s  - \int_t^T\int_AU^n_s(a)\,\mu (ds\,da), \;\;\;   t\in [0,T],  \nonumber
\end{eqnarray} 
Then we apply the It\^o formula to $|Y^n_t|^2$:
\[
|g(X_T)|^2- |Y^n_t|^2=2\int_{(t,T]}Y^n_{s-}\,dY^n_s +
\int_t^T|Z^n_s|^2\,ds + \sum_{t<s\le T} |\Delta Y^n_s|^2,
\]
where $\Delta Y^n_t:=Y^n_t-Y^n_{t-}=
\int_A U_t^n(a)\, \mu(\{t\},da)$.  By the discrete character of the measure $\mu$ we have
\[
 \sum_{t<s\le T} |\Delta Y^n_s|^2= 
\int_{(t,T]}\int_A |U_s^n(a)|^2\,\mu(ds\,da).
\]
It follows that
\begin{align}\nonumber
|g(X_T)|^2&= |Y^n_t|^2
-2\int_{t}^TY^n_{s}f(s,X_s,I_s)\,ds 
-2\int_{(t,T]}Y^n_{s-}\,dK^n_s +
\int_t^TY^n_{s-}Z^n_s \,dW_s 
\\
&+\int_t^T|Z^n_s|^2\,ds +\int_{(t,T]}\int_A (2Y^n_{s-}U^n_s(a)+|U^n_s(a)|^2)\,\mu(ds\,da).\label{itopenalizz}
\end{align}
The process
$\int_0^\cdot Y^n_{s-} Z^n_s\, dW_s$ is a  local martingale; to show that it is a martingale, and in particular that it has zero expectation, by the Burkholder-Davis-Gundy and the Cauchy-Schwarz inequality it is enough to note that
\[
\E\left[ \left(\int_0^T|Y^n_{s-}Z^n_s|^2\,ds\right)^\frac{1}{2}\right]
\le
(\E\,\sup_{s\in [0,T]}|Y^n_s|^2)^\frac{1}{2}
\left(\E \int_0^T|Z^n_s|^2\,ds\right)^\frac{1}{2} <\infty,
\]
since $Y^n$ $\in$ $\cals^2(\F^{W,\mu}) $ and
$Z^n$ $\in$ $L_W^2(\F^{W,\mu}) $.
Similarly, the process
\[\int_{(0,t]}\int_A (2Y^n_{s-}U^n_s(a)+|U^n_s(a)|^2)\,\mu(ds\,da)
-
\int_{(0,t]}\int_A (2Y^n_{s-}U^n_s(a)+|U^n_s(a)|^2)\,\lambda(da)\,ds, \quad t\in [0,T],
\]
is a  martingale, by 
Proposition \ref{propcompensator}-(vi), since the integrand process is $\F^{W,\mu}$-predictable (the integrability condition required in 
Proposition \ref{propcompensator} 
follows from 
$Y^n$ $\in$ $S^2(\F^{W,\mu}) $ and $U^n$ $\in$ $L_\lambda^2(\F^{W,\mu}) $);  in particular
\[
\E \left[\int_{(t,T]}\int_A (2Y^n_{s-}U^n_s(a)+|U^n_s(a)|^2)\,\mu(ds\,da)\right]=\E\left[ 
\int_{t}^T\int_A (2Y^n_{s-}U^n_s(a)+|U^n_s(a)|^2)\,\lambda(da)\,ds\right].
\]
Taking expectation in \eqref{itopenalizz} and rearranging we obtain
\begin{align*}&
\E\,|Y^n_t|^2 +\E\int_t^T|Z^n_s|^2\,ds 
+\E  
\int_{t}^T\int_A  |U^n_s(a)|^2 \,\lambda(da)\,ds 
\\
\qquad
&=\E\, |g(X_T)|^2
-2\E\int_{t}^TY^n_{s}f(s,X_s,I_s)\,ds 
-2\E\int_{(t,T]}Y^n_{s-}\,dK^n_s 
&+2\E\int_{t}^T\int_A Y_{s}^nU^n_s(a) \,\lambda(da)\,ds. 
\end{align*}

Let $B$ denote a constant such that
$\E\,[\sup_{s\in[0,T]}|Y_s^n|^2]\le B$ for all $n$. Then, for $\alpha>0$,
\[
2\E\left|\int_{(t,T]}Y^n_{s-}\,dK^n_s\right|\le
2\E\left[\sup_{s\in[0,T]}|Y_s|\,|K_T^n-K^n_t|\right]\le \alpha\, \E\left[ |K_T^n-K^n_t|^2\right]
+
\frac{1}{\alpha} \,B,
\]
where for the last passage we have used the numerical inequality $2ab\le \alpha a^2 + b^2/\alpha$. 
Using the same inequality we have, for every $\beta>0$,
\begin{align*}
\E\,|Y^n_t|^2 +\E\int_t^T|Z^n_s|^2\,ds 
+\E  
\int_{t}^T\int_A  |U^n_s(a)|^2 \,\lambda(da)\,ds \le \E\, |g(X_T)|^2+ \alpha\, \E\left[ |K_T^n-K^n_t|^2\right]+
\frac{1}{\alpha} \,B
\\
\;
+
\E\int_{t}^T|Y^n_{s}|^2\,ds 
 +\E\int_{t}^T|f(s,X_s,I_s)|^2\,ds 
+\frac{1}{\beta}  
\E\int_{t}^T|Y^n_{s}|^2\,ds +
\beta \lambda(A)\E\int_{t}^T\int_A |U^n_s(a) |^2\,\lambda(da)\,ds. 
\end{align*}
We can estimate 
$\E\left[ |K_T^n-K^n_t|^2\right]$
from \eqref{penalizzconkappa} obtaining
\begin{eqnarray} \label{stimakappan}
\E\left[ |K_T^n-K^n_t|^2\right]\le B'\bigg(
\E\,|Y_t^n|^2 +  \E\,| g(X_T) |^2 + \E\int_t^T |f(s,X_s,I_s)|^2\, ds   
\\
+  \E \int_t^T |Z^n_s|^2\,ds  + \E\int_t^T\int_A|U^n_s(a)|^2\,\lambda (da)\,ds, \bigg)\;\;\;   t\in [0,T],  \nonumber
\end{eqnarray} 
for some constant $B'$ that  depends only on $T$ and $\lambda(A)$. 
Substituting in the previous inequality we obtain
\begin{align*}
(1-\alpha B')\,\E\,|Y^n_t|^2 +(1-\alpha B')\, \E\int_t^T|Z^n_s|^2\,ds 
+(1-\beta \lambda(A)-\alpha B')\,\E  
\int_{t}^T\int_A  |U^n_s(a)|^2 \,\lambda(da)\,ds 
\\
\le (1+\alpha B')\,\E\, |g(X_T)|^2+ 
\frac{1}{\alpha} \,B
+\left(1+ \frac{1}{\beta} \right)\,
\E\int_{t}^T|Y^n_{s}|^2\,ds 
+(1+\alpha B')\,
 \E\int_{t}^T|f(s,X_s,I_s)|^2\,ds .
\end{align*}
Next we fix $\alpha,\beta>0$ so small that all the terms in round brackets are strictly positive and we apply the Gronwall Lemma to $t\mapsto \E\,|Y^n_t|^2$ obtaining
\begin{align*}
\sup_{t\in [0,T]}\,\E\,|Y^n_t|^2 + \E\int_t^T|Z^n_s|^2\,ds 
+\E  
\int_{t}^T\int_A  |U^n_s(a)|^2 \,\lambda(da)\,ds 
\\
\le C\left(B+\E\, |g(X_T)|^2+
 \E\int_{t}^T|f(s,X_s,I_s)|^2\,ds \right).
\end{align*}
This proves the required bound on
$ \|Z^n\|^2_{L_W^2(\F^{W,\mu})}$ and $\|U^n\|^2_{L_{\lambda}^2(\F^{W,\mu})}$. The required bound on 
$\|K^n\|^2_{\cals^2(\F^{W,\mu})}=\E\left[ |K_T^n|^2\right]$
now follows from
\eqref{stimakappan} and   the proof is finished. 
\qed

In order to conclude the  proof of Theorem \ref{Thm:RandomizedFormula} we also record the following simple technical lemma.

\begin{Lemma} \label{convdeboleint}For every $\F^{W,\mu}$-stopping time $\tau$ with values in $[0,T]$
    the mappings
    \begin{align}\label{intcontforte}
      &L_\lambda^2(\F^{W,\mu})\to L^2(\Omega, \calf_\tau^{W,\mu},\P) ,
      \qquad
      &&U\mapsto \int_0^\tau\int_AU_t(a)\,\mu(dt\,da),
      \\ \label{intcontfortedue}
    &L_W^2(\F^{W,\mu})\to L^2(\Omega, \calf_\tau^{W,\mu},\P) ,
      \qquad
      &&Z\mapsto \int_0^\tau Z_s\,dW_t, 
    \end{align}
    are continuous and also weakly continuous.
\end{Lemma}

\noindent {\bf Proof.}
Both maps being linear it is enough to prove strong continuity. Suppose that $ U \in L_\lambda^2(\F^{W,\mu})$. Then
\begin{align*}
\E\left| \int_0^\tau \int_A  U_s(a)\,\mu(da\,dt)\right|^2
&=
\E\left| \int_0^T \int_A  U_s(a)\,1_{\{t\le\tau\}}\mu(da\,dt)\right|^2
\\&=
\E \int_0^T \int_A | U_s(a)|^2\,1_{\{t\le \tau\}}|\,\lambda(da)\,dt 
 \le \|U\|^2_{L_\lambda^2(\F^{W,\mu})}<\infty,
\end{align*}
where the second equality follows from predictability of the integrand and from a form of the It\^o isometry (see e.g. \cite{bremaud2} Theorem 5.1.33). This shows that
  the map \eqref{intcontforte} is well defined. Suppose now that $U^k\to U$ in $L_\lambda^2(\F^{W,\mu})$. Then, by the same arguments, 
\begin{align*}
\E\left| \int_0^T \int_A (U^k_s(a)- U_s(a))\,1_{\{t\le\tau\}}\mu(da\,dt)\right|^2
&=
\E \int_0^T \int_A |U^k_s(a)- U_s(a)|^2\,1_{\{t\le \tau\}}\lambda(da)\,dt 
\\& \le \|U^k-U\|^2_{L_\lambda^2(\F^{W,\mu})}\to 0.
\end{align*}
The well posedness and strong continuity of \eqref{intcontfortedue} is proved similarly.
\qed

We are now in a position to complete the proof of Theorem \ref{Thm:RandomizedFormula}.

\vspace{3mm}

\noindent {\bf End of the proof of Theorem \ref{Thm:RandomizedFormula}: existence}.
 Let us 
consider again the  solution
$(Y^n,Z^n,U^n)$
to the penalized equation \eqref{bsdepenalized}.
By the monotonicity property \eqref{ynmonotono}, with the exception of a $\P$-null set the limit
\[
Y_t=\lim_{n\to\infty}Y^n_t, \qquad t\in [0,T],
\]
exists, so that the process $Y$ is defined up to indistinguishability and is $\F^{W,\mu}$-progressive. 
From \eqref{ynboundedmean} and the Fatou lemma we have $ 
\E\,[\sup_{t\in[0,T]}|Y_t|^2 ]<\infty$ so that, by the dominated convergence theorem,
\begin{eqnarray} \label{convyennetau}
\E\int_0^T |Y_t^n-Y_t|^2\,dt\to 0,
\qquad
\E |Y_\tau^n-Y_\tau|^2\,dt\to 0,
\end{eqnarray} 
for every stopping time $\tau$ with values in $[0,T]$.  So far, it is not yet clear that $Y$ is c\`adl\`ag.

Define the increasing continuous process   
$K_t^n= n \int_0^t \int_A U_s^n(a)^+ \lambda(da) \,ds$ as before and recall
\eqref{penalizzconkappa} that we write in the form
\begin{eqnarray} \label{kappaenne}
K_t^n= -Y_t^n+Y_0^n - \int_0^t f(s,X_s,I_s)\, ds  + \int_0^t Z^n_s\,dW_s  +\int_0^t\int_AU^n_s(a)\,\mu (ds\,da), \qquad   t\in [0,T]. \nonumber
\end{eqnarray} 
By  Lemma  \ref{boundsuznun} the sequence 
$(Z^n,U^n)$ is bounded in the Hilbert space $ {\Sc^2(\F^{W,\mu})}\times{L_W^2(\F^{W,\mu})}$ and therefore there exists a subsequence, denoted $(Z^{n'},U^{n'})$, weakly convergent to a limit $(Z,U)$ in this space.
Define an $\F^{W,\mu}$-optional process $K$ setting
\begin{eqnarray} \label{kappalimite}
K_t= -Y_t+Y_0 - \int_0^t f(s,X_s,I_s)\, ds  + \int_0^t Z_s\,dW_s  +\int_0^t\int_AU_s(a)\,\mu (ds\,da), \qquad   t\in [0,T].
\end{eqnarray} 
%Next we claim that $K$ is increasing, namely that  $\P$-a.s., 
%\begin{eqnarray} \label{Kcrescenteqc}
%K_s\le K_t, \qquad {\rm for \;all}
%\qquad 0\le s\le t\le T. 
%\end{eqnarray} 
%To see this 
Note the for every $[0,T]$-valued stopping time $\tau$ we have
\begin{eqnarray} \label{kappaennetau}
K_\tau^n= -Y_\tau^n+Y_0^n - \int_0^\tau f(s,X_s,I_s)\, ds  + \int_0^\tau Z^n_s\,dW_s  +\int_0^\tau\int_AU^n_s(a)\,\mu (ds\,da), \qquad  \P-a.s. \nonumber
\\
K_\tau= -Y_\tau+Y_0 - \int_0^\tau f(s,X_s,I_s)\, ds  + \int_0^\tau Z_s\,dW_s  +\int_0^\tau\int_AU_s(a)\,\mu (ds\,da), \qquad  \P-a.s. \nonumber
\end{eqnarray} 
and from Lemma \ref{convdeboleint} and
\eqref{convyennetau}
it follows that 
$K_\tau^{n'}\to K_\tau$ weakly in 
$L^2(\Omega, \calf_\tau^{W,\mu},\P)$.

It is also easy to see that $K^{n'}\to K$ weakly in $L_W^2(\F^{W,\mu})$. The latter is a space of equivalence classes. We say that a class is predictable if it contains a representative which is $\F^{W,\mu}$-predictable (equivalently, $\calp(\F^{W,\mu})$-measurable). It is easily verified that $\F^{W,\mu}$-predictable elements of $L_W^2(\F^{W,\mu})$ form a closed subspace, hence also a weakly closed set. 
Since $K^n$ are continuous adapted they are predictable and so the limit process $K$ is predictable as well.

Moreover, if $\tau$ and $\sigma$ are two stopping times satisfying $\sigma\le \tau\le T$ then, recalling that $K^n$ is increasing, we have $K^{n'}_\sigma\le K^{n'}_\tau$ $\P$-a.s. and therefore, passing to the limit, also
$K_\sigma\le K_\tau$ $\P$-a.s. This implies that the limit process is increasing as well: we refer the reader to \cite{BaCoCo19}
Theorem 5.6, pages 23-24,
for a detailed verification of this fact.

From \eqref{kappalimite} we see that $Y$, defined as a monotonic limit of the c\`adl\`ag  processes $Y^n$, equals a c\`adl\`ag process (given by the three stochastic integrals) minus the increasing process $K$. By a real-analytic lemma, it follows that $Y$ also has c\`adl\`ag trajectories, and so $K$ as well. We refer to \cite{Peng99} Lemma 2.2 for this lemma, which is short and elementary, but nontrivial.

Since $Y^n_T=g(X_T)$ we also have $Y_T=g(X_T)$ and it is easy to see that \eqref{kappalimite} can be rewritten as our equation \eqref{BSDEconstrained}. 

Next we prove that the constraint $U_t(a)\le 0$ is verified. To this end consider the functional $G:L^2_\lambda(\F^{W,\mu})\to \R$ defined by
\[
G(U)=   \E\left[\int_0^T \int_A U_s(a)^+ \lambda(da) \,ds\right],
\qquad U\in L^2_\lambda(\F^{W,\mu}). 
\]
%%%%%%%%%%%%%%%%%
Since $G$ is a convex functional on the Hilbert space $L^2_\lambda(\F^{W,\mu})$ and it is continuous
with respect to the norm introduced in 
\eqref{normpoissonlambda}, it is also weakly continuous.
%%%%%%%%%%%%%%%%%
 Since $n\cdot G(U^n)=\E\left[ K_T^n\right] \le \E\left[ |K_T^n|^2\right]^{1/2}=\|K^n\|_{\cals^2(\F^{W,\mu})}$ is bounded by Lemma \ref{boundsuznun}, we have $G(U^n)\to 0$.  Since we know that $U^{n'}\to U$ weakly, it follows that $G(U)\le\liminf_{n'} G(U^{n'})=0$ and so $G(U)=0$, which is equivalent to $U_t(\omega,a)= 0$,  $\P(d\omega)\lambda(da)dt$-almost surely.

 We have now proved that our tuple $(Y,Z,U,K)$ is a solution to $\eqref{BSDEconstrained}$. To prove that it is the minimal solution, suppose that $(Y',Z',U',K')$ is another solution. Since $K'$ is increasing and $U'\le 0$, removing two nonnegative terms in the equation we obtain
\[ Y'_t \ge  g(X_T)  + \int_t^T  f(s,X_s ,I_s) ds - \int_t^TZ_s\,dW_s.
\]
 For any $n\ge1$ let us take $\nu\in \calv^n =\{\nu\in\calv\,:\, \nu\le n\}$ and consider the corresponding probability $\P^\nu$. Proceeding  as in Lemma 
\ref{bsdepenlizaedandcontrol} one can show   that 
$\int_0^\cdot Z_s dW_s$ is a  $\P^\nu$-martingale, so taking the conditional expectation $\E^\nu$ given $\calf_t$   we obtain:
\begin{align}\label{Y'ineqvalue}
    Y'_t  \ge \E^\nu \Big[ \int_t^T f(s,X_s,I_s)\,ds + g(X_T) \big| \calf_t\Big], \qquad \P-a.s. 
\end{align}
Taking $\mathop{\rm ess\, sup}$ over all $\nu\in\calv^n$ and recalling 
Lemma 
\ref{bsdepenlizaedandcontrol} we conclude  that, for every $t$, 
$Y'_t\ge Y_t^n$ $\P$-a.s. Since these processes are c\`adl\`ag this inequality holds for every $t$, $\P$-a.s. Since our solution $Y$ is the increasing limit of $Y^n$ we conclude that $Y'_t\ge Y_t$ for every $t$, $\P$-a.s. Now minimality is proved and the proof of 
 Theorem
\ref{Thm:RandomizedFormula} is finished.
\qed

\begin{Remark}
    {\rm
    In the previous proof we have shown convergence of   a subsequence $(U^{n'},Z^{n'})$ to $(U,Z)$, weakly in  the space $ {\Sc^2(\F^{W,\mu})}\times{L_W^2(\F^{W,\mu})}$.
It may be useful to know that in fact the entire sequence $(U^n,Z^n)$ converges weakly. Indeed, suppose this is not the case. Then there exists another  
    subsequence $(U^{n''},Z^{n''})$  weakly convergent to a different limit, denoted $(U',Z')$. 
Setting, similar to formula 
\eqref{kappalimite}, 
\begin{eqnarray} \nonumber
K'_t= -Y_t+Y_0 - \int_0^t f(s,X_s,I_s)\, ds  + \int_0^t Z'_s\,dW_s  +\int_0^t\int_AU'_s(a)\,\mu (ds\,da), \qquad   t\in [0,T],
\end{eqnarray} 
we would conclude that 
$(Y,U,Z,K)$ and $(Y,U',Z',K')$ are both minimal solution, so by uniqueness we get the contradiction 
  $(U,Z)=(U',Z')$.
    }\qed
\end{Remark}

\noindent {\bf Proof of 
 Theorem
\ref{Thm:RandomizedFormulabis}.} 
Fix $t\in [0,T]$. During the proof of 
 Theorem
\ref{Thm:RandomizedFormula}, in \eqref{Y'ineqvalue}, we have shown that for 
any solution $(Y',Z',U',K')$
to $\eqref{BSDEconstrained}$ (even non-minimal) we have
\begin{align*}
    Y'_t  \ge \E^\nu \Big[ \int_t^T f(s,X_s,I_s)\,ds + g(X_T) \big| \calf_t\Big], \qquad \P-a.s. 
\end{align*}
where 
 $\nu\in \calv^n $, $\P^\nu$ denotes  the corresponding probability   and  $n\ge1$ is arbitrary.
  Now choose any $\nu\in \calv $. Since the elements of $\calv$ are bounded and
 $  \calv^n =\{\nu\in\calv\,:\, \nu\le n\}$, for large $n$ we have 
 $\nu\in \calv^n $ and we may apply the previous inequality to our minimal solution 
$(Y,Z,U,K)$ obtaining
\begin{align*}
    Y_t  \ge \E^\nu \Big[ \int_t^T f(s,X_s,I_s)\,ds + g(X_T) \big| \calf_t\Big], \qquad \P-a.s. 
\end{align*}
 It follows that $\P$-a.s.,
 \begin{align}
     \label{ytvalue}
 Y_t \ge  \esssup_{\nu\in\calv} \E^\nu\bigg[ \int_t^T   f (s,X_s,I_s )\,ds + g(X_T) \,\bigg|\, \calf_t^{W,\mu}\bigg].
 \end{align}
On the other hand we have, for all $n$ $\geq$ $1$, $\P$-a.s.,
\[
Y_t^n  =  \mathop{\rm ess\, sup}_{\nu\in\calv^n}  \E^\nu \Big[ \int_t^T f (s,X_s,I_s) \,ds +  g(X_T) \big| \calf_t \Big]\le 
\mathop{\rm ess\, sup}_{\nu\in\calv}  \E^\nu \Big[ \int_t^T f (s,X_s,I_s) \,ds +  g(X_T) \big| \calf_t \Big],
\]
where the   equality is the content of 
 Lemma \ref{bsdepenlizaedandcontrol}
and the inequality follows from the inclusion 
$\calv^n $ $\subset$ $\calv $.
Since $Y^n_t\to Y_t$ we obtain  the inequality opposite to \eqref{ytvalue}
and this proves \eqref{RandomizedFormula}. 
Finally, we note that \eqref{RandomizedFormuladet} follows immediately from \eqref{RandomizedFormula} setting $t=0$ since the $\sigma$-algebra $\calf_0^{W,\mu}$ is trivial.
\qed

We end this Section proving the following generalization of formula \eqref{RandomizedFormula}.

\begin{Theorem}\label{T:DPPRandomizedProblem}
\quad\\
For all $t\in[0, T]$ we have
\begin{eqnarray} \label{DPPRandomizedProblem}
Y_t &=& \esssup_{\nu\in\calv} \esssup_{\tau\in\Tc} \E^\nu\bigg[ \int_t^\tau f(r,X_r,I_r) \,dr + Y_\tau \,\bigg|\, \calf_t^{W,\mu}\bigg] \notag \\
&=& \esssup_{\nu\in\calv} \essinf_{\tau\in\Tc} \E^\nu\bigg[ \int_t^\tau f(r,X_r,I_r) \,dr + Y_\tau \,\bigg|\, \calf_t^{W,\mu}\bigg],
\end{eqnarray} 
where $\mathcal T$ denotes the class of $[t,T]$-valued $\F^{W,\mu}$-stopping times.

In particular, setting $t=0$, we have
\begin{eqnarray} \label{DPPRandomizedProblembis}
Y_0 &=& \esssup_{\nu\in\calv} \esssup_{\tau\in\Tc} \E^\nu\bigg[ \int_0^\tau f(r,X_r,I_r) \,dr + Y_\tau\bigg] \notag \\
&=& \esssup_{\nu\in\calv} \essinf_{\tau\in\Tc} \E^\nu\bigg[ \int_0^\tau f(r,X_r,I_r) \,dr + Y_\tau \bigg],
\end{eqnarray} 
where now $\mathcal T$ denotes the class of $[0,T]$-valued $\F^{W,\mu}$-stopping times.
\end{Theorem}

\noindent
\textbf{Proof.}
For every $n$, proceeding along the same lines as in the proof of Lemma \ref{bsdepenlizaedandcontrol}, we obtain
\begin{align}
Y_t^n \ &= \ \E^\nu\bigg[Y_\tau^n + \int_t^\tau f(r,X_r,I_r) \,dr + \int_t^\tau\int_A \big[n\big(U_r^n(a)\big)^+ - \nu_r(a) U_r^n(a)\big]\lambda(da)dr\,\bigg|\, \calf_t^{W,\mu}\bigg] \label{DPPRandomizedProblem4} \\
&\geq \ \E^\nu\bigg[Y_\tau^n + \int_t^\tau f(r,X_r,I_r) \,dr\,\bigg|\, \calf_t^{W,\mu}\bigg], \qquad \text{for all }\nu\in\Vc_n,\,\tau\in\Tc. \notag
\end{align}
Recalling that $Y\geq Y^n$, we find
\[
Y_t \ \geq \ \E^\nu\bigg[Y_\tau^n + \int_t^\tau f(r,X_r,I_r) \,dr\,\bigg|\, \calf_t^{W,\mu}\bigg], \qquad \text{for all }\nu\in\Vc_n,\,\tau\in\Tc.
\]
Noting that $\nu\in\Vc$ belongs to 
$\Vc_n$ for all $n$ sufficiently large, and 
letting $n\rightarrow\infty$,  we end up with
\[
Y_t \ \geq \ \E^\nu\bigg[Y_\tau + \int_t^\tau f(r,X_r,I_r) \,dr\,\bigg|\, \calf_t^{W,\mu}\bigg], \qquad \text{for all }\nu\in\Vc,\,\tau\in\Tc.
\]
The above inequality yields
\begin{align*}
Y_t \ \geq \ \esssup_{\nu\in\calv} \esssup_{\tau\in\Tc} \E^\nu\bigg[ \int_t^\tau f(r,X_r,I_r) \,dr + Y_\tau \,\bigg|\, \calf_t^{W,\mu}\bigg], \\
Y_t \ \geq \ \esssup_{\nu\in\calv} \essinf_{\tau\in\Tc} \E^\nu\bigg[ \int_t^\tau f(r,X_r,I_r) \,dr + Y_\tau \,\bigg|\, \calf_t^{W,\mu}\bigg].
\end{align*}
It remains to prove the reverse inequalities. As in the proof of Lemma 
\ref{bsdepenlizaedandcontrol}, for every $n$ and $\varepsilon \in(0,1)$, we define $\nu_r^{\varepsilon ,n}(a)=n\,1_{\{U_r^n(a)\geq0\}} + \varepsilon \,1_{\{-1<U_r^n(a)<0\}} - \varepsilon \,(U_r^n(a))^{-1}\,1_{\{U_r^n(a)\leq1\}}$. Then, $\nu^{\varepsilon ,n}\in\Vc_n$ and
\[
n\big(U_r^n(a)\big)^+ - \nu_r^{\varepsilon ,n}(a) U_r^n(a) \ \leq \ \varepsilon , \qquad \text{for all }r\in[0,T].
\]
Therefore, from equality \eqref{DPPRandomizedProblem4}, we find
\begin{align*}
Y_t^n \ &= \ \E^{\nu^{\varepsilon ,n}}\bigg[Y_\tau^n + \int_t^\tau f(r,X_r,I_r) \,dr + \int_t^\tau\int_A \big[n\big(U_r^n(a)\big)^+ - \nu_r(a) U_r^n(a)\big]\lambda(da)dr\,\bigg|\, \calf_t^{W,\mu}\bigg] \\
&\leq \ \E^{\nu^{\varepsilon ,n}}\bigg[Y_\tau^n + \int_t^\tau f(r,X_r,I_r) \,dr\,\bigg|\, \calf_t^{W,\mu}\bigg] + \varepsilon \,\lambda(A)\,T, \qquad \text{for all }\tau\in\Tc.
\end{align*}
Then, we obtain the two following inequalities:
\begin{align*}
Y_t^n \ &\leq \ \esssup_{\tau\in\Tc}\E^{\nu^{\varepsilon ,n}}\bigg[Y_\tau^n + \int_t^\tau f(r,X_r,I_r) \,dr\,\bigg|\, \calf_t^{W,\mu}\bigg] + \varepsilon \,\lambda(A)\,T, \\
Y_t^n \ &\leq \ \essinf_{\tau\in\Tc}\E^{\nu^{\varepsilon ,n}}\bigg[Y_\tau^n + \int_t^\tau f(r,X_r,I_r) \,dr\,\bigg|\, \calf_t^{W,\mu}\bigg] + \varepsilon \,\lambda(A)\,T.
\end{align*}
As a consequence, we get (we continue the proof with ``$\text{ess\,inf}$'' over $\tau\in\Tc$, since the proof with ``$\text{ess\,sup}$'' can be done proceeding along the same lines)
\[
Y_t^n \ \leq \ \esssup_{\nu\in\calv_n}\essinf_{\tau\in\Tc}\E^\nu\bigg[Y_\tau^n + \int_t^\tau f(r,X_r,I_r) \,dr\,\bigg|\, \calf_t^{W,\mu}\bigg] + \varepsilon \,\lambda(A)\,T.
\]
Using the arbitrariness of $\varepsilon $, and recalling that $\Vc_n\subset\Vc$ and $Y^n\leq Y$, we obtain
\[
Y_t^n \ \leq \ \esssup_{\nu\in\calv}\essinf_{\tau\in\Tc}\E^\nu\bigg[Y_\tau + \int_t^\tau f(r,X_r,I_r) \,dr\,\bigg|\, \calf_t^{W,\mu}\bigg].
\]
The claim follows letting $n\rightarrow\infty$.
\qed

\subsection{Randomized dynamic programming}\label{subs:randomdynprogr}

The purpose of this paragraph is to show that the result of 
Theorem \ref{T:DPPRandomizedProblem}
can be recast in a form very close to the classical dynamic programming principle. This result, that we call  randomized dynamic programming principle and is stated below in Proposition \ref{randDPP},   can be used to prove that the value function of the original classical control problem is a viscosity solution to the Hamilton-Jacobi-Bellman equation.
We will not introduce viscosity solutions formally, but for the interested reader we will indicate some references related to the randomization method.

We first formulate an optimal control problem and its randomized version for a system starting at time $t$ and evolving on the interval $[t,T]$. This is done mainly to introduce notation: since the constructions are entirely similar to 
Subsections \ref{Primal} and 
\ref{randomizedformulation} we omit many obvious details.

We assume that
$A,b,\sigma,f,g$
 are given and
satisfy the assumptions {\bf (A1)}. Let $(\Omega,\calf,\P)$ be a complete probability space equipped with
 a right-continuous and $\P$-complete filtration $\F=(\calf_t)_{t\ge 0}$, and
$W$   an $\R^{d}$-valued
standard Wiener process with respect to $\F$ and $\P$. For every $t\ge0$,
let  $\F^{W,t}=(\calf^{W,t}_s)_{s\ge t}$ be the right-continuous and $\P$-complete filtration generated by $(W_s-W_t)_{s\ge t}$.  For $t\in[0,T]$, denote by 
  $\Ac^{W,t}$ the set of   $\F^{W,t}$-progressive processes $A:\Omega\times [t,T]\to A$; it is the set of  admissible controls for the controlled problem started at time $t$.
Given 
 $t\in [0,T]$, $\alpha\in \Ac^{W,t}$ and a starting point $x\in\R^n$, the controlled equation has the form
\begin{equation*}
    \left\{
\begin{array}{lll}
  dX_s^{t,x,\alpha} & = & b(s, X_s^{t,x,\alpha}, \alpha_s)\,ds +
\sigma(s, X_s^{t,x,\alpha}, \alpha_s)\,dW_s, \qquad   s\in [t,T],
\\
X_t^{t,x,\alpha}&=&x,
\end{array}
\right.
\end{equation*} 
and the value function is defined by 
\begin{equation}\label{primalvaluefromt}
v(t,x) \ = \ \sup_{\alpha\in\Ac^{W,t}}\E\left[\int_t^Tf(s,X_s^{t,x,\alpha},\alpha_s)\,ds+g(X^{t,x,\alpha}_T)\right].
\end{equation}
  
Now take a finite Borel measure 
$ \lambda(da)$ on $A$, with full topological support. Assume that on 
$(  \Omega,   \calf,  \P)$ one can define a Poisson random measure
 $ \mu=\sum_{n\ge 1}\delta_{(  S_n,  \eta_n)}$
 with intensity $\lambda(da)$, independent of $W$; in fact, this involves no loss of generality, since this can be achieved up to a product extension, as described in  Remark  \ref{remchoice}.
For every $t\in [0,T]$ and every starting point $a\in A$
 let us define the $A$-valued, piecewise constant process $I^{t,a}$ as follows:
\begin{equation}
\label{Ifromt}
  I_s^{t,a}  =   \sum_{n\ge 0,\,S_n\le t}  a\,1_{[S_n,  S_{n+1})}(s)+  \sum_{n\ge 1,\,S_n>t}  \eta_n\,1_{[  S_n,  S_{n+1})}(s), \qquad s\in [t,T],
\end{equation}
where $S_0:=0$. This is defined in such a way that 
 $I^{t,a}_s=a$ for $t\le s<S_{n^*} $,
where $S_{n^*}$ is the first time  $S_n$ strictly greater that $t$; then 
$I^{t,a}_s=\eta_{n^*}$ for $s\in [S_{n^*}, S_{n^*+1})$ and so on.
Let $  X^{t,x,a}$ be the solution to the equation
\begin{equation}\label{randdynfromt}
    \left\{
\begin{array}{lll}
  dX_s^{t,x,a} & = & b(s, X_s^{t,x,a}, I_s^{t,a})\,ds +
\sigma(s, X_s^{t,x,a}, I_s^{t,a})\,dW_s \qquad   s\in [t,T],
\\
X_t^{t,x,a}&=&x,
\end{array}
\right.
\end{equation} 
 We define the filtration
$\F^{  W,  \mu,t}=(\calf^{  W,  \mu,t}_s)_{s\ge t}$
setting
\begin{equation*}
    \calf^{  W,  \mu,t}_s=\sigma \Big(  W_r-W_t,\,
  \mu((t,r]\times C)\,:\, r\in  [t,s],\, C\in\calb(A)\Big)
\vee \caln,
\end{equation*}
where $\caln$ denotes the family of $  \P$-null sets of $  \calf$.
We denote by $\calp(\F^{  W,  \mu,t})$ the corresponding predictable
$\sigma$-algebra on $\Omega\times [t,\infty)$.

We can now define the randomized optimal control problem via a change of probability measure of Girsanov type in the usual way.
We define the set 
of admissible controls from time $t$:
\begin{align*}
       \calv^t=\{
      \nu:\Omega\times [t,\infty)\times A\to \R,\;   \Pc(\F^{  W,  \mu,t})\otimes \calb(A){\rm  -measurable},\;
0< \nu_s(\omega,a)\le \sup  \nu<\infty\}.
\end{align*}
and the   exponential  $\F^{  W,  \mu,t}$-martingale
\begin{equation*}
    \kappa_s^{ \nu,t} \ 
= \ \exp\left(\int_t^s\int_A (1 -   \nu_r(a))\lambda(da)\,dr
\right)\prod_{t<  S_n\le s}\nu_{  S_n}(  \eta_n),\qquad s\ge t.
\end{equation*} 
We define a new probability on $( \Omega, \calf)$ setting
$ \P^{ \nu,t}(d \omega)=\kappa_T^{ \nu,t}( \omega)\, \P(d \omega)$ 
and the value function of the randomized control problem
\begin{eqnarray}  \label{defJrandomizedfromt}
v^\Rc (t,x,a)&=& \sup_{\nu\in\calv^t}  \E^{ \nu,t}
\left[\int_t^Tf(t,  X_s^{t,x,a},  I_s^{t,a})\,ds+g(  X_T^{t,x,a})\right].
\end{eqnarray}  

By Theorem \ref{MainThm},
 applied to the interval $[t,T]$, we have
\[
v(t,x)=v^\Rc (t,x,a),
\]
so that the latter function does not depend on $a$. 

Next we introduce a constrained BSDE on the time interval 
$[t,T]$:
\begin{equation}\label{BSDEconstrainedfromt}
\begin{cases} \dis Y_s \ = \ g(X_T^{t,x,a})  + \int_s^T  f(r,X_r^{t,x,a} ,I_r^{t,a}) dr + K_T - K_s - \int_s^TZ_r\,dW_r - \int_s^T\!\int_A U_r(a)\,\mu(dr\,da), \\
\dis U_s(a) \ \le \ 0.
\end{cases}
\end{equation}
By Theorem
\ref{Thm:RandomizedFormula}
there exists a unique minimal solution $(Y,Z,U,K)$ $\in$ $\Sc^2(\F^{W,\mu,t})\times L_W^2(\F^{W,\mu,t})\times L_{\lambda}^2(\F^{W,\mu,t})\times\Kc^2(\F^{W,\mu,t})$ to  \eqref{BSDEconstrainedfromt}. 
 The definition of these spaces  and the definition of  minimal solution is 
 the same as  
Definition \ref{BSDEdef}, replacing processes on $[0,T]$ by processes on $[t,T]$ and the filtration 
$ \F^{W,\mu}$ by $ \F^{W,\mu,t}$.
The minimal solution will be denoted
$(Y_s^{t,x,a},Z_s^{t,x,a},U_s^{t,x,a}(a),K_s^{t,x,a})$ $(s\in [t,T]$, $a\in A$$)$ to stress dependence on the parameters ${t,x,a}$.
  By 
Theorem
\ref{Thm:RandomizedFormulabis} we also conclude that for every $s\in [t,T]$,
\begin{equation*}
Y_s^{t,x,a} = \esssup_{\nu\in\calv^t} \E^{\nu,t}\bigg[ \int_s^T   f (r,X_r^{t,x,a},I_r^{t,a} )\,dr + g(X_T^{t,x,a}) \,\bigg|\, \calf_s^{W,\mu,t}\bigg]
\end{equation*}
and in particular, setting $s=t$,
\begin{align}\label{RandomizedFormuladetfromt}
    Y_t^{t,x,a}=v^\Rc(t,x,a)=v(t,x).
\end{align}
This formula shows that the value function of the original and the randomized control problem can be represented by the constrained BSDE \eqref{BSDEconstrainedfromt}.
Finally, we can apply
 Theorem \ref{T:DPPRandomizedProblem} (and more precisely formula 
\eqref{DPPRandomizedProblembis}, replacing the interval $[0,T]$ by $[t,T]$) and we obtain
\begin{eqnarray} \label{DPPRandomizedProblembisfromt}
Y_t^{t,x,a} &=& \esssup_{\nu\in\calv^t} \esssup_{\tau\in\Tc^t} \E^{\nu,t}\bigg[ \int_t^\tau f(r,X_r^{t,x,a},I_r^{t,a}) \,dr + Y_\tau^{t,x,a}\bigg] \notag \\
&=& \esssup_{\nu\in\calv^t} \essinf_{\tau\in\Tc^t} \E^{\nu,t}\bigg[ \int_t^\tau f(r,X_r^{t,x,a},I_r^{t,a}) \,dr + Y_\tau^{t,x,a} \bigg],
\end{eqnarray} 
where $\mathcal T^t$ denotes the class of $[t,T]$-valued $\F^{W,\mu,t}$-stopping times.

As a final step toward a dynamic programming principle we note that from its definition 
\eqref{Ifromt} the process satisfies the flow property: for $0\le t\le s\le  T$ we have $\P$-a.s.
\[
I_r^{t,a}= I_r^{s,I_s^{t,a}},\qquad   r\in [s,T].
\]
By standard arguments, based on  uniqueness for equation 
\eqref{randdynfromt}, it follows that for $0\le t\le s\le  T$ we have $\P$-a.s.
\[
X_r^{t,x,a}= X_r^{s, X_s^{t,x,a},I_s^{t,a}}, \qquad r\in [s,T].
\]
Since the pair $(I,X)$ satisfies this flow property,  the first component of minimal solution to \eqref{BSDEconstrainedfromt} satisfies a similar property: 
for $0\le t\le s\le  T$ we have $\P$-a.s.
\[
Y_r^{t,x,a}= Y_r^{s, X_s^{t,x,a},I_s^{t,a}}, \qquad r\in [s,T].
\]
so that setting $r:=s$ and recalling 
\eqref{RandomizedFormuladetfromt}
we obtain: for $0\le t\le s\le  T$ we have $\P$-a.s.
\[
Y_s^{t,x,a}= 
v^\Rc(s,X_s^{t,x,a},I_s^{t,a})= 
v(s,X_s^{t,x,a}).
\]
A rigorous proof proceeds by first establishing the corresponding property for the penalized equation (analogous to \eqref{bsdepenalized}) and then passing to the limit; we also refer the reader to 
\cite{BaCoCo19} Theorem 6.5,
where more details are included.
In \cite{BaCoCo19} Theorem 6.6 point 2), as a consequence of continuity properties of the function $v$, it is even shown that 
for $0\le t \le  T$ we have $\P$-a.s. 
\[
Y_s^{t,x,a}= 
v(s,X_s^{t,x,a}), \qquad s\in [t,T],
\]
and it follows that for any stopping time $\tau$ with values in $[t,T]$ we have $
Y_\tau^{t,x,a}= 
v(\tau,X_\tau^{t,x,a})$ $\P$-a.s. From 
\eqref{DPPRandomizedProblembisfromt}
we finally obtain the following conclusion, which holds under  
 the  assumptions and notations introduced above in this paragraph.
 
\begin{Proposition} \label{randDPP}
(randomized dynamic programming principle).
   We have
\begin{eqnarray} \label{DPPRandomizedProblemterfromt}
v(t,x) &=& \esssup_{\nu\in\calv^t} \esssup_{\tau\in\Tc^t} \E^{\nu,t}\bigg[ \int_t^\tau f(r,X_r^{t,x,a},I_r^{t,a}) \,dr + v(\tau,X_\tau^{t,x,a})\bigg] \notag \\
&=& \esssup_{\nu\in\calv^t} \essinf_{\tau\in\Tc^t} \E^{\nu,t}\bigg[ \int_t^\tau f(r,X_r^{t,x,a},I_r^{t,a}) \,dr + v(\tau,X_\tau^{t,x,a}) \bigg],
\end{eqnarray} 
where $\mathcal T^t$ denotes the class of $[t,T]$-valued $\F^{W,\mu,t}$-stopping times.
\end{Proposition}

As mentioned at the beginning of this paragraph, this result can be used to prove that  the value function $v(t,x)$ of the original control problem is a viscosity solution to the fully non linear Hamilton-Jacobi-Bellman equation
\eqref{hjbintro}. The arguments are similar to the ones classically used to deduce the viscosity property from the standard dynamic programming principle. For a full exposition we refer to \cite{BaCoCo19} Proposition 6.12, where the case of a Hilbert state space is treated with essentially the same arguments. A more complicated situation is treated in \cite{BCFP15}
in the case of viscosity solutions on the Wasserstein space, for a partially observed control problem (see also the discussion in Subsection \ref{subs:partialobs}
below). 
We will not report these proofs since introducing viscosity solutions formally  would take us too far.

\section{Other applications of the randomization method and possible future developments}
\label{S:appli}

The purpose of this section is to give a description of the large number of contexts where the randomization method has been  successfully applied. We will review applications to many stochastic optimization problems and we will try to give a complete picture of the existing results. Due to the great variety of addressed topics we keep an informal style and we only give references that are related to the randomization method, while we assume that the reader is familiar with the particular topic under examination. Moreover, for brevity reasons, some topics are just shortly mentioned but the corresponding results are not presented, even informally.

The final section concerns some indications for possible future research directions.

\subsection{Non-Markovian systems}

So far the controlled system has been defined as a stochastic equation of the form 
\begin{equation*}
    dX_t^\alpha = b(t,X_t^\alpha, \alpha_t)\,dt +
\sigma(t,X_t^\alpha, \alpha_t)\,dW_t, \qquad X_0^\alpha=x_0,
\end{equation*}
on the interval $[0,T]$, driven by a Brownian motion $W$. We note that the coefficients at time $t$ only depend on the current state $X_t$ and the current control action $\alpha_t$.
Equations of this form are often called of Markovian type: when the control is constant, or absent, the solutions to the equation may be used to define a Markov process in the usual way. However,  more general coefficients can be allowed in the system dynamics:
\begin{equation*}
    dX_t^\alpha = b_t(X^\alpha, \alpha)\,dt +
\sigma_t( X^\alpha, \alpha )\,dW_t, \qquad X_0^\alpha=x_0,
\end{equation*}
The notation indicates that
  the coefficients
$b,\sigma$ depend on the whole trajectory of $X^\alpha$ and $\alpha$. The dependence
will be non-anticipative, in the sense that  at any  time $t$ the coefficients $b_t,\sigma_t$ take values that depend on the trajectories
$X_s^\alpha$ and $\alpha_s$ for $s\in [0,t]$.  The non-anticipative requirement can be expressed in different ways, often as a measurability condition with respect to  canonical filtrations on the space of  paths.
Similar considerations also apply  to the payoff functionals $f,g$. 
 Thus, this formulation includes path-dependent (or hereditary) systems, and it allows for the presence of memory
effects both on the state and the control. 

In this non-Markovian situation the value function is not expected to satisfy a PDE of the usual form \eqref{hjbintro}, but methods based on BSDEs remain valid. In particular, it is worth noting that the randomization method carries over to this framework with only minor modifications. For this reason, non-Markovian cases are addressed in most of the papers we cite. The same applies to the various optimization problems that we are going to outline in the paragraphs below, although we will not indicate this at each point and we will only present the basic Markovian situation.

\subsection{Infinite horizon and other variants of the basic optimal control problem}

Let us consider a controlled equation on $[0,\infty)$ with coefficients that do not depend on time:
\begin{equation*}
    dX_t^\alpha = b(X_t^\alpha, \alpha_t)\,dt +
\sigma(X_t^\alpha, \alpha_t)\,dW_t, \;\;t\ge 0,\qquad X_0^\alpha=x,
\end{equation*}
with  value function defined as
\[ 
v^\beta(x)=    \sup_{  \alpha}
  \E 
\left[\int_0^\infty e^{-\beta t}f(  X_t^\alpha,  \alpha_t)\,dt\right],
\]
where the reward functional also depends on a  discount factor $\beta>0$. 
The randomization method has been applied to this infinite horizon optimal control problem in 
\cite{CoCoFu2019}, where the authors introduce an associated randomized problem, they prove equality of the corresponding values and represent them by a suitable constrained BSDE on the time interval $[0,\infty)$. In addition, in the Markovian case, they use the BSDE to construct a viscosity solution to the associated Hamilton-Jacobi-Bellman equation, which is fully non linear and of elliptic type:
\[
 \beta\, v^\beta(x)   =   \sup_{a\in A }\left[
\call^av^\beta( x) + f( x,a)
\right], \qquad\,x\in\R^n,
\]
where, for a regular function $\phi:\R^n\to\R$, we set
\[
\call^a\phi(x)= \frac{1}{2}\,{\rm Trace}\left[ D^2\phi(x)\,\sigma\sigma^T(x,a)\right] + 
D\phi(x)\cdot b(x,a).
\] 

Usually these results require the discount factor $\beta$ to be sufficienty large. Under additional assumptions, however, the function $v^\beta$ is defined for all $\beta>0$ and then the so-called ergodic problem can be studied. 
It may be 
defined in various ways and we limit ourselves to some generic indications. The ergodic problem may be viewed as the study of the limit $\lim_{\beta\to 0}\beta\,v^\beta(x)$. Under some assumptions, the limit exist - up to a subsequence - and it equals a constant $\lambda$; moreover, $v^\beta(x)-v^\beta(0) $ also admits a limit $w(x)$ that can be shown to be a solution to the ergodic equation
\[
\lambda  =   \sup_{a\in A }\left[
\call^aw( x) + f( x,a)
\right], \qquad\,x\in\R^n,
\]
Finally, the constant $\lambda$ and the function $w$ play a role in the study of the large time limit of the parabolic Hamilton-Jacobi-Bellman equation. 
In               
\cite{CossoFuhrmanPham16} it is shown that these results, as well as some generalizations, can be recovered by the randomization method, and more precisely using  the representation of the function $v^\beta$ by means of a constrained BSDEs.
Another ergodic problem is studied in the context of   optimal switching  in       \cite{BaCoPh2017}: see also paragraph 
 \ref{subs:optimalswitching}           below.
Finally, in 
\cite{CoGuaTe2019} an ergodic problem
is studied for a system evolving in a Hilbert space, again using the randomization method.

We also mention
\cite{ChoukrounCosso16}, where a constrained BSDE representation is given for the value function of a stochastic control problem, on finite interval $[0,T]$, where the noise in the equation also depends on a jump term, whose intensity is also controlled. The associated Bellman equation is therefore of integro-differential PDE.

A jump term in the controlled equation is also contained in the system studied in \cite{BaCoCo19}, in the case of an infinite-dimensional state space. The authors prove various results and construct a solution to the Bellman equation using the randomization method.

\subsection{Optimal control with partial observation}
\label{subs:partialobs}

Let us consider  a system driven by a Brownian motion $B$ and of the form  
\begin{equation*}
    dX_t^\alpha = b(X_t^\alpha, \alpha_t)\,dt +
\sigma(X_t^\alpha, \alpha_t)\,dB_t, \qquad X_0^\alpha=x_0,
\end{equation*}
on the interval $[0,T]$, where the initial condition  is now, more generally, a random variable $x_0$ independent of $B$. Here  the coefficients may depend explicitly on time but we take them to be time-invariant in order to to simplify notation. 
 We suppose that the process $B$ driving the controlled equation is of the form $B=(V,W)$ for two independent Brownian motions $V,W$ (in general, both multidimensional). 
In the partial observation problem formulated in 
\cite{BCFP16c}, \cite{BCFP15}  
the authors assume that the controller can not directly observe the full noise sources $x_0,V,W$ nor the controlled process $X$, but   he/she must instead base the choice of control actions  on observation  of the component $W$ alone. This is modeled by requiring that   the control process $\alpha$ should be progressive with respect to the filtration $\F^W$ generated by $W$. So the reward is defined as
\[ 
{\text{\Large$\upsilon$}}_0  \;=\;   \sup_{  \alpha}
  \E 
\left[\int_0^Tf(  X_t^\alpha,  \alpha_t)\,dt+g(  X_T^\alpha)\right],
\]
where $\alpha$ ranges over $\F^W$-progressive processes.
\begin{Remark}
    \emph{ 
In classical, but more general, partial observation models, the controlled  can only observe 
another stochastic process $Y^\alpha$ of the form
\begin{equation*}
    dY_t^\alpha = h(X_t^\alpha, \alpha_t)\,dt +
\,dW_t,\qquad Y_0^\alpha= 0,
\end{equation*}
(possibly with more general coefficients) 
and the control is required to be progressive with respect to the filtration generated by $Y^\alpha$. It can be shown, see \cite{BCFP16c}, that this problem can be recast in the form presented before.
}\qed 
\end{Remark}
The associated randomized control process is introduced as the solution to 
\begin{equation*}
  d  X_t=  b(   X_t,  I_t)\,dt + \sigma(  X_t,  I_t)\,d  B_t, \qquad 
 X_0 =x_0,
\end{equation*} 
where the $A$-valued process $I$ is, as before, the piecewise-constant process 
\[
  I_t \ = \ \sum_{n\ge 0}  \eta_n\,1_{[  S_n,  S_{n+1})}(t), \qquad t\ge 0,
\]
corresponding to an independent Poisson random measure 
$ \mu=\sum_{n\ge 1}\delta_{(  S_n, \eta_n)}$  with intensity measure $\lambda(da)$ being finite and with full support on $A$. 

Next we define an optimization problem introducing a space of control strategies   $\calv$ 
whose elements are bounded
 random fields 
\[\nu=\nu_t(\omega,a):\Omega\times [0,\infty)\times A\to (0,\infty).
\]
Then, for every $\nu\in\calv$, one defines as before a martingale $\kappa^\nu$ and Girsanov change of measure leading to a  probability $\P^\nu$. 
The reward functional is  defined as
\[ 
{\text{\Large$\upsilon$}}_0^\Rc \;=\;   \sup_{  \nu\in  \calv}
  \E^{ \nu}
\left[\int_0^Tf(  X_t,  I_t)\,dt+g(  X_T)\right],
\]
where $\E^{ \nu}$ denotes the expectation under $\P^{ \nu}$.
We note that the process $X$ is adapted to the filtration $\F^{\mu,x_0,B}=\F^{\mu,x_0,V,W}$
generated by the random measure $\mu$, by $x_0$ and by $B=(V,W)$. Instead, to take into account the measurability requirement imposed on the control by the feature of partial 
observation, the random fields $\nu\in\calv$ are required to the predictable with respect to $\F^{\mu,W}$, the filtration generated by $\mu$ and $W$ only. 
 With this modification one proves that 
$
{\text{\Large$\upsilon$}}_0=
{\text{\Large$\upsilon$}}_0^\Rc$
and that this common value can be represented using the unique minimal solution  of a constrained BSDE (solved under the filtration 
 $\F^{\mu,W}$). 

In my opinion this result, contained in  \cite{BCFP16c}, is a clear indication to  the effectiveness of the randomization method: indeed, at least to  my knowledge, this is the first time that a representation formula was given for the value of a general partially observed control problem in terms of a solution to a BSDE.

The result in  \cite{BCFP16c} deals with more general path-depending coefficients;   the Markovian case,
presented above, is further investigated in \cite{BCFP15}. In the latter paper the authors consider the case of a control problem starting at a generic time $t$ from a generic random variable $x_0$. They show that the corresponding value function does not depend on the specific variable $x_0$, but it is only a function of the law $\rho$ of $x_0$ (as well as the coefficients $b,\sigma,f,g$) and so it can be seen as a function $v(t,\rho)$ defined for $t\in [0,T]$ and $\rho\in\calp_2(\R^n)$, the Wasserstein space of probabilities on $\R^n$ with finite second moment. 
Next, a fully non linear Hamilton-Jacobi-Bellman equation is derived  and it is shown that  the value function is the unique viscosity solution. 
Thus, even for this class of PDEs on the Wasserstein space, a constrained BSDE has been introduced to construct solutions. It should be noted that, when $\rho$ admits a density $f$ with respect to the Lebesgue measure, belonging to some appropriate Hilbert space $H$ of functions defined on $\R^n$, the value function can be considered as a function $v(t,f)$ 
defined for $t\in [0,T]$ and $f\in H$ and there is a vast literature on Bellman's equation as a PDE on $H$: we refer the reader to \cite{gozswi15}
for further information. Using the randomization method it is possible to deal with the general case $v(t,\rho)$.
Although in \cite{BCFP15} there are some details missing in the proof of uniqueness of viscosity solution, the topic of uniqueness has been further investigated by other authors: we refer the reader to 
\cite{BaEkZh23} or \cite{BaChEkQiTaZh25}
and the references therein.

\subsection{Optimal switching}
\label{subs:optimalswitching}

Optimal switching problems are historically important for the randomization method, since the first formulation of the randomized system, due to B. Bouchard   \cite{bou09}, was given in this context. 

In optimal switching problems one consider a controlled equation 
\begin{equation*}
    dX_t^\alpha = b(X_t^\alpha, \alpha_t)\,dt +
\sigma(X_t^\alpha, \alpha_t)\,dW_t
\end{equation*}
on the interval $[0,T]$ with initial condition $X_0^\alpha=x_0\in\R^n$.
In the standard case the space of control actions $A$ is a finite set, that we identify with $\{1,2,\ldots,m\}$; its elements are called modes or regimes.
The control process is piecewise constant and has the form
\begin{equation}
    \label{formaswitch}
\alpha_t=\xi_0\,1_{[0,\tau_1)}(t)+\sum_{n\ge 1}\xi_n\,1_{[\tau_n,\tau_{n+1})}(t),
\qquad t\in [0,T],
\end{equation}
and is therefore determined by the initial regime $\xi_0$ and by a strategy which is a double sequence
$$
\alpha=(\tau_n,\xi_n)_{n\ge 1}.
$$
Let   $\F^W=(\calf^W_t)_{t\ge 0}$ denote the Brownian filtration; then
$(\tau_n)$ is required to be an  increasing sequence of  $\F^W$-stopping times diverging to $\infty$, and each $\xi_n$ is an $\calf^W_{\tau_n}$-measurable, $A$-valued random variable. Thus, the dynamics of the controlled system can be chosen within a finite set of drifts and volatilities $b(\cdot,i)$, $\sigma(\cdot,i)$, each corresponding to a possible regime $i\in A$,  and  the controller chooses the regimes and the switching times $\tau_n$ from regime $\xi_{n-1}$ to $\xi_n$. 
The functional to be maximized has the form 
\begin{align*}  
J(\alpha)   =   \E\Big[\int_0^Tf(X^\alpha_t,\alpha_t)\,dt+g(X^\alpha_T,\alpha_T)\Big]
- 
   \E\Big[\sum_{n\ge 1}1_{\tau_n\le T}\,c_{\tau_n}(X^\alpha_{\tau_n},\xi_{n-1},\xi_n)\Big].
\end{align*}
The first term is the standard running and terminal cost, the latter depending on the terminal regime as well.
The second term represents the switching cost: the nonnegative function $c_t(x,i,j)$ is the cost incurred to switch from regime $i$ to regime $j$ at time $t$ if the current state is $x$. Usually the cost functions are bounded away from $0$, so that infinitely many switchings would give $J=-\infty$: this is the reason to choose control processes of the form  \eqref{formaswitch}.
We define the value of the optimal switching problem taking the supremum  over all admissible strategies:
\begin{equation}\label{primalvalueswitching}
{\text{\Large$\upsilon$}}_0 = \sup_{\alpha} J(\alpha).
\end{equation}

It is possible to represent the value using  BSDEs, more precisely solving a system of reflected BSDEs with interconnected obstacles, indexed by $i\in A$: one looks
  for unknown adapted processes $(\bar Y^{i}_t, \bar Z^{i}_t, 
\bar K^{i}_t)_{t\in[0,T]}$, satisfying suitable conditions, such that
\begin{equation}\label{reflected_intro}
\left\{
\begin{array}{l}\dis
\bar Y^{i}_t+ \int_t^T \bar Z^{i}_s\,dW_s =
g(\bar X^{i}_T) + \int_t^Tf(\bar X^{i}_s,i)\,ds +\bar K^{i}_T-\bar K^{i}_t,
\\
\dis \bar Y^{i}_t\ge \max_{j\neq i} [\bar Y^{j}_t- c_t(\bar X^{i}_t,i,j)],
\\\dis
\int_0^T\Big[\bar Y^{i}_t- \max_{j\neq i} [\bar Y^{j}_t- c_t(\bar X^{i}_t,i,j) ]\Big]\,d\bar K^{i}_t=0,
\end{array}
\right.
\end{equation}
where, in particular, $\bar K^{i}$ are non decreasing processes, $\bar K^{i}_0=0$,
and $\bar X^{i}$  are defined by the equations
$$
d\bar X^{i}_t =  b( \bar X_t^{i}, i)\,dt +
\sigma( \bar X_t^{i}, i)\,dW_t, \quad t\in [0,T],
\qquad \; \bar X_0^{i}=x_0.
$$
Under suitable conditions this system is well-posed and one has
a probabilistic representation for the value: if the initial regime is $\xi_0=i$ then ${\text{\Large$\upsilon$}}_0 =\bar  Y^{i}_0$.

Another possibility is to use a dynamic programming approach and to characterize the value using a system of  PDEs. Let us consider the switching problem starting at a generic time $t\in [0,T]$  with initial state $x\in \R^n$ and regime $i\in A$. Let $v(t,x,i)$ denote the corresponding value. 
Then we expect   these functions to be solution to
the system:
for $i=1,\ldots, m$
\begin{equation}\label{QVIsystem_intro}
\left\{
\begin{array}{l}\dis
\min\left\{
-\partial_tv(t,x,i)-
\call^iv(t,x,i) - f(x,i),
v(t,x,i)-\dis\max_{j\neq i} [v (t,x,j)- c_t(x,i,j)]
\right\}=0,
\\
v(T,x,i)=g(x,i),  \qquad x\in\R^n,\;t\in [0,T),
\end{array}
\right.
\end{equation}
where
$$
\call^iv (t,x,a)=\frac12{\,\rm Trace}\,[\sigma(x,i)\sigma(x,i)^TD^2_xv(t,x,i)]
+D_xv(t,x,i)b(x,i)
$$
is the Kolmogorov operator corresponding
to the controlled coefficients $b(x,i)$, $\sigma(x,i)$.

In particular, a system of BSDEs entirely similar to 
 \eqref{reflected_intro} provides a representation formula for the solution functions $v(t,x,i)$, since they both coincide with the value of the same control problem.

The randomization method provides a third possibility to represent the value. 
It works as follows.
Take a measure $\lambda$ on $A$, identified with a vector $(\lambda(i))_{i\in A}$, and assume $\lambda(i)>0$ for all $i$. Next, as in Subsections \ref{randomizedformulation} or \ref{subsub;randomized}, 
consider a Poisson random measure on $A$ with intensity $\lambda$, independent of $W$, of the form
 $ \mu=\sum_{n\ge 1}\delta_{( S_n, \eta_n)}$.
We also define the $A$-valued, piecewise constant process associated to $ \mu$ as follows:
\[
 I_t \ = \ \sum_{n\ge 0} \eta_n\,1_{[ S_n, S_{n+1})}(t), \qquad t\ge 0,
\]
where   $ S_0=0$ and $  I_0=\eta_0$ is equal to the initial regime. 
Let $ X$ be the solution to the equation
\[ 
d X_t = b( X_t, I_t)\,dt + \sigma( X_t, I_t)\,d W_t,\;\;t\in [0,T], 
\qquad  X_0 = x_0.
\]
 We introduce the filtration
$\F^{ W, \mu}=(\calf^{ W, \mu}_t)_{t\ge 0}$ generated by the Wiener process and the Poisson random measure and we denote by $\calp(\F^{ W, \mu})$ the corresponding predictable
$\sigma$-algebra. We recall that the $\F^{ W, \mu}$-compensator of $ \mu$ is 
$\lambda(i)\,dt$.
We can now define the randomized optimal control problem via a change of probability measure of Girsanov type.  
We define the set 
of admissible controls
\begin{align*}
      \calv=\{
      \nu:\Omega\times [0,\infty)\times A\to \R,\;   \Pc(\F^{ W, \mu})\otimes \calb(A)-{\rm  measurable},\;
0< \nu_t(\omega,i)\le \sup  \nu<\infty\}.
\end{align*}
We define
  the   exponential process
\begin{align*}
\kappa_t^{ \nu} \ 
= \ \exp\left(\int_0^t\int_A (1 -  \nu_s(i))\lambda(di)\,ds
\right)\prod_{0< S_n\le t}\nu_{ S_n}( \eta_n),\qquad t\ge 0, 
\end{align*}
where the integral over $A$ is in fact a finite sum. $\kappa^{ \nu} $
is a martingale with respect to $ \P$ and $\F^{ W, \mu}$
and we   define a new probability   setting
$ \P^{ \nu}(d \omega)=\kappa_T^{ \nu}( \omega)\, \P(d \omega)$.  
Under $ \P^{ \nu}$
the $\F^{ W, \mu}$-compensator of $ \mu$ on the set
$[0,T]\times A$ is the random measure $ \nu_t(i)\lambda(i)dt$. 
We finally introduce the gain functional of the randomized control problem
\begin{align*}  
J^\Rc( \nu) =  \E^{ \nu}
\left[\int_0^Tf( X_t, I_t)\,dt+g( X_T,I_T)\right]-
  \E^{ \nu}
\Big[\sum_{n\ge 1}1_{S_n\le T}
\,c_{  S_n}(  X_{  S_n},  \eta_{n-1}, \eta_n)\Big]
,
\end{align*}
and the corresponding value   $
{\text{\Large$\upsilon$}}_0^\Rc \;=\;   \sup_{ \nu\in \calv} J^\Rc( \nu)$.

The first result is the equivalence of this problem with the original switching problem, in the sense that one can prove that the values are the same:  $
{\text{\Large$\upsilon$}}_0  =
{\text{\Large$\upsilon$}}_0^\Rc$.

The values can be represented by means of an associated BSDE with constrained jumps. For the switching problem it takes the form
\begin{equation}\label{BSDEconstrainedswitching}
\begin{cases}
\vspace{2mm} \dis Y_t \ = \ g( X,  I_T)
+ \int_t^T f(X_s,I_s)\, ds + K_T - K_t - \int_t^TZ_s\,dW_s
- \int_{(t,T]}\!\int_A U_s(i)\,\mu(ds\,di), \\
\dis U_t(i) \ \le \ c_t(X_t, I_{t-},i),
\end{cases}
\end{equation}
(here, again, 
  the integral over $A$ is in fact a finite sum). 
The BSDE is solved under  $\P$ and  in the filtration $\F^{W,\mu}$,  
on the time interval $[0,T]$. As before the
  unknown is a quadruple
$(Y_t,Z_t,U_t(i),K_t)$, $(t\in [0,T]$, $i\in A)$ where, in particular, $Y$ and $Z$ are progressive, the process $K$ and the random field $U$ are predictable, $K$ is non-decreasing and $K_0=0$. 
The main difference is the constraint on $U$, which is not simply the requirement that  $U_t(i)\le 0$ but it involves the switching costs.

A second result is well-posedness of the BSDE, in the sense that it admits a unique minimal solution (minimality is defined as before). 
Finally, we have the representation formula:
\[Y_0=  
{\text{\Large$\upsilon$}}_0  =
{\text{\Large$\upsilon$}}_0^\Rc.
\]
Minor variations allow to represent the value functions $v(t,x,i)$ introduced above by means of a constrained BSDE of similar form on the time interval $[t,T]$.

\vspace{2mm}

These results on randomization of the switching problem are contained in 
\cite{bou09}, \cite{EKa}, 
\cite{EKb}. There one may also find formulae that connect the solution to the constrained BSDE \eqref{BSDEconstrainedswitching} and the reflected system of BSDEs \eqref{reflected_intro}.
In \cite{EK2010} the authors also introduce constrained BSDEs of similar form to represent solutions to system of PDEs with interconnected obstacles that generalize \eqref{QVIsystem_intro}; they also address some issues concerning the numerical approximation of the solution. We finally mention \cite{BaCoPh2017} where the authors use randomization techniques to investigate the long-time behavior of some switching control problem in robust form. Finally, using the randomization method, one can consider the case when the number of modes is infinite, possibly uncountable, so that the control action space can be very general: see \cite{FuhMor2020}.

\subsection{Optimal stopping}

This optimization problem is usually  easier than the previous one and could even be seen as a special case. However, an explicit treatment has been given in 
\cite{FuhrmanPhamZeni16}. The reader will find formulae connecting the solution to a constrained BSDE, typical of the randomization method, with the standard reflected BSDE associated to optimal stopping.

\subsection{Impulse control}

In impulse control problems the controller chooses an intervention strategy 
$$
\alpha=(\tau_n,\xi_n)_{n\ge 1}
$$
where, as in the case of optimal switching,
 $(\tau_n)$ is  an  increasing sequence of  $\F^W$-stopping times diverging to $\infty$, and each $\xi_n$ is an $\calf^W_{\tau_n}$-measurable $A$-valued random variable,    $\F^W=(\calf^W_t)_{t\ge 0}$ being  the Brownian filtration. The system starts at time $0$ from a given state $x_0\in \R^n$ and evolves on the interval $[0,\tau_1)$ according to the dynamics
\begin{equation}\label{dinlibera}
    dX_t = b(X_t)\,dt +
\sigma(X_t)\,dW_t.
\end{equation}
At time $\tau_1$ the state of the system is switched to
\[
X_{\tau_1}=X_{\tau_1-}+\gamma(X_{\tau_1-},\xi_1)
\]
for a given function $\gamma(x,a)$ defined on  $\R^n$ $\times $ $A$. Then the system evolves freely on $[\tau_1,\tau_2)$ according to equation \eqref{dinlibera} and at time $\tau_2$ the state is given by
\[
X_{\tau_2}=X_{\tau_2-}+\gamma(X_{\tau_2-},\xi_2).
\]
This process continues up to a time horizon $T$ and the functional to be maximized, over the choice of all possible $\alpha$, has the form
\[
\E\left[ 
\int_0^Tf(X_s)\,ds + g(X_T) - \sum_{\tau_i\le T} c(X_{\tau_i-},\xi_i)
\right]
\]
for given real functions $f$ and $g$, representing running and terminal gain, and for given nonnegative function $c(x,a)$ representing the cost to choose an action $a\in A$ when the system is going to enter state $x\in\R^n$.

When the system starts at a generic time $t$ from a generic state $x$ let us  denote by $v(t,x)$ the corresponding value function. It is expected that $v$ is a solution to the so-called quasi-variational inequality
\begin{equation}\label{quasiQVIsystem_intro}
\left\{
\begin{array}{l}\dis
\min\left\{
-\partial_tv(t,x)-
\call v(t,x) - f(x),
v(t,x)-\calh v(t,x)
\right\}=0,
\\
v(T,x)=g(x),  \qquad x\in\R^n,\;t\in [0,T),
\end{array}
\right.
\end{equation}
where
$$
\call v (t,x)=\frac12{\,\rm Trace}\,[\sigma(x )\sigma(x )^TD^2_xv(t,x )]
+D_xv(t,x )\cdot b(x )
$$
is the Kolmogorov operator corresponding 
  to  $b$, $\sigma$ and $\calh$ is the nonlocal operator
\[
\calh v(t,x)= 
\sup_{a\in A} [v (t,x+\gamma(x,a))+ c (x,a)].
\]

In the paper
\cite{KMPZ10} the authors introduce a suitable BSDE with constrained jumps to construct a viscosity solution to the equation \eqref{quasiQVIsystem_intro}, thus establishing a connection with the optimal impulse problem as well. This is one of the first results on representation of a fully non linear PDE by means of constrained BSDEs.

\subsection{Applications to stochastic game theory}

The randomization method has been occasionally used in this context. We mention
\cite{CCH15} for a specific form of a two-player game, and 
the paper  \cite{BaCoPh2017} where a robust control problem (in optimal ergodic switching) has been modeled as a two-player game as well.

\subsection{Applications to PDEs}

As explained in Section \ref{s;informalintro}, one of the merits of the randomization method is to show that some classes of PDEs, even fully non linear, admit a representation in terms of an associated BSDEs. 
In the present paper we have mostly focused on the randomization method as a tool for representing the value of a control problem by means of the BSDE, but one can also establish a direct connection between the constrained BSDE and a PDE of interest. Usually the class of such PDEs is a generalization of some class of Hamilton-Jacobi-Bellman equation. Correspondigly, the BSDE also needs to have a more general form.

For instance, in \cite{KP12} 
the authors use a constrained BSDE to construct a viscosity solution to an integro-partial differential equation that generalizes the basic equation \eqref{hjbintro}, where the leading operator also contains nonlocal terms.  
In 
\cite{CPH15} the constrained BSDE is used to obtain viscosity solutions to a {H}amilton-{J}acobi-Bellman equation with improved growth conditions on the coefficients, while in 
\cite{CossoFuhrmanPham16} the long-time behaviour of such PDEs is studied.
In 
\cite{EK2010} a system of PDEs with interconnected obstacles is studied, that generalizes \eqref{QVIsystem_intro}.

\subsection{Numerical aspects}

Since constrained BSDEs can represent - or be used to construct - solutions to nonlinear PDEs, appropriate methods for simulation of those stochastic equations may lead to efficient numerical approximation techniques for the solutions to PDEs, even in cases where other more classical numerical methods suffer from limitations and difficulties. Numerical treatments of constrained BSDEs can be found in 
\cite{KLP14} and
\cite{KLP15}. Numerical approximation of the value function of optimal switching problems, based on randomization techniques,  is addressed in \cite{EK2010} and
\cite{DePhWa25}.

\subsection{Optimal control of McKean-Vlasov equations}

This name is given to a class of stochastic differential equations of the form
\begin{equation*}
    dX_t  = b(X_t,\P_{X_t}, \alpha_t)\,dt +
\sigma(X_t, \P_{X_t},\alpha_t)\,dW_t
\end{equation*}
namely where the coefficients also depend on the law $\P_{X_t}$ of the controlled process $X_t$. These models typically arise as limits of stochastic models for 
$N$ interacting agents when $N\to\infty$ and the interaction is of mean-field type, see for instance \cite{CaDe_book_I}-\cite{CaDe_book_II} 
as a general reference.

The randomization method is implemented in this framework in 
\cite{BaCoPh2018}. In 
\cite{DeKhPh24} the authors also consider the case of so-called  common noise.

\subsection{Infinite-dimensional state space}

In the optimization problems mentioned so far the  controlled state process  evolves in a finite-dimensional Euclidean space. There exists a huge literature on stochastic optimal control problems for systems in infinite-dimensional spaces, especially Hilbert spaces, see for instance \cite{gozswi15}.
These models also encompass control problems for stochastic partial differential equations. The corresponding Hamilton-Jacobi-Bellman equation is written for a value function depending on a space variable which typically ranges in a Hilbert space $H$, although more general cases may be considered. Besides the intrinsic difficulty of dealing with PDEs on $H$, additional complications are given by the occurrence of an unbounded linear operator (only defined on a subspace of $H$)  in the drift of the state equation: these operators are needed to allow immediate  applications to stochastic PDEs.  

The randomization method has been used to improve some results in this area as well. In 
\cite{BaCoCo19} the authors introduce a constrained BSDE to represent the value of the optimal control problem for an  $H$-valued stochastic equation and to construct a viscosity solution to the Hamilton-Jacobi-Bellman equation on $H$. In 
\cite{CoGuaTe2019} the constrained BSDE is introduced in connection with an ergodic problem where appropriate analytic results are missing for the corresponding ergodic PDE on $H$.

\subsection{Optimal control of non-diffusion processes}

There exist many stochastic models where controlled processes are not constructed as solutions to differential equations. An important case is the optimal control theory for continuous-time Markov processes, see for instance  
\cite{GuHe_book_2009}. 
It shares many technical aspects with the theory corresponding to classical controlled stochastic equations, in particular one can introduce the Hamilton-Jacobi-Bellman equation where, in many cases, the leading operator is of integral type and it is non-local.

In \cite{BandiniFuhrman15} it is shown that
the randomization method makes it possible to find a BSDE representation for the solution to such equations, on very general state spaces, and so for the value function of the corresponding optimal control problems. 

A more difficult related problem is the optimal control of piecewise deterministic Markov processes introduced by M.H.A. Davis, see 
\cite{DavisPDMP93}. In this case the nonlinear operator in the Hamilton-Jacobi-Bellman equation also contains first order derivatives with respect to space variables.  The results in 
 \cite{Bandini15} show that
the solution, which corresponds to the value function, can be represented by a constrained BSDE by another application of the randomization method.

\subsection{Indications for possible future research}

In this section we have tried to give a fairly complete picture of the literature on the randomization method. It is apparent that the applications to various areas of stochastic optimization are rather unequal. For instance, applications to control of  infinite-dimensional systems are rather few. In particular, there are no examples where the method has been applied to specific classes of stochastic PDEs, while optimal control theory has now a broad range of available results, see for instance \cite{gozswi15}.
Below we try to indicate some other research directions that seem particularly promising or important.

Concerning the randomization method itself, some effort should be devoted to simplifying the arguments needed to prove that the value of the original control problem equals the value of the randomized problem: this would greatly help in specific applications. Another important point would be to establish a more direct connection between the randomized BSDE and the starting control problem. More precisely, it would be of some importance to show how the constrained  BSDE may be used to construct optimal controls (or approximate optimal controls) for the original system. Finally, the numerical treatment of the constrained BSDE should be further refined and numerical   methods should also be developed that may help to find numerical approximations of the optimal control of the original problem.

On the side of possible applications, we note that stochastic game problems have been seldom studied using the randomization method.  Even rather classical results on zero-sum two-player games have not been studied using this method. This study may open the possibility of extending investigations to mean-field games: these stochastic models share some similarity with McKean-Vlasov systems shortly described above and have been the subject of great interest in the past few years, see \cite{CaDe_book_I}, \cite{CaDe_book_II}.

Another field in stochastic optimization that has not been studied so far using the randomization method is singular control, where the control process is chosen in a class having increasing paths (componentwise). Besides its intrinsic interest, for example in the modeling of irreversible investments, one may recall that optimal stopping problems are sometimes studied introducing generalized stopping times that are in fact defined using increasing control processes. The randomization method may help to improve existing results in these areas as well.

We finally mention another possible use of the randomization method that has not been explored so far. We note that the method does not require any nondegeneracy assumption on the noise acting on the system. In fact, the noise may be completely absent and the starting control problem can be of deterministic type. The auxiliary randomized control problem, however, is always stochastic and leads to the randomized BSDE. Thus, this BSDE can represent the value even of a deterministic optimal control problem. In this particular case one may try to obtain more precise results, possibly with simpler proofs. In particular,  numerical approximation of the solution to the BSDE may help in approximating the value function of a deterministic optimal control: this is in general a difficult task, since the regularity of the value function can be very low.

It is our hope that the open problems around the randomization method may attract the interest of other researchers so that it may become a more powerful tool in the analysis of stochastic optimization problems.

\appendix

\section{Appendix}

\setcounter{Theorem}{0}
\setcounter{equation}{0}

This section is devoted to the proof of Proposition \ref{extensionapproximation} below, which was used in the proof of
Theorem \ref{MainThm}. 
It consists in a point process construction that was stated without proof as Proposition 4.1 in \cite{BCFP16c} and contained  in the preprint version
\cite{BCFP16b}
of the same paper. Here we essentially follow the exposition in the appendix to
\cite{BCFP16b}, that was not published before.

We assume that $A$ is a Borel space,
and that $\lambda$ and $a_0$
 are given and satisfy the assumption {\bf (A2)}. Our starting point is also a probability  space $(\Omega,\calf,\P)$, with a filtration
 $\G=(\calg_t)_{t\ge 0}$.
In the rest of this section we do not need to have completed filtrations.

Recall that 
for any pair  $\alpha^1,\alpha^2:\Omega\times [0,T]\to A$ of   $\G$-progressive (more generally, measurable) 
  processes we have introduced the distance $\tilde\rho(\alpha^1,\alpha^2)$ setting
\begin{eqnarray*}
\tilde \rho(\alpha^1,\alpha^2) &=& \E \Big[\int_0^T\rho(\alpha^1_t,\alpha^2_t)\,dt \Big],
\end{eqnarray*}
where $\rho$ is an arbitrary metric in $A$ satisfying $\rho<1$.

Below we will use an auxiliary
probability space denoted
 $(\Omega',\calf',\P')$. This can be taken as an arbitrary  probability space
 where appropriate random objects  are defined.
For integers $m,n,k\ge1$,  we assume that real random variables  $ U_n^m $, $V^m_n$
and random measures
 $\pi^k$ are defined on  $(\Omega',\calf',\P')$ and satisfy the following conditions:
\begin{enumerate}
\item every $ U_n^m $ is uniformly  distributed on $(0,1)$;
\item  each   $ V_n^m$ has exponential distribution with parameter $\lambda_{nm}$
and $\sum_{n=1}^\infty\lambda_{nm}^{-1}=1/m$ for every $m\ge 1$;
\item  every $\pi^k$ is a Poisson random measure on $(0,\infty)\times A$, admitting compensator $k^{-1}\lambda(da)\,dt$ with respect to its natural filtration;
 \item the random elements  $U^m_n,V^h_j$, $\pi^k$    are all independent.
 \end{enumerate}

The role of these random elements   will become clear in the constructions that
follow. Notice that for the construction of the space
$(\Omega',\calf',\P')$ only the knowledge of the measure $\lambda$ is required.

\begin{Remark}
    {\rm 
At first sight the space $(\Omega',\calf',\P')$ may look rather complex. However, in order to construct the required random elements we only need the existence of a sequence of independent random variables having prescribed distributions. 
    By a classical construction,
 see for instance \cite{Zabczyk96} Theorem 2.3.1, it is enough to take $\Omega'=(0,1)$,
 $\calf'=\calb((0,1))$ and $\P'$ the Lebesgue measure.
    }
    \qed
\end{Remark}
Next we define
\begin{eqnarray*}
\hat\Omega \; = \; \Omega\times \Omega',
\qquad
\hat \calf \; = \; \calf\otimes \calf',
\qquad
\Q \; = \; \P\otimes \P'
\end{eqnarray*}
and note that
the filtration $\G$ can be canonically extended to a filtration $\hat\G=(\hat\calg_t)_{t\geq 0}$
in $(\hat\Omega,\hat\calf)$ setting $\hat\calg_t=\{A\times \Omega'\,:\, A\in\calg_t\}$.
Similarly, any process $\alpha$ in $(\Omega,\calf)$ admits an extension $\hat\alpha$
to $(\hat\Omega,\hat\calf)$ given by $\hat\alpha_t(\hat\omega)=\alpha_t(\omega)$,
where $\hat\omega=(\omega,\omega')$.
The metric $\tilde\rho$ can also be extended to
  any pair  $\beta^1,\beta^2:\hat\Omega\times [0,T]\to A$ of   $\hat\G$-progressive
  processes  setting
\begin{eqnarray*}
\tilde \rho(\beta^1,\beta^2) &=& \E^\Q \Big[\int_0^T\rho(\beta^1_t,\beta^2_t)\,dt \Big].
\end{eqnarray*}
We use the same symbol $\tilde\rho$ to denote  the extended metric as well.

\vspace{1mm}

Our aim in this section is to prove the following result. We note at the outset that the requirement that $\lambda(da)$ has full topological support is needed to define the kernels $q^m(b,da)$
in \eqref{kernelsuA}.

\begin{Proposition}
\label{extensionapproximation}
Let  $A$ be a Borel space, and let $\lambda$ and $a_0$
 satisfy   {\bf (A2)}.
Let  $(\Omega,\calf,\P)$ be  any probability
 space with a filtration
 $\G=(\calg_t)_{t\ge 0}$ and let
$(\hat\Omega,
\hat \calf,\Q)$ be the product space defined above. Then
 for any $\G$-progressive  $A$-valued
 process $\alpha$, and for any $\delta >0$,
 there exists a marked point process $(\hat S_n,\hat \eta_n)_{n\ge 1}$ defined in
$(\hat\Omega,
\hat \calf,\Q)$
satisfying the following conditions:
\begin{enumerate}
\item
setting
$$
\hat S_0=0,\qquad\hat \eta_0=a_0,
\qquad
\hat I_t=\sum_{n\ge 0}\hat  \eta_{n}1_{ [\hat S_n,\hat S_{n+1})}(t),
$$
 the process
$\hat I$
satisfies
\begin{equation}\label{distkrylovdelta}
    \tilde\rho (\hat I, \hat \alpha)=
\E^\Q\left[\int_0^T  \rho (\hat I_t,\hat \alpha_t)\,dt\right]<\delta;
\end{equation}
\item denoting
$\hat \mu=\sum_{n\ge1}\delta_{(\hat S_n,\hat \eta_n)}$ the random measure
associated to $(\hat S_n,\hat \eta_n)_{n\ge 1}$,
$\F^{\hat\mu}=(\calf_t^{\hat\mu})_{t\geq 0}$
 the natural filtration of
$\hat\mu$ and $\hat\G\vee \F^{\hat\mu}=(\hat\calg_t\vee\calf_t^{\hat\mu})_{t\geq 0}$,
then
the $\hat\G\vee \F^{\hat\mu}$-compensator of    $\hat \mu$ under $\Q$ is absolutely continuous
 with respect to $\lambda(da)\,dt$ and
it can be written in the form
\begin{equation}\label{compensatoreBmu}
 \hat\nu_t(\hat\omega,a)\, \lambda(da)\,dt
\end{equation}
 for  some  $\calp(\hat\G\vee \F^{\hat\mu})
 \otimes \calb(A)$-measurable  function $\hat\nu$ satisfying
 \begin{equation}\label{compensatoreBmubounds}
 0<
 \inf_{\hat\Omega\times [0,T]\times A}\hat\nu\le
  \sup_{\hat\Omega\times [0,T]\times A}\hat\nu<\infty.
  \end{equation}
\end{enumerate}

\end{Proposition}
\textbf{Proof.}
Fix $\alpha$ and $\delta$ as in the statement of the Proposition.
It can be proved that there exists an $A$-valued process  $\bar\alpha$     such that
$\tilde\rho (\alpha,\bar\alpha)<\delta/3$ and $\bar\alpha$
has the form $\bar\alpha_t=\sum_{n= 0}^{N-1} \alpha_{n}1_{ [t_n,t_{n+1})}(t)$,
where  $0=t_0<t_1<\ldots t_N=T$ is a deterministic subdivision of $[0,T]$,
$\alpha_0,\ldots,\alpha_{N-1}$ are $A$-valued
random variables that take only a finite number of values, and each $\alpha_n$ is $\calg_{t_n}$-measurable:
this follows from
 Lemma 3.2.6 in \cite{80Krylov} (which was stated above as  
 Lemma \ref{densitykrylov}) where  it is proved that the set of admissible
controls $\bar\alpha$ having the form specified in the lemma are dense in
 the set of all  $\G$-progressive $A$-valued processes with respect to the metric $\tilde \rho$.

We can  (and will)  choose $\bar\alpha$ satisfying $\alpha_0=a_0$
($a_0$ is the same as in {\bf (A2)}). Indeed
this additional requirement can be fulfilled by adding, if necessary, another point $t'$
close to $0$ to the subdivision
and modifying $\bar\alpha$ setting $\bar\alpha_t=a_0$  for   $t\in [0,t')$.
This modification is as close as we wish to the original process with respect
to the metric $\tilde \rho$, provided $t'$ is chosen sufficiently small.

Finally, we further extend   $\bar\alpha$ to a function
defined on $\Omega\times [0,\infty)$ in a trivial way setting
$\bar\alpha_t=\sum_{n= 0}^{\infty} \alpha_{n}1_{ [t_n,t_{n+1})}(t)$ where
$\alpha_n =  \alpha_{N-1}$ for $n\ge N$ and $t_n=t+n-N$ for
$n>N$.
This way $\bar\alpha$ is associated to the marked point process $(t_n,\alpha_n)_{n\ge 1}$ and $\bar \alpha_0=a_0$.

\vspace{1mm}

Next recall  the spaces
 $(\Omega',\calf',\P')$ and
 $(\hat\Omega,\hat\calf,\Q)$ and the filtration $\hat\G$
 introduced before the statement
 of Proposition
\ref{extensionapproximation}. Fron now on the symbol $\E$ denotes the expectation with respect to $\Q$.
We  extend the processes $\alpha$ and $\bar\alpha$ to $\hat\Omega\times
[0,\infty)$ and denote $\hat\alpha$ and $\hat{\bar{\alpha}}$ the corresponding
extensions. We note that clearly
\begin{eqnarray} \label{rhotildeuno}
\tilde \rho(\hat\alpha,\hat{\bar{\alpha}})
 &=& \tilde \rho(\alpha,{\bar{\alpha}})<\delta/3.
\end{eqnarray} 

The next step of the proof consists in constructing a sequence of random
measures $\kappa^m$ whose associated piecewise constant
trajectories, denoted $\hat\alpha_t^m$,  approximate $\hat   \alpha$
in the sense of the metric $\tilde\rho$.
The construction will be carried out in such a way that $\kappa^m$
admits a compensator absolutely continuous with respect to the measure
$\lambda(da)\,dt$.

For every   $m\ge 1$, let $\bfB(b,1/m)$ denote the open ball  of radius $1/m$, with respect to the metric $\rho$, centered at $b\in A$.
Since $\lambda(da)$ has full support, we have $\lambda(\bfB(b,1/m))>0$ and we can  define a transition kernel $q^m(b,da)$ in $A$ setting
\begin{eqnarray} \label{kernelsuA}
q^m(b,da)&=& \frac{1}{\lambda(\bfB(b,1/m))}\, 1_{\bfB(b,1/m)}(a) \lambda(da).
\end{eqnarray} 
To the kernel $q^m$ we will apply the following result.

\begin{Lemma}\label{transkernel}
Let $q(b,da)$ be a transition kernel on the Borel space $A$ (i.e., $q(\cdot,B):A\to A$ is Borel measurable for every $B\in \calb(A)$ and $q(b,\cdot)$ is a probability on $(A,\calb(A))$ for every $b\in A$).
Then there exists a function $\bar q:A\times (0,1)\to A$, measurable with respect to $\calb(A \times  (0,1))$ and $\calb(A)$, such that for every
$b\in A$,  $q(b,\cdot)$ is the image of the Lebesgue measure  on $(0,1)$ under the mapping
$u \mapsto \bar q(b,u)$; equivalently,
$$
\int_A k(a)\,q(b,da) \; = \; \int_0^1k( \bar q(b,u))\,du,
$$
for every nonnegative measurable function $k$ on $A$.
\end{Lemma}
\noindent {\bf Proof.} 
The existence of the function $\bar q$
 is well known when $A$ is a
 separable complete metric space: more precisely, for fixed $b$, the construction of a measurable function $\bar q(b,\cdot): (0,1)\to A$ such that
 the image of the Lebesgue measure  coincides with the given measure $q(b,\cdot)$ is performed as part of the proof of Skorohod's Representation Theorem and an inspection  of the proof shows that the resulting function $\bar q(b,u)$ is Borel measurable as a function of the pair $(b,u)$; the reader may consult 
for instance \cite{Zabczyk96} Theorem 3.1.1 for details. The case when $A$ is a general Borel space follows immediately, because it is known that  any Borel  space is either finite or countable (with the discrete topology)
or isomorphic, as a measurable space, to the real line  $\R$: see e.g.
\cite{BertsekasShreve78}, Corollary 7.16.1.
\qed

\vspace{2mm}

According to this Lemma, there exists a Borel measurable function $\bar q^m:A\times (0,1)\to A$
such that for every $b\in A$ the measure $ q^m(b,\cdot)$ is the image of the Lebesgue measure  on $(0,1)$ under the mapping
$u\mapsto \bar q^m(b,u)$
Thus, if  $U$ is a random variable defined on some probability space and having uniform law on $(0,1)$ then, for fixed $b\in A$, the $A$-valued random variable
$\bar q^m(b,U)$ has law $q^m(b,da)$. From now on write $q^m(b,u)$ instead of $\bar q^m(b,u)$, for simplicity.

\vspace{1mm}

For fixed   $m\ge 1$,
define    $V^m_0=R^m_0=S^m_0=0$ and
 $$
R_n^m=t_n+V^m_1+\ldots +V^m_n,
\quad
S_n^m=t_n+V^m_1+\ldots +V^m_{n-1},
\quad \beta_n^m=q^m(\alpha_n,U^m_n),
\qquad \quad  n\ge 1.
$$

Since we assume $t_n<t_{n+1}$ and since $V^m_n>0$
we see that $(R^m_n,\beta^m_n)_{n\ge1}$
is a marked point process in $A$. Also note that
$R^m_{n-1}<S^m_n<R^m_n$ for $n\geq 1$.
Let
$$\kappa^m=\sum_{n\ge1}\delta_{(R^m_n,\beta^m_n)},
\qquad
\hat\alpha_t^m \; = \; \sum_{n\ge 0} \beta^m_{n}1_{ [R^m_n,R^m_{n+1})}(t),
$$
(with the convention $\beta^m_0=a_0$)
denote the corresponding random measure and the associated
trajectory. We claim that
\begin{equation}\label{compensatorofmodifiedppmbis}
\tilde\rho (\hat{\bar{\alpha}}, \hat\alpha^m)\to 0
\end{equation}
as $m\to\infty$.
Indeed, since $0=t_0<t_1<\ldots t_N=T$ we have
\begin{equation}\label{decompdist}
\tilde\rho (\hat{\bar{\alpha}}, \hat\alpha^m)=
\sum_{n=0}^{N-1}\E\int_{t_n}^{t_{n+1}} \rho(\hat{\bar{\alpha}}_t, \hat\alpha^m_t)\,dt.
\end{equation}
Note that $t_n<R_n^m$, and whenever $ R_n^m\le t<t_{n+1}<R_{n+1}^m$
we have $\hat{\bar{\alpha}}_t=\alpha_n$, $\hat\alpha^m_t=\beta_n^m$ and so
$\rho(\hat{\bar{\alpha}}_t, \hat\alpha^m_t)=\rho(\alpha_n,\beta^m_n)<1/m$
since, for every $b$ $\in$ $A$ , $q^m(b,da)$ is supported in $\bfB(b,1/m)$.
If $R_n^m<t_{n+1}$
then, recalling that $\rho <1$,
\begin{eqnarray*}
  \int_{t_n}^{t_{n+1}} \rho(\hat{\bar{\alpha}}_t, \hat\alpha^m_t)\,dt
  &=&
\int_{t_n}^{R_n^m}   \rho(\hat{\bar{\alpha}}_t, \hat\alpha^m_t)\,dt
+\int_{R_n^m}^{t_{n+1}} \rho(\hat{\bar{\alpha}}_t, \hat\alpha^m_t)\,dt
\\
&\le&
(R_n^m-t_n) + \frac1m (t_{n+1}- R_n^m)
\\
&\le&
V_1^m+\ldots +V_n^m + \frac1m (t_{n+1}- t_n).
\end{eqnarray*}
If $R_n^m\ge t_{n+1}$
then the same inequality still holds since we even have
$$
\int_{t_n}^{t_{n+1}} \rho(\hat{\bar{\alpha}}_t, \hat\alpha^m_t)\,dt
\le
t_{n+1}- t_n \le R_n^m - t_n = V_1^m+\ldots +V_n^m.
$$
Substituting in \eqref{decompdist} and
computing the expectation of the exponential random variables $V_n^m$
we arrive at
$$
\tilde\rho (\hat{\bar{\alpha}}, \hat\alpha^m)=
\sum_{n=0}^{N-1}\left(
\lambda_{1m}^{-1}+\ldots +\lambda_{nm}^{-1} + \frac1m (t_{n+1}- t_n)\right)
\le
\sum_{n=1}^\infty\lambda_{nm}^{-1} +\frac{T}{m}\le \frac1m +\frac{T}{m}
$$
 which proves the claim \eqref{compensatorofmodifiedppmbis}.
From now on we fix  a value of $m$ so large that
\begin{eqnarray} \label{rhotildedue}
\tilde\rho (\hat{\bar{\alpha}}, \hat\alpha^m)&<&\delta/3.
\end{eqnarray} 

Let
$\F^{\kappa^m}=(\calf^{\kappa^m}_t)$ denote the natural filtration of $\kappa^m$ and
 set
\[
\H^m = (\calh_t^m)_{t\ge 0}=
(\hat\calg_t\vee \calf^{\kappa^m}_t)_{t\ge 0}.
\]
We have the following technical result that describes the compensator
$\tilde\kappa^m$ of $\kappa^m$ with respect to the filtration $\H^m$.

 \begin{Lemma}\label{MPPperturbed}
With the previous assumptions and notations, the compensator of the random measure
$\kappa^m$ with respect to $\H^m$ and $\Q$ is given by the formula
 \begin{eqnarray*}
 \tilde\kappa^m(dt,da) &=& \sum_{n\ge 1}  1_{(S_n^m,R^m_n]}(t)\,
 q^m(\alpha_n,da) \lambda_{nm}.
 \end{eqnarray*}
\end{Lemma}
{\bf Proof of Lemma \ref{MPPperturbed}.} To shorten notation, we drop all the sub- and superscripts
$m$ and write $\tilde\kappa$, $S_n$, $R_n$,
$ q$, $ \lambda_{n}$,
$\F^{\kappa}$,
$\H = (\calh_t)=
(\hat\calg_t\vee \calf^{\kappa}_t)$
instead of
$\tilde\kappa^m$, $S_n^m$, $R^m_n$,
$ q^m$, $ \lambda_{nm}$, $\F^{\kappa^m}$, $\H^m $ $=$ $ (\calh_t^m)$ $=$
$(\hat\calg_t\vee \calf^{\kappa^m}_t)$.

Let us first check that $\tilde\kappa(dt,da)$, defined by the  formula above, is an $\H$-predictable random measure. The variables $R_n$ are clearly
$\F^\kappa$-stopping times and hence
$\H$-stopping times and therefore
 $S_n= R_{n-1}+t_n-t_{n-1}$  are also $\H$-stopping times.
 Since
  $\alpha_n$    are $\calf_{t_n}$-measurable and
  $\calf_{t_n}\subset  \calh_{t_n}\subset \calh_{S_n}$, $\alpha_n$
  are also
$\calh_{S_{n}}$-measurable. It follows that for every $C\in\calb(A)$ the process
$
1_{(  S_{n},R_n]}(t)\, q(\alpha_n,C) \lambda_n
$
is $\H$-predictable and finally that  $\tilde\kappa(dt,da)$ is an $\H$-predictable random measure.

To finish the proof  we need now to verify that for every  positive $\calp(\H)\otimes \calb(A)$-measurable
 random field $H_t(\omega,a)$ we have
 $$
 \E \Big[\int_0^\infty\int_AH_t(a)\,\kappa(dt\,da) \Big]  \;= \;
 \E \Big[ \int_0^\infty\int_AH_t(a)\,\tilde\kappa(dt\,da)\Big].
 $$
 Since $\calh_t=\calf_t\vee \calf^\kappa_t$, by a monotone class argument
 it is enough to consider $H$ of the form
 $$H_t(\omega,a)=H_t^1(\omega)H_t^2(\omega)k(a),
 $$
 where $H^1$ is a positive $\hat\G$-predictable random process,
$H^2$ is a positive $\F^\kappa$-predictable random process
 and $k$ is a positive $\calb(A)$-measurable function.
 Since $\F^\kappa$ is the natural filtration of $\kappa$,
 by a known result (see e.g. \cite{ja} Lemma (3.3) or 
Proposition \ref{propnaturalfiltration}-(iv) above)
 $H^2$ has the following form:
 \begin{eqnarray*}
    H^2_t &=& b_1(t)1_{(0,R_1]}(t)+
 b_2(\beta_1,R_1,t)1_{( R_1,R_2]}(t)
 \\
 &&+
 b_3(\beta_1,\beta_2 ,R_1,R_2,t)1_{(R_2,R_3]}(t)+ \ldots
 \\
     && +
 b_n(\beta_1,\ldots,\beta_{n-1},R_1,\ldots, R_{n-1},t)1_{(R_{n-1},R_n]}(t)+
 \ldots,
 \end{eqnarray*}
where each $b_n$ is a positive measurable deterministic function of
$2n-1$ real variables.
Since
$$
\E \Big[ \int_0^\infty\int_AH_t(a)\,\kappa(dt\,da) \Big] \; = \;
\E \Big[  \sum_{n\ge 1} H_{R_n}(\beta_n) \Big]
 $$
to prove the thesis  it is enough to check that for every $n\ge1$ we have the equality
$$
\E \big[ H_{R_n}(\beta_n) \big]  =
 \E \Big[ \int_0^\infty \int_A H_t(a)\,  q(\alpha_n,da) \lambda_n\,
 1_{  S_{n}<t\le R_n }\,dt \Big]
 $$
which can also be written
$$
\begin{array}{l}\dis
 \E\big[H^1_{R_n}  b_n(\beta_1,\ldots,\beta_{n-1},R_1,\ldots, R_{n-1},R_n)    k(\beta_n) \big]
 =
 \\\dis
   \E \Big[ \int_0^\infty \int_A H_t^1  b_n(\beta_1,\ldots,\beta_{n-1},R_1,\ldots, R_{n-1},t)k(a)
 q(\alpha_n,da)
\lambda_n\,  1_{  S_{n}<t\le R_n }\,dt \Big].
\end{array}
$$
We use the notation
$$
K_n(t) \; = \; H_t^1\,
 b_n(\beta_1,\ldots,\beta_{n-1},R_1,\ldots, R_{n-1},t)
 $$
 to reduce the last equality to
\begin{equation}\label{thesisrewritten}
    \E\, [K_n(R_n)
    k(\beta_n)]
=
 \E\Big[\int_0^\infty \int_A K_n (t)\,k(a)\,
 q(\alpha_n,da) \lambda_n\,
 1_{  S_{n}<t\le R_n }\,dt\Big].
\end{equation}
By the definition of $R_n$ and $\beta_n$, we have
$\E[K_n(R_n)k(\beta_n)]$ $=$  $\E[K_n(S_n+V_n)   k(q(\alpha_n,U_n))] $.
As noted above,   since $U_n$ has uniform law on $(0,1)$,  the random variable
$q(b,U_n)$ has law $q(b,da)$ on $A$,   for any fixed $b\in A$.
We note that $R_1,\ldots, R_{n-1}, S_n$ are measurable with respect to
$\sigma (V_1,\ldots, V_{n-1})$, that $\beta_1,\ldots, \beta_{n-1}$ are
measurable with respect to
$\hat\calg_\infty\vee\sigma (U_1,\ldots, U_{n-1})$ and therefore that
 the random elements $U_n $,
 $S_n$  and $(\hat\calg_\infty, \beta_1,\ldots,\beta_{n-1},R_1,\ldots, R_{n-1},S_n)$ are all  independent.
 Recalling that $V_n$ is exponentially distributed with parameter $\lambda_n$
 we obtain
\begin{equation}\label{S_nindep}
        \E\, [K_n(R_n)k(\beta_n)]
    \; = \;
    \E \Big[ \int_0^\infty\int_A\, K_n(S_n+s)\,    k(a)\,q(\alpha_n,da)\,\lambda_n
    e^{-\lambda_ns}\,ds\Big].
\end{equation}
 Using  again the independence of $V_n$  and
$(\hat\calg_\infty, \beta_1,\ldots,\beta_{n-1},R_1,\ldots, R_{n-1},S_n)$ we also have
\begin{eqnarray*}
& &   \E \Big[ \int_0^\infty\int_A\, K_n(S_n+s)\,    k(a)\,q(\alpha_n,da)\,
    \lambda_n\,1_{V_n\ge s}\, ds \Big] \\
 &=&     \E \Big[ \int_0^\infty\int_A\, K_n(S_n+s)\,    k(a)\,q(\alpha_n,da)\,  \lambda_n\,
     \P(V_n\ge s)\, ds\Big]
 \end{eqnarray*}
and since  $\P(V_n\ge s)=e^{-\lambda_ns}$, this coincides with the
right-hand side of \eqref{S_nindep}. By a change of variable we arrive at
equality  \eqref{thesisrewritten}:
\[
\begin{array}{lll}\dis
        \E\, [K_n(R_n)k(\beta_n)]
    & = &\dis
    \E \Big[ \int_{S_n}^\infty\int_A\, K_n(t)\,    k(a)\,q(\alpha_n,da)\,
\lambda_n\,    1_{V_n\ge t-S_n}\, dt \Big]
    \\
& =&\dis
    \E \Big[ \int_{0}^\infty\int_A\, K_n(t)\,    k(a)\,q(\alpha_n,da)\,
  \lambda_n\,   1_{S_n<t\le R_n}\,dt\Big].
    \end{array}
\]
This concludes the proof of Lemma \ref{MPPperturbed}.    \qed

\vspace{3mm}

It follows from this lemma that
the $\H^m$-compensator of   $\kappa^m$ under $\Q$
is absolutely continuous with respect to $\lambda(da)\,dt$ and it can be written in the form
  $$
 \tilde\kappa^m(dt,da) \; = \;  \phi^m_t(a)\, \lambda(da)\,dt
 $$
for a suitable nonnegative $\calp(\H^m)\otimes \calb(A)$-measurable  function $\phi^m$
which is bounded on $\hat\Omega\times [0,T]\times A$. Indeed,
from the choice of the kernel $q^m(b,da)$ we obtain
\begin{eqnarray*}
\phi^m_t(a) &=& \sum_{n\ge 1}  1_{(S^m_{n},R^m_n]}(t)\,
\frac{1}{\lambda(\bfB(\alpha_n,1/m))}\, 1_{\bfB(\alpha_n,1/m)}(a) \lambda_{nm}.
\end{eqnarray*}
which is bounded on $\hat\Omega\times [0,T]\times A$
since each $\alpha_n$ takes only a finite number of values and
$S^m_{N}> t_N=T$, so that the values of $\phi_t^m$ on $[0,T]$
only depend on the first $N-1$ summands.

In the final step of the proof we will modify the random
measure $\kappa^m$ by adding an independent Poisson process $\pi^k$
with ``small'' intensity. This will not affect too much the
$\tilde\rho$-distance between the corresponding trajectories and will
produce a random measure whose compensator
remains absolutely continuous with respect to the measure
$\lambda(da)\,dt$ and has a bounded density which,  in addition, is bounded
away from zero.

Recall that on the space
$(\Omega',\calf',\P')$ we assumed that for every integer $k\ge 1$ there
existed    a Poisson random measure $\pi^k$  on $(0,\infty)\times A$, admitting compensator $k^{-1}\lambda(da)\,dt$ with respect to its natural filtration.
We will consider $\pi^k$  as defined in $(\hat\Omega,\hat\calf)$.
Each
$\pi^k$ has the form
\begin{eqnarray*}
\pi^k &=& \sum_{n\ge 1}\delta_{(T_n^k , {\xi_n^k })},
\end{eqnarray*}
for a marked point process $(T_n^k , {\xi_n^k })_{n\ge1}$
on $(0,\infty)\times A$,
and we denote  $\F^{\pi^k}$ $=$ $(\calf^{\pi^k}_t)$ its natural filtration.
Let us define another random measure setting
$$
\mu^{km} =\kappa^m+\pi^k.
$$
Note that the jumps times $(R_n^m)_{n\ge1}$  are independent
of the jump times $(T_n^k )_{n\ge1}$, and the latter have
absolutely continuous laws. It follows that, except possibly
on a set of $\Q$ probability zero, their graphs  are
disjoint, i.e. $\kappa^m$ and $\pi^k $ have no common jumps. Therefore,
  the random measure $\mu^{km} $
and its associated pure jump process (denoted $ I^{km}$)
    admit a representation
$$
\mu^{km} =\sum_{n\ge 1}\delta_{(S_n^{km} , {\eta^{km} _n})},
\qquad
 I^{km} _t =
 \sum_{n\ge 0} \eta^{km} _{n}1_{ [S_n^{km} ,S^{km} _{n+1})}(t), \;\;\; t \geq 0,
$$
where $\eta^{km} _{0}=a_0$,
$(S_n^{km} ,\eta^{km} _n)_{n\ge 1}$ is a marked point process,  each $S_n^{km} $ coincides with one of the times $R_n^m$ or one of the times
$T_n^k $, and each $\eta_n^{km} $ coincides with one of the random variables $\xi_n^k $ or one of the random variables $\beta_{n}^m$.
We claim that, for large $k$, $I^{km}$ is close to $\hat\alpha^n$ with respect to the metric $\tilde \rho$,
namely that
\begin{equation}\label{approxvalue}
\tilde \rho(I^{km},\hat\alpha^m)\to 0
\end{equation}
as  $k\to\infty$.
To prove this claim it suffices
to prove that $I^{km}\to \hat\alpha^m$ in $dt\otimes d\Q$-measure.
Recall that the jump times of $\pi^k $ are denoted $T_n^k $.
Since  $T_1^k $ has exponential law with parameter
$ \lambda(A)/k$ the event $B_k =\{T_1^k  >T\}$ has probability
$e^{-\lambda(A)T/k}$, so that $\Q(B_k)\to 1$ as $k\to \infty$.
Noting that  on the set $B_k$, we have $\hat\alpha^m_t =
\alpha_0=a_0=\eta_0^{km}=
I^{km} _t$ for all $t\in [0,T]$,
 the claim
 \eqref{approxvalue} follows immediately. We will fix from now on an integer
$k$ so large that
\begin{eqnarray} \label{rhotildetre}
\tilde\rho (\hat\alpha^m,I^{km})<\delta/3.
\end{eqnarray} 

Having fixed both $m$ and $k$ we now define, for $n\ge 0$,
$$
\hat S_n= S_n^{km},\quad
\hat \eta_n =\hat \eta_n^{km},
\quad
\hat \mu=\sum_{n\ge1}\delta_{(\hat S_n,\hat \eta_n)},
\quad
\hat I_t=\sum_{n\ge 0}\hat  \eta_{n}1_{ [\hat S_n,\hat S_{n+1})}(t),
$$
so that the random measure $\hat\mu$ and the associated process
$\hat I$ coincide with
$\mu^{km}$ and
$I^{km}$ respectively.
The inequalities \eqref{rhotildeuno}, \eqref{rhotildedue}, \eqref{rhotildetre}
imply that $\tilde\rho (\hat\alpha,\hat I)<\delta$, which gives
\eqref{distkrylovdelta}.

To finish the proof
it remains to prove \eqref{compensatoreBmu}-\eqref{compensatoreBmubounds}.
We first note that,
since $\kappa^m$ and $\pi^k $ are independent, it is easy to prove
 that
 $\hat\mu=\mu^{km} $ has compensator  $(\phi^m_t(a)+k^{-1})\,\lambda(da)\,dt$ with respect
 to the filtration $\H^m\vee\F^{\pi^k}$ $:=$ $(\calh^m_t\vee\calf^{\pi^k }_t)_{t\ge0}
 =(\hat\calg_t\vee \calf^{\kappa^m}_t\vee \calf^{\pi^k }_t)_{t\ge0}$.
Let $\F^{\hat\mu}=(\calf_t^{\hat\mu})_{t\geq 0}$ denote
 the natural filtration of
$\hat\mu$ and let $\hat\G\vee \F^{\hat\mu}$ be the filtration $(\hat\calg_t\vee\calf_t^{\hat\mu})_{t\geq 0}$, which is smaller
than $\H^m\vee\F^{\pi^k}$.
We wish to compute the compensator of $\hat\mu$ with respect to $\hat\G\vee \F^{\hat\mu}$
under $\Q$.
To this end, consider the measure space
$([0,\infty)\times\Omega\times A, \calb([0,\infty))\otimes\calf\otimes \calb(A),dt\otimes\Q(d\omega)\otimes\lambda(da))$.
Although this is not a probability space, one can define in a standard way the conditional expectation
 of any positive measurable
function, given an arbitrary sub-$\sigma$-algebra.
Let us denote by  $\hat\nu_t(\hat\omega,a)$ the conditional expectation
of the  random field $ \phi_t^m(\hat \omega,a)+k^{-1}$
 with respect to the $\sigma$-algebra $\calp(\hat\G\vee \F^{\hat\mu})
 \otimes \calb(A)$. It is then easy to verify that
the compensator of $\hat\mu $ with respect to
$\hat\G\vee \F^{\hat\mu}$ coincides with $\hat\nu$.
Moreover, since  $ \phi_t^m(\hat \omega,a)$ is nonnegative and bounded on
$\hat\Omega\times [0,T]\times A$,
we can take a version of $\hat\nu$ satisfying
$$
 k^{-1}\le \inf_{\hat\Omega\times [0,T]\times A}\hat\nu\le
  \sup_{\hat\Omega\times [0,T]\times A}\hat\nu<\infty.
  $$
Now \eqref{compensatoreBmu}-\eqref{compensatoreBmubounds}
are proved and the proof of Proposition
\ref{extensionapproximation}
is finished.
\qed

\section*{Acknowledgments}
I wish to end this paper with my warmest thanks to all the friends and colleagues that have worked with me on this subject and have given an invaluable contribution to its development: Elena Bandini, Fulvia Confortola, Andrea Cosso,  Idris Kharroubi, Marie-Amélie Morlais,
Huyên Pham, Federica Zeni. 

A special thank goes to Jiongmin Yong, Editor in chief of {\it Numerical Algebra, Control and Optimization},
for inviting me to prepare this review.

\vspace{5mm}

\small
\bibliographystyle{plain}
\bibliography{biblio}

\end{document}